\documentclass[a4paper,twoside,10pt]{article}

\usepackage{amssymb}
\usepackage{amsmath}
\usepackage{epsfig}

\newtheorem{The}{Theorem}
\newtheorem{Lem}{Lemma}
\newtheorem{Cor}[Lem]{Corollary}
\newtheorem{Pro}[Lem]{Proposition}

\newtheorem{Rem}[Lem]{Remark}
\newtheorem{Defi}{Definition}
\newtheorem{Exem}{Example}

\def\a{\alpha}

\def\C{\mathbf C}
\def\D{\Delta}
\def\d{\delta}

\def\f{\phi}
\def\g{\gamma}

\def\l{\lambda}

\def\o{ o}
\def\O{\mathbf O}
\def\p{\pi}

\def\Q{\mathbf Q}
\def\R{\mathbf R}
\def\r{\rho}

\def\s{\sigma}
\def\t{\tau}

\def\w{\omega}
\def\Z{\mathbf Z}
\def\z{\zeta}

\title{Toric embedded resolutions \\ of quasi-ordinary hypersurface singularities}
\author{Pedro D. Gonz\'alez P\'erez}
\date{ }
\begin{document}

\maketitle

\begin{quote}
{\bf Abstract.}\footnotetext[1]{ Math.\ Subject Class.: 32S15, 32S45, 14M25, 14E15, 32S25}
We build two toric embedded resolutions procedures 
of a reduced  quasi-ordinary hypersurface
singularity $(S,0)$ of dimension $d$.
The first one
provides an embedded resolution as hypersurface of $(\C^{d+1}, 0)$
as a composition  of toric morphisms 
which depend only on the 
{\it characteristic monomials} associated to a quasi-ordinary projection 
$(S,0) \rightarrow (\C^{d}, 0)$.
This gives a positive answer to a  a question of  Lipman (see \cite{Lipman3}).
The second method applies to the analytically irreducible case; if
 $g \geq 1 $ denotes the number of characteristic monomials
we re-embed the germ $(S, 0)$  in the affine space $(\C^{d+g}, 0)$
by using certain approximate roots of a suitable Weierstrass 
polynomial defining the embedding $(S,0) \subset  (\C^{d+1}, 0)$.
We build a toric morphism 
which is a  simultaneous toric embedded resolution 
of the irreducible germ  $(S, 0) \subset  (\C^{d+g}, 0)$, 
and of an affine toric variety  $Z^{\Gamma}$
obtained from $(S, 0) \subset  (\C^{d+g}, 0)$ by specialization
and defined by a rank $d$ semigroup ${\Gamma}$ 
generalizing the classical semigroup of a plane branch.
Finally we compare both resolutions 
and we prove that
the first one is the restriction of the second to a 
smooth $(d+1)$-variety containing the strict transform of $S$.

\end{quote}


\section*{Introduction}

A germ  of  complex analytic variety is {\it quasi-ordinary} if there exists
a finite projection, called quasi-ordinary,
to the complex affine space  $({\mathbf C}^d, 0)$
with discriminant locus contained in  a normal crossing divisor
(for instance, the singularities of complex analytic curves are quasi-ordinary).
These singularities appear classically in Jung' strategy to obtain
the resolution of singularities of  surfaces from the embedded 
resolution of plane curves (see \cite{Jung}, 
 \cite{Walker} and \cite{LipmanR}).
Some properties of complex analytic curve singularities generalize to
quasi-ordinary hypersurface singularities:
for instance, Jung-Abhyankar's theorem guarantees the existence
of fractional power series parametrizations generalizing the classical Newton-Puiseux 
parametrizations of the plane curve case; by comparing 
these parametrizations  we obtain a finite set of 
{\it distinguished} or {\it characteristic monomials}
which generalize 
the notion of 
characteristic exponents in the plane branch case.

The results on quasi-ordinary hypersurface singularities concern 
mainly the analytically irreducible case: 
Lipman builds a non embedded resolution procedure of a quasi-ordinary surface 
where only quasi-ordinary singularities occurs
and uses it to prove the
the analytical invariance properties of the characteristic monomials
(see \cite{Lipman0}, and  \cite{Lipman1}); 
another proof of this result was given by  Luengo in \cite{Luengo}; 
more generally Gau 
proved that
the characteristic monomials, suitably normalized by an inversion formulae of Lipman \cite{Lipman0},  
define a 
complete invariant of the {\it embedded topological type} of the 
quasi-ordinary hypersurface singularity
(see  \cite{Gau});
Gau's proof involves Lipman's description of the divisor class group of the 
singularity in terms of the characteristic monomials
(see  \cite{Lipman2}).

An important step to establish the relations between the topological type and the 
embedded resolutions
of a hypersurface singularity, which are well-known in the case of
 plane curve singularities 
 (see \cite{Zequi}, \cite{Equi} and \cite{Reeve}), 
is to determine if
the characteristic monomials 
of a
hypersurface quasi-ordinary singularity determine a 
procedure of embedded resolution. 
This is the content of Lipman's  open problem 5.1
 (see \cite{Lipman3}) which is stated in the context of 
the generalizations of {\it equisingularity}, in particular by using Zariski's 
work on the {\it dimensionality type}
 with respect to the classification by  {\it equiresolution}.
In the case of an
analytically irreducible quasi-ordinary surface germ 
Ban and McEwan (see \cite{Ban}) 
have found a such a procedure 
following 
the  algorithm  of resolution of Bierstone and Milman, developed from
the work of Hironaka.
Villamayor has
given a solution to Lipman's problem for any 
quasi-ordinary hypersurface singularity  (see \cite{Villamayor}).
Villamayor's approach studies the abelian branch covering of 
the affine space obtained by taking suitable roots of the 
regular parameters defining the components of the discriminant locus.
By the Jung-Abhankar's Theorem the equation of the quasi-ordinary hypersurface under this extension splits in a product of Weierstrass polynomials of degree one. 
The singularity obtained is a non transversal intersection of smooth 
hypersurfaces, whose embedded resolution 
requires the simplest combinatorial part of 
Hironaka's method. The important point that he proves is that this resolution  
procedure is Galois equivariant, in such a way that when 
taking the quotients by 
the Galois action the local constructions glue up together 
defining a modification of the embedded quasi-ordinary hypersurface.
The ambient space obtained in this way has only {\it toric quotient singularities}
and a canonical resolution of these singularities 
(see \cite{Villamayor2}) provides an embedded 
desingularization of the quasi-ordinary hypersurface.
The desingularization  obtained is not 
necessarily an isomorphism 
outside the singular locus of the quasi-ordinary hypersurface.

In this paper we give another solution to Lipman's problem in two different ways.

In the first one
we  build an  embedded resolution of a reduced quasi-ordinary hypersurface germ
$(S, 0) \subset (\C^{d+1}, 0)$ as a
composition of toric morphisms which depend only on 
the characteristic monomials (see Theorem \ref{hres}).
The first toric morphism we build is  defined by the 
{\it dual Newton diagram} of a suitable Weierstrass polynomial $f \in \C \{ X_1,\dots, X_d \} [Y]$ 
defining the embedding $(S,0) \subset (\C^{d+1}, 0)$.
Suitable here means that $Y$ is a good coordinate: the Newton polyhedron of $f$, 
whose compact faces have dimension at most one, 
is canonically determined by the characteristic monomials. 
We study the strict transform $S'$ of $S$ by this modification: we show that 
the restriction $\p_1 : S' \rightarrow S$ is a finite map.
The germ of $S'$  at any of the finitely many points of the fiber $ \p_1^{-1} (0) $  
is a {\it toric quasi-ordinary singularity}, defined as finite branched covering of
a normal affine toric variety unramified over its torus (see \cite{GP}).
From this it follows that it is more natural to
build the resolutions 
for {\it toric} quasi-ordinary hypersurfaces by generalizing to this case 
the notions of characteristic monomials and many of their properties
in the {\it classical} quasi-ordinary case.
At any point of $ \p_1^{-1} (0) $,
the strict transform $S'$  
has less characteristic monomials, with respect to a projection canonically determined from the fixed quasi-ordinary projection of $S$, and 
 we determine them from the given characteristic monomials of $(S, 0)$. 
By iterating we obtain, in a canonical manner from the fixed quasi-ordinary projection of $S$, a
{\it partial embedded resolution}: a normal variety of dimension $d+1$ 
with only toric singularities (not necessarily quotient singularities)
and 
a modification $\p = \p_1 \circ \dots \circ \p_k$ such that 
the strict transform of $S$ is a $d$-dimensional section
transversal to the exceptional fiber $\p^{-1} (0)$  (which is of dimension one). 
This implies that any toric resolution of the 
ambient space is an embedded resolution of the strict transform 
and provides a fortiori an embedded resolution of $S$.
It follows also that the restriction of $\p$ to the 
the strict transform of $S$ is the {\it normalization map} of $S$. This implies that 
the restriction of any of these embedded resolutions to the strict 
transform of $S$ is an isomorphism outside the singular locus of $S$. 
In the case of a plane curve germ 
we show
that our procedure, with respect to a transversal projection, 
leads to the minimal resolution of the curve
and we compare our method with those  given by L\^e, Oka and A'Campo  
(see \cite{Le}, \cite{Oka}, and \cite{A'C}).

The second method builds embedded resolutions of 
an analytically irreducible 
quasi-ordinary hypersurface germ $(S, 0)$ by
generalizing the method of Goldin and Teissier 
for plane branches (see \cite{Rebeca}).
The approach and results of this part are also inspired those obtained by Lejeune and Reguera   
in the case of sandwiched surface singularities (see \cite{LR}) and 
sketched  for plane branches in \cite{LR2}.
If $g \geq 1 $ denotes the number of characteristic exponents 
we re-embed the germ $(S, 0)$  in the affine space $(\C^{d+g}, 0)$, 
by using certain {\it approximate roots} of a suitable Weierstrass 
polynomial defining the embedding $(S,0) \subset  (\C^{d+1}, 0)$.
These approximate roots have {\it maximal contact} in the sense that, 
at each step of the partial resolution there is one approximate root 
whose strict transform   
defines a good coordinate for the strict transform of $S$.
We define a 
toric modification  $p: Z \rightarrow  \C^{d+g}$
 depending only on  a rank $d$  semigroup $\Gamma$, 
which is a  {\it partial  embedded resolution} 
of the irreducible germ  $(S, 0) \subset  (\C^{d+g}, 0)$, 
and of an affine toric variety  $Z^{\Gamma} \subset \C^{d+g}$
obtained from $(S, 0) \subset  (\C^{d+g}, 0)$ by specialization
and defined by the semigroup  ${\Gamma}$
(see Theorem \ref{Bernard}).
This semigroup, which generalizes the classical semigroup of a plane branch, 
does not depend on the quasi-ordinary projection and defines 
a complete invariant of the embedded topological type of $S$, as characterized by Gau 
(see \cite{Tesis} or \cite{semi}).
As in the first  method any toric resolution of singularities of the ambient space $Z$
provides an embedded resolution of $S$. 

We compare the partial resolutions $\p$ and $p$: we prove in Theorem \ref{mismo} that
 $\p$  is the restriction of  $p$   to a 
$(d+1)$-dimensional smooth variety of $Z$ containing the strict transform of $S$.

One of the technical tools common to both methods is  
the construction of toric embedded resolutions of 
non necessarily normal  affine toric varieties equivariantly embedded,
a result obtained in collaboration with Teissier
(see Proposition \ref{binomial}, 
Proposition 6.4 of \cite{Teissier}, and  \cite{cras}).

One important contribution of our approach is a better understanding of the 
structure of the exceptional divisor 
of these resolutions.
The ambient space of the partial resolution $\p$, which is canonical and factors any of these embeded resolutions, 
is built with a toroidal embedding structure such that the associated conic polyhedral complex
with integral structure  (see \cite{TE}) is built explicitly from the characteristic monomials.
This description allows us to re-embed this complex as a fan in an affine space of 
bigger dimension, a technical lemma which is essential to compare 
the partial resolutions $p$ and $\p$ (see Propositions \ref{combo} and \ref{hedra}).
The toric resolutions of the ambient space are defined  by certain regular subdivisions of this fan
(which always exists, see \cite{Cox} and \cite{TE}). These  
regular subdivisions determine many features of the geometry of the exceptional divisor
which are  very useful for the applications:

\begin{itemize} 
\item
In collaboration with N\'emethi and McEwan we have shown that 
the zeta function of the geometric monodromy of the germ $(S, 0)$ 
coincides with the zeta function of the plane curve germ obtained from 
$(S,0)$ 
by intersection with $d-1$ coordinate hyperplanes, 
which are determined by the quasi-ordinary projection
(see \cite{zeta1} and  \cite{zeta2}).
\item
In collaboration with 
Garc\' \i a Barroso
we analyse in \cite{paquetes} the strict transform of the polar hypersurfaces of $(S, 0)$ 
under the partial resolution of $(S, 0)$ and we obtain a decomposition theorem which provides in 
the case of plane curve germs a simple algebraic proof of a Theorem of 
L\^e, Michel and Weber (\cite{LMW}). 
\end{itemize}


The proofs are written in the analytic case. They provide also two 
embedded resolutions of 
quasi-ordinary hypersurface singularities in the algebroid case
(over an algebraically closed field of zero characteristic).

\noindent
{\bf Acknowledgements.} I am grateful to B. Teissier, M. Lejeune-Jalabert,  A. N\'emethi 
for their suggestions and to {\em Universidad de La Laguna} and 
{\em Institut de Math\'ematiques de 
Jussieu} 
for their hospitality. 
The author has been supported by a grant of {\it DGUI del Gobierno de Canarias} and by 
a {\it Marie Curie Fellowship} of the European Community 
Programm ``Improving Human Research Potential and the Socio-economic Knowlegde Base'' 
under contract number HPMF-CT-2000-00877.

\section{Toric maps, Newton polyhedra  and  partial resolution of singularities}   \label{TE}

We introduce the notations and basic definitions of toric geometry
and 
we build 
embedded resolutions of non necessarily normal affine toric varieties.

\subsection{A reminder  of toric geometry} \label{toricg}

We give some definitions and notations (see \cite{Fulton}, \cite{Ewald} and \cite{Oda} for proofs).
If $N \cong \Z^{d+1}$ is a lattice we denote by 
$N_\R $ the real vector space $N \otimes_\Z \R$ spanned by $N$ and 
by $M$ the dual lattice.
A {\it rational convex polyhedral cone} $\s$  in  $N_\R$ is the set non
negative linear combinations of vectors $a^1, \dots, a^s  \in N$.
In what follows a {\it cone} will mean a  rational convex polyhedral cone.
The cone $\s$ is {\it strictly convex} if 
$\s$ contains no linear subspace of dimension $>0$;
the
cone $\s$ is {\it regular\index{regular cone}} if 
the {\it primitive integral vectors} defining the $1$-dimensional faces 
belong to a basis of the lattice  $N$. 
We denote by $\stackrel{\circ}{\s}$ the { \it relative interior} of a cone $\s$.
The {\it dual} cone  $\s^\vee$ (resp. {\it orthogonal} cone $\s^\bot$) of $\s$ is the set
$ \{ w  \in M_\R/ \langle w, u \rangle \geq 0,$  (resp. $ \langle w, u \rangle = 0$)  
$ \; \forall u \in \s \}$. 
A {\it fan\index{fan}} $\Sigma$ is a family of {\it strictly convex 
 cones\index{convex rational polyhedral cone}} in $N_\R$ 
such that any face of such a cone is in the family and 
the intersection of any two of them is a face of each.
The {\it support\index{support of a fan}} of the fan $\Sigma$ is the set 
$ \bigcup_{\s \in \Sigma} \s \subset  N_\R$.
The $i$-skeleton $\Sigma^{(i)}$ is the subset of 
$i$-dimensional 
cones of $\Sigma$. 
 The fan $\Sigma$ is {\it regular} if all its  cones are regular. 


 Any non necessarily normal affine toric
 variety  over the field  $\C$ of complex numbers is of the form
$Z^{\Lambda} = \makebox{Spec} \, \C [ \Lambda ]$ where $\Lambda $ is a {\it monoid}, i.e., 
a sub-semigroup of finite type of a lattice $- \Lambda + \Lambda$ which generates it as a group.
The closed points of $Z^\Lambda$ correspond to homomorphisms of semigroups 
$\Lambda \rightarrow \C$ where $\C$ is considered as a semigroup with
respect to multiplication. The torus embedded in $Z^\Lambda$ is the group
of homomorphisms of semigroups $\Lambda \rightarrow \C- \{ 0 \}$ and acts
naturally on the closed points of  $Z^\Lambda$.
The {\it normalization} of $Z^{\Lambda}$  is obtained from  the
inclusion $\Lambda  \rightarrow  \R_{\geq 0} \Lambda \cap (- \Lambda + \Lambda) $ where $ \R_{\geq 0} \Lambda $
is the 
cone spanned by the elements of  $\Lambda$ (see \cite{TE}). 
The action of the torus has a fixed point
if and only if the
cone  $ \R_{\geq 0} \Lambda $  is strictly convex,
then this point
is defined by the ideal 
$(X^u / u \in \Lambda - \{ 0 \})$ of $\C [ \Lambda ]$
and coincides with the $0$-dimensional orbit; the analytic
algebra $\C \{ \Lambda \}$ of $Z^\Lambda$ at this point
can be viewed as
a subring of the ring $\C [[ \Lambda ]]$ of formal complex power series with exponents 
in the semigroup $\Lambda$ (see \cite{GP} lemme 1.1).

In particular, if  $\s$ is a cone in the fan $\Sigma$ 
the semigroup $\s^\vee  \cap M$  is of finite type, 
it spans the lattice $M$ 
and the variety $Z^{\s^\vee \cap M}$, which 
we denote also by  $Z_{\s, N}$ or by $Z_\s$ when the lattice is clear from the context, 
is {\it normal}.

If $\s \subset \s'$ are cones in the fan $\Sigma$ then
we have an open immersion $Z_\s \subset Z_{\s'} $;  
the affine varieties $Z_\s$  corresponding to cones in a fan
$\Sigma$ 
glue up to define the {\it toric variety\index{toric variety}} 
$ Z_\Sigma$. The torus,   $(\C^*)^{d+1} $,   
is embedded as an open dense subset $Z_{ \{ 0 \} }$ of  
$ Z_\Sigma$, which acts  on each chart $Z_\s$; these actions
paste to an action on $ Z_\Sigma$ that extends
the  product on the torus.
General toric varieties are defined by this property,
the toric varieties which can be defined using fans are precisely the
normal ones (see \cite{TE}).
The toric variety $Z_\Sigma$ is non singular if and only
if the fan $\Sigma$ is regular.

We describe the orbits of the action of the torus on the variety $Z_\Sigma$.
The orbit $\O_{\s, N}$ (which we denote also  by $\O_{\s}$) 
is the Zariski closed subset of   $Z_\s$  defined by the ideal 
$(X^w/ w  \in (\s^\vee - \s^\bot) \cap M
)$ of $\C[\s^\vee  \cap M ]$.
This orbit is a torus for $0 \leq \dim \s <\makebox{rk} \, N$,
 since the associated coordinate ring is the $\C$-algebra of the sub-lattice
$M(\s) :=  M \cap \s^\bot$ of $M$ of codimension  equal to $\dim \s$. 
On the closed orbit $\O_\s$ we consider the {\it special point} $\o_\s $
defined by $X^u (\o_\s) = 1 $ for all $u \in M(\s) $. 
If $\dim \s = \makebox{rk} \, N$ the orbit $\O_\s$ is reduced to the special point.
If $\dim \s < \, \makebox{rk} \, N$ we have an exact sequence of lattices:
\[
0 \rightarrow M(\s) \rightarrow   M  \stackrel{j}{\rightarrow}
M_\s \rightarrow 0.
\]
If 
$ 0 \rightarrow N_\s  \stackrel{j^*}{\rightarrow} N  \rightarrow
N(\s) \rightarrow 0$
 is the dual exact sequence the lattice  $N_\s $ spanned by $ \s \cap N $ is of dimension
equal to $\dim \s$  and the semigroup $ \s^\vee_{N_\s} $ associated to the cone
$\s$ with respect to the lattice $N_\s$ is isomorphic to  $ j(\s^\vee  \cap M)$. 
If we choose a splitting $M \cong M(\s) \oplus  M_\s$ we obtain 
a semigroup isomorphism $\s^\vee  \cap M  \cong M(\s) \oplus
 (\s^\vee_{N_\s}   \cap M_\s) $ inducing an isomorphism of $\C$-algebras
$\C[\s^\vee  \cap M] \cong \C[M(\s)] \otimes_\C \C  [ \s^\vee_{N_\s}
\cap M_\s ]$
which defines  (non canonically) the product structure
\begin{equation} \label{product} 
Z_{\s, N} \cong \O_{\s, N} \times Z_{\s, N_\s}. 
\end{equation}

The map that sends a  cone $\s$ in $\Sigma$  to the orbit
$\O_\s \subset Z_\Sigma$ is a bijection between  the fan $\Sigma$ and the set of orbits.
If $\s $ is a face of $\t$ then  $Z_\s$ is an open subset of $
Z_\t$ and 
the orbit $\O_\t$ is contained in
the closure of $\O_\s$ in $ Z_\t$ since $\t^\bot \subset \s^\bot$,
thus the closure of the orbit of $\s$ in
$Z_\Sigma$ is $\overline{ \O _\s } = \bigcup \O_\t$
where $ \t$ runs through the cones of $\Sigma$ which have $\s$ as a face.

The orbit closures are normal toric varieties by themselves with
respect to the lattice $N(\s)$. The cones
of the fan 
associated to $\overline{ \O_\s }$ are of the form  $\t + (N_\s)_\R \subset N_\R / (N_\s)_\R $
for $\t \in \Sigma$  containing  $\s$ as a face.

\begin{Rem} \label{sinloc}
The singular locus of $Z_\Sigma$ is the union of those orbits 
$\O_\s$ for $\s$ a non regular cone.
\end{Rem}
\noindent
This follows from formula (\ref{product}) by noticing that the 
orbit $\O_\s$ is contained in the singular locus of $Z_\s$ if and only if 
$\o_{\s, N_\s}$ is a singular point of $Z_{\s, N_\s}$ if and only
the cone $\s$ is not a regular cone.

\begin{Defi} \label{sobrio}
A fan $\Sigma'$ is a {\it subdivision\index{fan subdivision}} of the fan
$\Sigma$ if both fans have the same support and if any cone of
$\Sigma'$ is contained in a cone of $\Sigma$.
The fan $\Sigma'$ is {\it regular subdivision} if $\Sigma'$ is a regular fan.
A regular subdivision $\Sigma'$ is a {\it  resolution of the fan $\Sigma$} if
any regular cone of $\Sigma$ belongs to $\Sigma'$.
\end{Defi}

\noindent
Associated to a
subdivision of fans there is a {\it
  modification\index{modification}} 
$ \p_\Sigma : Z_{\Sigma '} \rightarrow   Z_\Sigma$ inducing  an
 isomorphism between their tori.


\begin{Exem}
Let  $\Sigma$ be a regular subdivision of the cone $\s :=\R^{d+1}_{\geq 0}$ 
with lattice $N:= \Z^{d+1}$. 
This subdivision 
defines a modification  $ \p_\Sigma : Z_\Sigma \rightarrow
Z_{\s} = \C^{d+1} $ which  we describe in detail:
\end{Exem}

\noindent
The variety $Z_\Sigma$ is non singular, 
for each cone $\s$ of maximal dimension the variety $Z_\s $ is
isomorphic to $\C ^{d+1}$ and the restriction  $\p_\s : Z_\sigma \rightarrow \C^{d+1} $ of the morphism  
$\p_\Sigma $ 
is induced by the semigroup inclusion $ \R^{d+1}_{\geq 0} \cap M  \rightarrow \s^\vee  \cap M$.
The set of primitive vectors in the $1$-skeleton  $\s$ is  a basis of $N$ and its dual basis 
is a minimal set of generators of the semigroup $\s^\vee  \cap M$.
These generators give us coordinates to describe the map  $\p_\s : Z_\sigma \rightarrow \C^{d+1} $ in the form:
\begin{equation}  \label{chart}
\begin{array}{ccc}
X_1 &=& U_1^{a^{1}_1} U_2^{a^{2}_1} \cdots  U_{d+1}^{a^{d+1}_1} \\
X_2 &=& U_1^{a^{1}_2} U_2^{a^{2}_2} \cdots  U_{d+1}^{a^{d+1}_2}\\
& \dots &\\
X_{d+1} &=&  U_1^{a^{1}_{d+1}} U_2^{a^{2}_{d+1}} \cdots  U_{d+1}^{a^{d+1}_{d+1}}   \\
\end{array}
\end{equation}
where $(a^{i}_1, a^{i}_2, \dots, a_{d+1}^{i})$ is the 
coordinate of the primitive vector $a^{i}$ in the  $1$-skeleton of $\s$,
for $i =1, \dots, d+1$. Since the fan $\Sigma$ is regular, it is easy
to see directly from formula (\ref{chart}) that 
the map $ \p_\Sigma$ is an isomorphism over the torus $X_1 \cdots
X_{d+1} \ne 0$ of $\C^{d+1}$.

A {\it resolution of singularities\index{resolution of singularities}} of a variety $Z$ is a smooth
variety $Z'$ with  a  modification  $ Z' \rightarrow Z$
which is an isomorphism outside the singular locus of $Z$.
The resolution of singularities of normal toric varieties is reduced to a combinatorial property
of faces (see \cite{TE}).
More precisely we have that: 
Given any fan $\Sigma$ there is a resolution  $\Sigma'$  of $\Sigma$ (see definition \ref{sobrio}). 
The associated toric morphism $ Z_{\Sigma '} \rightarrow   Z_\Sigma$
is a resolution of singularities of the variety  $Z_\Sigma$
 (see \cite{Cox}, Theorem 5.1).

We describe now the exceptional locus associated to a subdivision
$\Sigma'$ of a fan $\Sigma$.
Taking away  the cone $\s$ from the fan of the cone $\s$ means geometrically to
take away the orbit $\O_\s$ from the variety $Z_\s$. It follows that (see \cite{GS-LJ} Proposition page 199):
\begin{equation} \label{exceptional}
\p^{-1}( \O_\s ) = \bigcup_{\t \in \Sigma', \stackrel{\circ}{\t} \subset 
  \stackrel{\circ}{\s} } \O_\t
\end{equation}
It follows from (\ref{exceptional}) that the {\it
  exceptional fibers},  i.e., the union of subvarieties of dimension $>1$  which are mapped to points, 
are given by:
$$ \bigcup_{\dim \s = \mbox{{\small rk}} N , \s \notin \Sigma' } \p^{-1} ( \O_\s)$$ and
that the {\it exceptional locus}, i.e., the subvarieties that are mapped on a
variety of smaller dimension, is:
$$ \bigcup_{\s \notin \Sigma' } \p^{-1} ( \O_\s) =  
\bigcup_{\t \mbox{ {\small minimal }} \in \Sigma' - \Sigma } \overline{\O_\t}.$$
The {\it discriminant locus}, i.e., the image of the exceptional locus, is equal to
\begin{equation} \label{discriminantlocus}
\bigcup_{\s \mbox{ {\small minimal }} \in \Sigma - \Sigma' } \overline{\O_\s}. 
\end{equation}

%

\subsection{Newton polyhedra and partial resolution of singularities}

The {\it Newton polyhedron\index{Newton polyhedron}} ${\cal N} (\f)$
of a non zero series 
$\f = \sum c_a X^{a} \in \C\{ X \} $ with $X=(X_1, \dots, X_{d+1})$
is the convex hull of the set 
$\bigcup_{c_a \ne 0} \  a  + \R_{\geq 0}^{d+1}$.
More generally the  Newton polyhedron of  any non-zero
germ $\f = \sum c_a X^{a}$ of
holomorphic function at the special point $\o_\r$ of a normal affine toric
variety $Z_\r = \makebox{Spec} \C [\r^\vee \cap M  ]$ (for a
strictly convex cone  $\r^\vee $)
is the convex hull of the subset 
$\bigcup_{c_a \ne 0} \  a  + \r^\vee$ of $M_\R$. We denote it by  ${\cal
  N}_\r (\f)$  or by  ${\cal N} (\f)$ if the cone $\r$ is clearly determined by the context.
Many of the properties associated with classical Newton polyhedra hold
in this case; for instance, if $0 \ne \f = \f_1 \cdots \f_s$  we have that ${\cal N} (\f) $ is the Minkowski sum 
${\cal N} (\f_1) +  \dots + {\cal N} (\f_1) $ since the series $
\f_i $ have coefficients in a domain. It follows from this property that:

\begin{Rem} \label{mdivisors}
If $0 \ne \f = \f_1 \cdots \f_s$ and  ${\cal N} (\f) $ has only one vertex 
the same holds for each of the Minkowski terms  ${\cal N} (\f_i) $ for $i=1, \dots, s$.
\end{Rem}

\noindent
The {\it face ${\cal F}_u$ of the polyhedron ${\cal N}_\r (\f)$  defined
  by a vector\index{face (of a polyhedron}} in 
$ u \in \r$
is the set of vectors $  v \in {\cal N}_\r (\f)$ 
such that  
$\langle u, v \rangle   = \inf_{v' \in {\cal N}_\r (\f)} \langle u, v' \rangle$.
All faces of the polyhedron ${\cal N}_\r (\f)$ can be recovered in this way.
The face of ${\cal N}_\r (\f)$  
defined by $ u$ is compact if and only if  $u \in \stackrel{\circ}{\r} $.

\noindent
The {\it cone $\s( {\cal F} ) \subset \r $  associated to the face $\cal F$}
of the polyhedron ${\cal N}_\r (\f)$ is 
$$\s( {\cal F} ) := \{ u \in \r  \; /
\; \forall v \in {\cal F}, \makebox{ we have } 
\langle u, v \rangle   = \inf_{v' \in {\cal N}_\r (\f)} \langle u, v'
\rangle \}.$$

\noindent
The cones $\s( {\cal F})$, for $\cal F$ running through the set
of faces of the polyhedron ${\cal N}_\r (\f)$, define a  subdivision $\Sigma ({\cal N}_\r (\f))$ of
the fan of the  cone $\r$ called the {\it dual Newton diagram}.
The relative
interiors of the cones in the fan $\Sigma ({\cal N}_\r (\f))$ are
equal to the equivalence classes of vectors in $\r$ by the relation: 
$ u \sim u'  \Leftrightarrow {\cal F}_u ={\cal F}_{u'}$.
We say that a fan $\Sigma$  supported on the cone $\r$ is 
{\it compatible\index{fan compatible with a polyhedron}} with a set of series
$\f_1, \dots, \f_s \in \C \{ \r^\vee \cap M \} $
if it  subdivides the fan $\Sigma ({\cal N}_\r (\f))$ with $\f = \f_1
\cdots \f_s$.
A cone in 
the fan $\Sigma ({\cal N}_\r (\f_1 \cdots \f_s))$ is  intersection
of cones of the fans $\Sigma ({\cal N}_\r (\f_i))$
therefore 
$\Sigma$ 
is compatible with all the polyhedra ${\cal N}_\r (\f_i)$.
If $\Sigma $ is compatible with  ${\cal N}_\r (\f)$ all vectors  in $\stackrel{\circ}{\s}$  
define the same face ${\cal F}_\s$ of  ${\cal N}_\r (\f)$, for $\s \in
\Sigma$. 

\begin{Defi} Let $0 \ne \f = \sum c_a X^a \in  \C \{ \r^\vee \cap M \}$. 
The {\it symbolic restriction\index{symbolic restriction}} $\f_{| \cal F}$ of $\f$ to
the compact face $\cal F$ of the polyhedron   ${\cal N}_\r (\f )$ is the polynomial 
$\f_{| \cal F} := \sum_{a \in {\cal F}}
c_a X^{a} \in \C[ \r^\vee \cap M] $.
The {\it Newton principal part} $  \f_{| \cal N}$ 
 of $\f$ is the sum of those terms of  $\f$ having exponents lying on the compact faces of the 
Newton polyhedron ${\cal N}_\r (\f)$ 
\end{Defi}
We follow here the terminology of \cite{Kou} and \cite{Oka}.
The Newton principal part  $  \f_{| \cal N} \in \C[ \r^\vee \cap M]$
does not change if we change the ring  $\C [[ \r^\vee \cap M ]]$ 
by extending the lattice $M$.

Let $\Sigma$ be any fan supported on $\r$ defining the
modification $\p_\Sigma: Z_\Sigma \rightarrow Z_\r$. 
Let ${\cal V}$ be a subvariety of $Z_\r$ such that the intersection of the 
discriminant locus of $\p_\Sigma$ with each irreducible component ${\cal V}_i$ of ${\cal V}$
is nowhere dense on  ${\cal V}_i$. For instance if ${\cal V}$ is irreducible this condition holds
if the torus is an open dense subset of ${\cal V}$. 
The {\it strict transform} ${\cal V}_\Sigma \subset Z_\Sigma $  is the subvariety of $\p_\Sigma^{-1} ({\cal V})$ 
such that the restriction $ {\cal V}_\Sigma  \rightarrow {\cal V}$ is a modification. 

If the fan  $\Sigma$ is  regular, 
the toric map $\p_\Sigma: Z_\Sigma \rightarrow Z_\r$
is a (toric) {\it embedded pseudo-resolution\index{toric embedded
    pseudo-resolution}} of ${\cal V}$ if the
restriction
 ${\cal V}_\Sigma \rightarrow {\cal V}$
is a modification such that the strict transform ${\cal V}_\Sigma$ is
non singular and transversal to the orbit stratification of the 
exceptional locus of 
$Z_\Sigma$.
The modification $\p_\Sigma $ is  a (toric) {\it
embedded resolution\index{embedded resolution}} of ${\cal
V} $ if the restriction  to the strict
transform $  {\cal V}_\Sigma  \rightarrow  {\cal V}$ is an isomorphism outside the
singular locus of ${\cal V}$ 
 (see  \cite{Rebeca}).  
If $\p_\Sigma $ is only a pseudo-resolution we can only guarantee that
the map $  {\cal V}_\Sigma  \rightarrow  {\cal V}$ is an isomorphism outside the intersection
of ${\cal V}$ with the discriminant
locus of $\p_\Sigma $. In this case, this set  contains the singular locus of ${\cal V} $
but it is not necessarily equal to it. 
\begin{Defi}
If $\Sigma$ is a (non necessarily regular) subdivision of $\r$
the toric morphism $\p_\Sigma: Z_\Sigma \rightarrow Z_\r$
is a {\it partial (toric) embedded resolution  of ${\cal V}$}  
if for any resolution  $\Sigma'$ of the fan  $\Sigma$ the map $\p_{\Sigma'} \circ \p_\Sigma $
is an embedded resolution  of ${\cal V}$.
\end{Defi}

Let $ {\cal V} \subset Z_\r$ an irreducible subvariety 
such that the intersection with the torus is an open dense subset.
Let $\Sigma$  be a subdivision $\r$ compatible with a set of generators  $\f_1, \dots,
\f_{s} $ of the ideal of $ {\cal V} \subset Z_\r$. 
We give a combinatorial condition on the Newton polyhedra of  $\f_1, \dots, \f_{s}$
for the intersection of the strict transform with the exceptional fiber being non empty.

\begin{Lem} \label{nondeg}
Let  $\s$ a cone in $ \Sigma$ such that
$\stackrel{\circ}{\s} \subset \stackrel{\circ}{\r}$.
If $\O_\s \cap {\cal V}_\Sigma \ne
\emptyset$ then  the face ${\cal F}_i$ of the Newton polyhedron
${\cal N} (\f_i)$ of $\f_i$ defined by $\s$ is  of dimension $\geq 1$ for  $1\leq i\leq s$.
\end{Lem}
{\em Proof.} We have that $ \f_{ i} -  \f_{ i {| \cal F}_i}  $ belongs to the ideal generated by
$\{ X^u /  u \in ({\cal N} (\f_i) - {\cal F}_i ) \cap M \}$. 
Since $\Sigma$ is compatible with the $\f_i$  the cone $\s^\vee$ contains
the cone spanned by elements in the polyhedron  $- u_0 + {\cal N}(\f_i)$ for any $u_0 \in {\cal F}_i$. 
Let  $u_i \in {\cal F}_i$ be a vertex  then 
we can factor in the ring  $\C [ \s^\vee \cap M ]$: 
 $$ \f_{ i} \circ \p_\s = X^{u_i} \psi_i \makebox{ and } 
 \f_{ i {| \cal F}_i} \circ \p_\s  = X^{u_i} \psi_{ i {| \cal F}_i}
\makebox{ with } \psi_{ i {| \cal F}_i} \in \C[\s^\bot \cap M]$$
in such a way that the exponent of a term appearing in 
$ X^{-u_i} ( \f_{ i} \circ \p_\s 
- \f_{ i {| \cal F}_i} \circ \p_\s )$ belongs to  $ (\s^\vee -
\s^\bot) \cap M$ and thus this term  vanishes on the orbit  $\O_\s$.
By definition the elements 
$X^{-u_i} \f  \circ \p_\s $ for $1\leq i\leq s$
belong to the ideal defining the strict transform of ${\cal V}$.
If the face ${\cal F}_i$ is a vertex
for some $i$ the ideal of $\O_\s \cap {\cal V}_\Sigma$ in $Z_\s$ is
equal to $(1)$ thus ${\cal V}_\Sigma \cap \O_\s$ is empty.
\hfill $\ {\diamond}$


The following lemma is an easy consequence of 
the implicit function theorem.

Let $\r \subset N_\R$ be a rational strictly convex cone of dimension equal to $\mbox{rk } N$.
We denote by $\Delta$ the cone 
$\r \oplus \R^g_{\geq 0} \subset (N_\D)_\R$ where $N_\D$ is the lattice 
$N \oplus \Z^g$ with dual lattice $M_\D$. 
The semigroup $\Delta^\vee \cap M_\D$ is of the form $(\r^\vee \cap M ) \oplus \Z^g_{\geq 0}$.
The monomial corresponding to $(\a , v) \in \Delta^\vee \cap M_\D$ is denoted 
by $X^\a U^v$ or by $X^\a U_1^{v_1} \dots U_g^{v_g}$.

\begin{Lem} \label{IFT}
If $\f_1, \dots, \f_g \in \C \{   \Delta^\vee \cap M_\D  \} $ verify that 
$
\f_1(\o_\s, U) = U_i, \mbox{ for } i=1, \dots, g
$
then 
there exist series 
$\epsilon_i \in \C \{ \r^\vee \cap M \} $ for $i=1, \dots, g$  
such that the ideals of $  \C \{   \Delta^\vee \cap M_\D  \} $
 generated by $\f_1, \dots, \f_g$
and $U_1 - \epsilon_1, \dots, U_g - \epsilon_g $
coincide.
\end{Lem}
{\em Proof.}
An homomorphism of semigroups 
$\Z^{s}_{\geq 0} \stackrel{\psi}{\longrightarrow}  \r^\vee \cap M $
extends to an homomorphism 
 $\Z^{s+ k }_{\geq 0} \stackrel{\psi \times \mbox{{\rm \small Id}}}{\longrightarrow}  \D^\vee \cap M_\D$.
If $\psi$ is surjective it defines an equivariant embedding   $Z_\r \subset \C^{s}$ which 
extends (by using the homomorphism
 $\psi \times \mbox{{\rm Id}}$) to an equivariant embedding $Z_\D = Z_\r \times \C^g  \subset \C^{s +g }$.
If $\varphi_1, \dots, \varphi_g$ are power series defining holomorphic functions at $(\C^{s +g }, 0)$   
representing $\f_1, \dots, \f_g$ 
the implicit function theorem 
guarantees the existence of power series $\varepsilon_i$
in $s$ variables such that 
the ideals 
$(\varphi_1, \dots, \varphi_g)$
and $(U_1 - \varepsilon_1, \dots, U_g - \varepsilon_g) $
coincide.
The  result follows by passing to the quotient   
by the  binomial ideal defining 
the embedding $Z_\D \subset \C^{s +g }$.
\hfill $\ {\diamond}$

\subsection{Embedded resolution of non necessarily normal toric varieties}

We build an embedded resolution of non necessarily normal affine  toric
variety $Z^\Lambda $ equivariantly embedded in a normal
affine toric variety $Z_\rho$ (for  $\r^\vee$ a strictly convex cone).
We build first a
{\it partial embedded resolution} which is a toric morphism
providing an {\it embedded normalization}
inside a normal toric ambient space.
Then any toric resolution of the
singularities of the ambient space, which always exists,  
provides
an embedded resolution.
The advantage is method is that the partial resolution is
completely determined by
the embedding $Z^\Gamma\subset Z_\rho$.
This result is the fruit of discussions with  Professor 
B. Teissier (see \cite{Teissier}, \S 6, Proposition 6.4 and \cite{cras}).

Let  $\Lambda$  be a  monoid.
An {\it equivariant embedding} of $Z^\Lambda$ in the  normal affine toric
variety  $ Z_\r$ is given by a surjective  homomorphism of  semigroups
$ \r^\vee \cap M \rightarrow \Lambda$
which extends to 
a lattice homomorphism
$\varphi:  M \rightarrow -
\Lambda + \Lambda $ and a vector space homomorphism
 $ \varphi_\R: M_\R \rightarrow 
(-\Lambda + \Lambda)_\R $.
The torus of $Z^\Lambda$ is equivariantly embedded in the torus of $ Z_\r$,
the embedding is obtained from the homomorphism $\varphi$.
The linear subspace $( \makebox{Ker} (\varphi_\R ))^\bot \subset N_\R$, denoted by $\ell$ in what follows, 
is of dimension equal to $\makebox{rk} \Lambda$ and the same holds for the cone 
$\s_0: = \ell \cap \r$. 
The ideal of the embedding  $Z^\Lambda \subset  Z_\r$ is generated by the binomials:
\begin{equation} \label{eqembedding}
X^u - X^v \in \C [ \r^\vee \cap M ] \mbox{ such that }  \varphi(u) = \varphi(v) \quad \mbox{ (see \cite{Sturmfels}, Chapter 4). } 
\end{equation}

\begin{Lem} \label{pseudo}
With the above notations suppose that the cone $\r^\vee$ is strictly convex.
Let $\Sigma$ be any fan compatible with a finite set of binomial equations 
$X^{u_i} - X^{v_i} =0$ for $i \in I$ defining the embedding $Z^\Lambda \subset Z_\r$. 
Then the fan $\Sigma$ is compatible with the linear subspace $\ell$. 
If $\s \in \Sigma $ and $\stackrel{\circ}{\s} \subset \stackrel{\circ}{\r}$
then $\O_\s \cap Z^\Lambda_\Sigma \ne
\emptyset$ implies that $\s \subset \ell$.
Moreover, if $\s \subset \ell$ and  $\dim \s = \dim
  \ell$ the intersection  ${Z}_\Sigma^\Lambda \cap \O_\s $  as schemes is the simple point $ o_\s $
and ${Z}_\Sigma^\Lambda \cap  Z_\s$ is isomorphic to $Z_{\s, N_\s}$.
If $\Sigma $ is regular the morphism 
$\p_\Sigma $ is an embedded pseudo-resolution of singularities of $Z^\Lambda \subset Z_\r$. 
\end{Lem}
{\em Proof.}
The cone $\s_0 = \r \cap \ell $ is associated to the Minkowski sum of compact edges of
 $ {\cal N} (X^{u_i} - X^{v_i}) $ for $i \in
I$ since  $\langle w, u_i \rangle = \langle w, v_i
\rangle, \, \forall i \in I$  if and only if $w \in \ell$.
Since the fan $\Sigma$ is compatible with the binomial equations of $Z^\Lambda$ it follows that 
a subdivision of $\s_0$ is contained in $\Sigma$, i.e., this fan is compatible with the linear subspace
$\ell$. 

We deduce by duality from the equality $\s_0 = \r \cap \ell$ that:
\begin{equation} \label{duali}
\s_0^\vee = \r^\vee + \ell^\vee = \r^\vee +  \ell^\bot =
\r^\vee + \mbox{Ker} (\varphi_\R) 
\end{equation}
Since the cone $\r^\vee$ is strictly convex formula
(\ref{duali}) implies that 
\begin{equation} \label{orthoduali}
\s_0^\bot = \mbox{Ker} (\varphi_\R) 
\end{equation}
and thus 
\begin{equation} \label{kerduali}
\mbox{Ker} (\varphi )  \subset \s_0 ^\vee \cap M
\end{equation}

Let  $\s \in \Sigma$ with $\stackrel{\circ}{\s} \subset \stackrel{\circ}{\r}$, 
since $\Sigma$ is compatible with the binomials $X^{u_i} - X^{v_i}$, 
the ideal generated by $1 - X^{u_i -v_i}$ (up to relabeling)
is contained in the ideal defining the strict transform ${Z}_\Sigma^\Lambda$ in the chart $Z_\s$.
Thus the variety $Z'$, defined by $  X^{u_i -v_i} - 1 =0$  for $i \in I$,  
contains ${Z}_\Sigma^\Lambda \cap Z_\s$.
Then we have:
\[
\begin{array}{c}
 Z' \cap \O_\s \ne \emptyset \Leftrightarrow \; \exists p \in
 Z_\s: X^{u_i - v_i} (p) =1 \; \forall i \in I,  X^{u} (p) =0  \;  \forall u
 \in (\s^\vee - \s^\bot) \cap M \Leftrightarrow \\
 u_i -v_i \in \s^\bot \cap M, \, \forall i \in I \Leftrightarrow
 \makebox{ Ker } (\varphi) \subset \s^\vee \Leftrightarrow \s \subset \r
 \cap \ell\\
\end{array}
\]
The chart  $Z_\s$ is isomorphic to  $\O_\s \times Z_{\s, N_\s}$ by formula (\ref{product}).

If $\s \subset \ell$ and  $\dim \s = \dim
  \ell$ we have that $\s^\bot = \s_0^\bot $ coincides with $ \mbox{Ker} (\varphi_\R )$ by (\ref{orthoduali}).
We deduce  an isomorphism 
\begin{equation} \label{wr}
Z' \cong \{ \o_\s \} \times Z_{\s, N_\s } \subset  Z_\s.
\end{equation}
from (\ref{product}) since the lattice $\s^\bot \cap M = \mbox{Ker} (\varphi) $ 
is generated by $\{ u_i - v_i \}_{i \in I}$. Therefore the variety $Z'$ is irreducible
and of dimension equal to $\mbox{rk } \Lambda$.  
We deduce from (\ref{wr}) that ${Z}_\Sigma^\Lambda$ 
intersects the orbit  $\O_\s$ transversally since the coordinate ring
of ${Z}_\Sigma^\Lambda  \cap \O_\s$ is $\C$.
Since ${Z}_\Sigma^\Lambda \cap Z_\s$ is a
subvariety 
of the irreducible variety  $Z'$ and both are of the same dimension
they coincide.

If $\Sigma$ is regular we deduce that 
${Z}_\Sigma^\Lambda  $ is smooth and intersects transversally the  orbit
stratification of the exceptional locus of  $Z_\Sigma$ 
thus $\p_\Sigma$ is an embedded  pseudo-resolution of $Z^\Lambda$.
\hfill $\ {\diamond}$

With the notations of lemma \ref{pseudo} we have:
\begin{Pro} \label{binomial}
Suppose that the cone $\r^\vee$ is strictly convex. 
Let $\Sigma$ be  a subdivision of $\r$ containing the cone $\s_0$.
\begin{enumerate}
\item
The strict transform $Z_{\Sigma}^\Lambda$ of $Z^\Lambda$ by the morphism $\p_{\Sigma}$ is 
isomorphic to $Z_{\s_0, N_{\s_0}}$ and the restriction 
$\p_{\Sigma} |  Z_{\Sigma}^\Lambda : Z_{\Sigma}^\Lambda  \rightarrow Z^\Lambda$
is the normalization map.
\item
The morphism 
$\p_{\Sigma}$ is a partial embedded resolution of $Z^\Lambda \subset Z_\r$.
\end{enumerate}
\end{Pro}
{\em Proof.} We keep notations of lemma \ref{pseudo}.
If we choose a splitting $M \cong \mbox{Ker} (\varphi) \oplus \mbox{Im} (\varphi)$
we obtain using (\ref{kerduali}) a semigroup isomorphism 
$$
\s_0 ^\vee \cap M \cong \mbox{Ker} (\varphi) \oplus \varphi (\s_0 ^\vee \cap M), 
$$
which corresponds geometrically to the isomorphism 
$Z_{\s_0} \cong \O_{\s_0} \times  Z_{\s_0, N_{\s_0}}$ of (\ref{product}).

We deduce from (\ref{duali}) that $\s_0^\vee = \varphi_\R^{-1} ( \varphi_\R (\r^\vee))$
and it follows that the semigroup
\begin{equation} \label{saturation}
\varphi (\s^\vee_0 \cap M ) = \varphi_\R (\r^\vee)  \cap \varphi (M) 
\end{equation}
is the saturated semigroup $ \R_{\geq 0} \Lambda \cap (- \Lambda + \Lambda)$ 
of $\Lambda$ in the lattice it spans; therefore
the variety $Z_{\s_0, N_{\s_0}} $ is isomorphic to the normalization of $Z^\Lambda$
(see \cite{TE}).

Let $\Sigma'$ be a subdivision of $\Sigma$ compatible with the equations of $Z^\Lambda$.
By lemma \ref{pseudo} if $\s \in \Sigma' $, $\stackrel{\circ}{\s} \subset \stackrel{\circ}{\r}$
 and $\O_\s \cap Z^\Lambda_{ \Sigma'} \ne \emptyset$ then we have $\s \subset \ell$.
A fortiori the same property holds replacing $\Sigma'$ by $\Sigma$ 
as a consequence of (\ref{exceptional}).
It follows that the strict transform of the germ $(Z^\Lambda, o_\r)$ is contained in the chart
corresponding to the cone $\s_0$. This implies that $Z^\Lambda_\Sigma \subset Z_{\s_0}$ since
the morphism $\p_{\Sigma}$ is equivariant and $Z^\Lambda$ is equivariantly embedded.
It follows also from the proof of lemma \ref{pseudo} that 
the restriction of $\p_\Sigma $ to  $Z^\Lambda_\Sigma \rightarrow Z^\Lambda$
corresponds algebraically to the inclusion of $\C[ \Lambda ]$ in its integral closure
thus it is the normalization map.

A resolution  $\Sigma'$ of the fan  $\Sigma$ is subdivided by a regular fan $\Sigma''$
which is compatible with the equations of $Z^\Lambda$.
By lemma \ref{pseudo} the map $\p_{\Sigma''} \circ \p_{\Sigma'} \circ \p_{\Sigma }$ 
is a pseudo-resolution of $Z^\Lambda$. A fortiori the same holds for $p_{\Sigma'} \circ \p_{\Sigma }$ 
by (\ref{exceptional}).
By definition if $\s' \in \Sigma$
is a regular cone then $\s' \in \Sigma'$, thus $Z_{\Sigma'} \rightarrow Z_{\Sigma}$ is an isomorphism over 
the points of the $\O_{\s'}$. 
By  remark  \ref{sinloc}
the singular locus of ${Z}_{\Sigma}^\Lambda$ is defined by the intersection of those orbits 
$\O_{\s'}$ for those cones ${\s'} $ running through the set of non regular faces of $\s_0$.
 This shows that 
${Z}_{\Sigma'}^\Lambda \rightarrow {Z}_{\Sigma}^\Lambda $ is a resolution of
singularities of the normalization $ {Z}_{\Sigma}^\Lambda  $ of ${Z}^\Lambda$. 
A fortiori the map ${Z}_{\Sigma}^\Lambda \rightarrow {Z}^\Lambda $
is a resolution of
singularities.
\hfill $\ {\diamond}$

\subsection{Equivariant branched coverings of normal toric varieties} \label{Galois}

Some branched coverings of normal toric varieties are equivariant.
Typically, 
if $\s$ is a  rational cone for $N$ it  is also rational for  $N' \subset N$ a sub-lattice of the same rank 
and  we have a homomorphism
of semigroups 
$ \s^\vee  \cap M  \rightarrow  \s^\vee  \cap M '$
where $M \subset M'$ is the inclusion of lattices  corresponding to  $N' \subset N$
by duality.  
This homomorphism defines an equivariant morphism
\begin{equation} \label{group}
Z _{\s, N'}  \rightarrow Z_{\s, N }
\end{equation}
extending  the homomorphism of tori 
$T' \rightarrow T$ defined by the lattice extension $M \subset M'$,   
which has kernel  a finite subgroup $H$ of $T'$.
Each $w \in H$ corresponds to a morphism
$Z _{\s, N'} \rightarrow Z _{\s, N'} $ given by the homomorphism 
$\C[ \s^\vee  \cap M '] \rightarrow \C[ \s^\vee  \cap M ']$ mapping 
$X^u \mapsto w(u) X^u$. The ring $\C[ \s^\vee  \cap M ]$ is the set of invariants
of  $\C[ \s^\vee  \cap M ']$ by the action of the group $H$ and the morphism
(\ref{group}) coincides with canonical projection of the  quotient of $Z_\s'$
with respect to the action of the group $H$ by Corollary 1.16 of 
\cite{Oda}.
If $\s$ is of maximal dimension the  $0$-orbit $\o_{\s}'$ of $Z _{\s, N'}$  projects
to the $0$-orbit  $\o_{\s}$ of $Z_{\s, N}$ and we have that $(Z_{\s, N'},  \o_\s')
\rightarrow (Z_{\s, N}, \o_\s)$ is a morphism of analytically irreducible
germs. The corresponding homomorphism of analytic algebras 
$\C \{ \s^\vee  \cap M \}  \rightarrow \C \{ \s^\vee  \cap M' \}$ extends
to a homomorphism $L \rightarrow L'$ of their fields of fractions
of degree equal to the cardinality of $H$, i.e., 
the index of $M $ as a subgroup of  $M'$.
This field extension is Galois and the Galois group is obtained from 
the automorphisms of  $\C \{ \s^\vee  \cap M' \}$ defined by the elements of
$H$ (see \cite{GP}).

Let $\nu_1, \dots, \nu_g \in M'$ and define from them a sequence of lattices and integers:
\begin{equation}  \label{ninu}
\left\{
\begin{array}{c} 
M_0 := M, M_i: = M_{i-1} + \nu_i \Z , \makebox{ for } i=1, \dots, g
\\
n_0 :=1, n_i = \# M_i / M_{i-1}        \makebox{ for } i=1, \dots, g
\end{array} 
\right.
\end{equation}
The lattices $M_i$ are all sub-lattices of finite index of $ M'$.
We have the inclusions of lattices  $N' \subset N_g \subset  \dots \subset N_1 \subset N_0 =N $ where  
$N_i$ denotes the dual lattice of $M_i$.

\begin{Lem} \label{field}
The  field of fractions of $  \C \{ \r^\vee \cap M_j \} $ is 
$ L[X^{\nu_1}, \dots, X^{\nu_j}]$.
If $\l \in \r^\vee \cap M'$ then   $X^\l \in \makebox{{\rm Fix}} (\makebox{{\rm Gal}} 
( L' / L [ X^{\nu_1}, \dots, X^{\nu_j}]))$  if and only if 
$\l \in \r^\vee \cap M_j$.
\end{Lem}
{\em Proof.}
The homomorphism of analytic algebras
$ \C \{ \r^\vee \cap M \} \rightarrow  \C \{ \r^\vee \cap M_j \}$
is finite and defines an extension of the corresponding fields of fractions of degree $n_1 \cdots n_j $ 
equal to
the order of the finite group $M_j/ M $.
We prove the first assertion by induction on $j$:
for $j=1$ the roots of the minimal polynomial  of $X^{\nu_1}$ over $L$ are
the different conjugates of $ X^{\nu_1}$ by the action of 
the elements of the Galois group of $L/L'$. 
We deduce from this that the minimal polynomial of $X^{\nu_1}$ is $Y^{n_1}- X^{n_1 \nu_1}$
where $n_1 = \# M_1 / M_0$ is  also the degree
of the extension $ L[ X^{\nu_1}]/L$. Since $ L[ X^{\nu_1}] $ is contained in the field of fractions of 
$\C \{ \r^\vee \cap M_1\}$ and both fields define extensions of $L$ of the same degree
they are equal. By induction hypothesis the  field of fractions of $  \C \{ \r^\vee \cap M_{j-1} \} $ is 
$ L[X^{\nu_1}, \dots, X^{\nu_{j-1}}]$ and we can replace $L $,  $\nu_1$ and $n_1$  in the previous argument by 
$L[X^{\nu_1}, \dots, X^{\nu_{j-1}}]$,  $\nu_j$ and $n_j$ respectively to obtain the assertion for $j$.

If $\nu \in \r^\vee \cap M_j$ 
it is clear that $\nu$ is fixed by any element of the Galois group of the extension
$L' / L [ X^{\nu_1}, \dots, X^{\nu_j}]$. The converse follows by the first assertion and 
Corollary 1.16 of \cite{Oda} applied to 
the inclusion of  semigroups $\r^\vee \cap M_j \subset \r^\vee \cap M'$.
\hfill $\ {\diamond}$


\subsection{A reminder on toroidal embeddings}

Let $ \cal X$  be  a normal variety of dimension $d+1$, and  let $E_i$ be a finite set of normal  hypersurfaces 
with complement  $\cal U$ in  $ \cal X$.
A {\it  toroidal embedding without self intersection \index{toroidal embedding}} is defined by requiring 
the triple $({ \cal X}, {\cal U}, x)$ at any point 
$ x \in { \cal X}$ to be formally isomorphic to $(Z_\s, T=(\C^*)^{d+1}, z)$ 
for $z $ a point in some toric variety $ Z_\s$.
This means that there is a formal isomorphism between the completions of the 
local rings at respective points which sends the ideal of ${ \cal X} - {\cal U}$ into the ideal
of $Z_\s - T$; (see \cite{TE}).
The variety $\cal X$ is naturally stratified, with 
strata $\bigcap_{i \in K} E_i - \bigcup_{i \notin K} E_i$ and the open stratum $\cal U$.

The {\it star of a stratum\index{star of a stratum}} $\mathfrak  S$, $\makebox{star} \, {\mathfrak S}$,  is the union of the  strata containing $ \mathfrak S$  
in their closure.
We associate to the stratum $ \mathfrak S$ the set $M^{{\mathfrak S}}$ of Cartier divisors supported on  
$\makebox{star} \, {\mathfrak S} -{\cal U}$. We denote by    $N^{{\mathfrak S}}$ the dual group 
$\makebox{ Hom} (M^{{\mathfrak S}}, \Z)$. The semigroup of effective divisors defines in 
the real vector space $M^{{\mathfrak S}}_\R := M^{{\mathfrak S}} \otimes \R$ a rational convex polyhedral cone and 
we denote its dual cone in $N^{{\mathfrak S}}_\R := N^{{\mathfrak S}} \otimes \R$  by $\r^{{\mathfrak S}}$.
If $ \mathfrak S'$ is a stratum in $\makebox{star} \, {\mathfrak S}$, 
we have a group  homomorphism  defined by restriction of Cartier divisors  
$M^{{\mathfrak S}} \rightarrow M^{{\mathfrak S'}} $ which is onto; by duality we obtain an inclusion
$N^{{\mathfrak S'}} \rightarrow N^{{\mathfrak S}} $ and the cone $\r^{{\mathfrak S'}}$ is mapped onto a face of
$\r^{{\mathfrak S}}$ (see \cite{TE}).
We can associate in this way to a toroidal embedding without self-intersection
a {\it conic polyhedral complex with integral structure\index{conic polyhedral complex with integral structure}} 
({\it   c.p.c.} in what follows)
see \cite{TE}.
 This generalizes  
the way of recovering from a normal toric variety the associated fan. 
This complex is {\it combinatorially isomorphic} to the cone over the  dual graph of 
intersection of the divisors $E_i$.
We have that the strata of the stratification are in one-to-one correspondence with 
the faces of the conic polyhedral complex.
For instance,  the conic polyhedral complex associated 
to the toroidal embedding defined by $Z_\Sigma$ and the normal hypersurfaces   
$ \{ \overline{\O_\s} \}_{\s \in \Sigma^{(1)} }$
is isomorphic to the 
conic polyhedral complex (with integral structure) $(\Sigma, N)$
defined by the fan $\Sigma$ and the lattice $N$.



We can define, in an analogous manner to the  case of a fan, a regular subdivision of 
a conic polyhedral complex.  Associated to a subdivision we 
have an induced { \it toroidal modification\index{toroidal modification}} (see \cite{TE} Th. 6* and 8*),  
i.e, a normal variety ${\cal X}'$
with a  toroidal embedding ${\cal U} \subset {\cal X}'$ and a modification 
${\cal X}' \rightarrow {\cal X}$
provided with a commutative diagram:
$$\begin{array}{ccc}
{\cal U} &  \rightarrow   & {\cal X}
\\
\downarrow & \nearrow &
\\
{\cal X}'  &    &  
\end{array}$$
The notion of toric partial embedded  resolution  
generalize easily in the toroidal case.

\section{Toric quasi-ordinary singularities} \label{Quasi} 
 
We introduce toric quasi-ordinary 
singularities and we extend
to this case many notions and properties of quasi-ordinary 
singularities.

Let $(S, o)$ be a germ of analytically irreducible complex variety of dimension $d$.
We denote by $R$ the associated analytic algebra.
 A sufficiently small representative $S \rightarrow S'$ of 
a finite map germ  $(S, o) \rightarrow (S', o')$
has 
finite fibers, its image is an open neighborhood of $o'$ 
and 
the maximal cardinality of the fibers  is equal to 
to the {\it degree} of the map.
The {\it discriminant locus\index{discriminant locus}}, i.e.,  the set of points of 
having fibers of cardinality less than the degree, is an analytical subvariety of $S'$. 
Outside the discriminant locus, the map is an {\it unramified\index{unramified covering}} covering.
  We can think of the discriminant locus  
as an analytic space or as a germ at $o'$.


\begin{Defi}
A germ of complex analytic variety $(S, o)$ is a {\it quasi-ordinary singularity\index{quasi-ordinary singularity}} 
if there exist a finite morphism $(S, o) \rightarrow (\C^d,0)$
(called a {\it quasi-ordinary projection)\index{quasi-ordinary projection}}
and some analytical coordinates $(X_1, \dots, X_d)$ at $0$, such that
the morphism is unramified over the torus $X_1 \dots X_d \ne 0$ in a
neighborhood of the origin.
\end{Defi}

The class of  quasi-ordinary singularities contains  all
curve singularities. 
The
Jung-Abhyankar Theorem guarantees that 
$R $ can be viewed as a subring of $\C \{ X_1^{1/m}, \dots, X_d^{1/m}\}$
for some suitable integer $m$
(see \cite{Jung} for a topological proof in the surface case and 
 \cite{Abhyankar}, Th. 3 for an algebraic proof).

The finite  map germ $(S, o) \rightarrow (S', o')$
corresponds algebraically to a local homomorphism $R' \rightarrow R$  of their analytic algebras 
which gives $R$ the structure of finite module over $R'$.
In particular if $R$ is generated over $R'$  by one element
there is a surjection $R'[Y] \rightarrow R$ which corresponds geometrically to an embedding 
$ (S, o)    \subset (S'  \times \C , (o', 0))$. 
We say that $ (S, o) $ is an {\it hypersurface 
relative to the base $(S', o')$}. 
We define toric quasi-ordinary singularities by replacing the base $(\C^d, 0)$ by 
the germ $(Z_\r, o_\r)$  of an affine  toric variety at its zero orbit
(for a   strictly convex cone  $\r^\vee$).

\begin{Defi} (see \cite{GP})
The germ  $(S, o)$ is a {\it toric quasi-ordinary singularity}
if there exist a finite morphism $(S, o) \rightarrow (Z_\r, \o_\r)$
unramified over the torus in a neighborhood of the zero-orbit $
\o_\r$ of a suitable normal affine toric variety $Z_\r$.
\end{Defi}

\begin{Rem}
The classical quasi-ordinary singularities are obtained when $(\r, M) = (\R^d_{\geq 0}, \Z^d)$.
\end{Rem}

By definition the analytic algebra $R$  of a toric quasi-ordinary singularity is 
a $\C \{ \r^\vee \cap M \}$-algebra of finite type. The germ  $(S, o)$ is an hypersurface 
relative to the toric base  if there exists $x \in R$ such that $R = \C \{ \r^\vee \cap M \} [x]$.
Then the $\C \{ \r^\vee \cap M \}$-algebra homomorphism 
$\C \{ \r^\vee \cap M \}[Y] \rightarrow R$ that maps $Y \mapsto x$ is surjective.
Its kernel is a principal ideal generated by a monic polynomial $f$ such that $f(\o_\r , Y)= Y^{\deg f}$ and 
$\deg f$ is equal to the degree of the map $(S, o) \rightarrow (Z_\r, \o_\r)$.
The polynomial $f$ is a {\it quasi-ordinary polynomial}, i.e.,  the discriminant 
 $\D_Y f$ of $f$ with respect to $Y$
is of the form:
$$ \D_Y f = X^\eta  H \makebox{ with $H(\o_\r) \ne 0$}.$$
Conversely each monic quasi-ordinary polynomial $f \in \C \{ \r^\vee \cap M \} [Y]$
such that $f(\o_\r , Y)= Y^{\deg f}$ defines a germ of toric quasi-ordinary hypersurface.
The  $\C \{ \r^\vee \cap M \}$-algebra homomorphism 
$\C \{ \r^\vee \cap M \}[Y] \rightarrow R$
defines an embedding $ S \subset Z_\r \times \C$ 
that maps $o \mapsto (\o_\r, 0)$. The quasi-ordinary 
projection of $(S,o)$ is induced by the first projection of the product $Z_\r \times \C$.

The product $Z_\r \times \C$ is the  toric variety $Z_\varrho$ defined 
by the cone $\varrho = \r \times \R_{\geq_0}$ with respect to the lattice  $N'$ dual to 
the lattice $M' := M \oplus  y \Z$. Then we have 
$ \varrho^\vee \cap M' \cong  (\r^\vee \cap M ) \oplus y \Z_{\geq 0}$. 
We denote the monomial corresponding to $u + s y \in (\r^\vee \cap M ) \oplus y \Z_{\geq 0}$ by $X^u Y^s$.



If $f$ is an irreducible quasi-ordinary polynomial the associated analytic algebra $R $ is 
the domain
$R = \C \{  \r^\vee \cap M  \} [Y]/(f)$.
There exists 
a fractional power
series $\z \in  \C\{     \r^\vee \cap \frac{1}{n} M     \}$  
which is a root  of $f$ where 
${n}$ is the degree of $f$ 
(see Th\'eor\`eme 1.1 and Remarque 1 of \cite{GP}). The inclusion  $\C \{  \r^\vee \cap M \}  \subset
\C\{ \r^\vee \cap  \frac{1}{n} M  \}$ corresponds to a 
branched covering of a normal affine toric variety and 
defines a Galois  extension $L \subset L_n $
of their corresponding fields of fractions
(see subsection \ref{Galois}).
The minimal polynomial of the root  $\z$ over the field  $L$   is equal to $f$, 
we have $R \cong \C  \{  \r^\vee \cap M  \} [\z]$  and the 
field of fractions of $R$ is  $L[\z]$ since $\z$ is finite over $L$.
The conjugates $\z^{(i)}$ of $\z$ by the action of the Galois group of $L \subset L_n$ 
define all the roots of $f$ since the extension $L[\z] \subset L_n$
is Galois.

We call {\it (toric) quasi-ordinary branches
\index{quasi-ordinary branches}}
the roots of (toric) quasi-ordinary polynomials.

If $f$ is a reduced quasi-ordinary polynomial of degree $n$ then it splits on $\C\{ \r^\vee \cap \frac{1}{n!} M \}$.
The difference $\z^{(s)} - \z^{(t)}$ of two different roots of $f$ 
divides the discriminant of $f$ on the ring $\C\{ \r^\vee \cap \frac{1}{n!} M \}$.
By remark \ref{mdivisors}, 
the Newton polyhedron of  $\z^{(s)} - \z^{(t)}$
has only one vertex therefore
$\z^{(s)} - \z^{(t)}$ is of the form 
$X^{\l_{st}}  H_{st} $ where $H_{st}$ is a unit in $\C\{ \r^\vee \cap \frac{1}{n!} M \}$.
It follows that the irreducible factors of $f$ are quasi ordinary polynomials.
The monomials $X^{\l_{st}} $ so obtained are called  {\it
characteristic monomials\index{characteristic monomials}} 
and the exponents  $ \l_{st} \in \r^\vee \cap \frac{1}{n} M $ are called the 
{\it characteristic exponents\index{characteristic exponents}}.
If $\mbox{{\rm rk}} M =1$ and if $f$ is irreducible the characteristic exponents 
correspond to the classical Puiseux characteristic exponents 
in arbitrary coordinates.
We do not need the classical argument 
to define the characteristic monomials which  
uses the factoriality of the ring $\C \{ X_1, \dots, X_d \}$ (see \cite{Lipman1}), a property 
which does not hold for the rings of the form $\C \{ \r^\vee \cap M \}$ in general.
The notion of characteristic monomials in the classical quasi-ordinary
case is already present in  Zariski's work (see \cite{Licei}); in the analytically irreducible hypersurface case 
many geometrical and  topological features of these singularities are determined in terms of the characteristic monomials 
by Lipman, Luengo, Gau and others (see \cite{Lipman0},  \cite{Lipman1}, \cite{Lipman2}, \cite{Luengo} 
and  \cite{Gau}).

We define a partial order $\leq_\r$ (or $\leq$ for short) on the cone $\r^\vee$:
$$ u \leq_\r u' \Leftrightarrow u' \in u + \r^\vee  \Leftrightarrow \, \forall w \in \r : \langle  u'- u, w \rangle \geq 0.$$ 
We can extend this partial ordering to a total one on the subset $\r^\vee \cap M$ by taking
 an {\it irrational} vector $\eta \in \r$, 
i.e.,  the coordinates of $\eta$  with respect to any base of the lattice $N$ 
are linearly independent
over $\Q$, and defining then $\leq_\eta$ by 
$ u \leq_\eta u' \Leftrightarrow \langle \eta, u - u' \rangle \leq 0$.

\begin{Lem} \label{rorder} (see  \cite{Licei} and \cite{Lipman2} in the classical case).
Let $f_1$ be an irreducible factor of the reduced toric quasi-ordinary polynomial $f$.
If $f_1 (\z^{(s_0)}) = 0 $ then we have:
$$
\left\{ \l_{s_0t} / \z^{(s_0)} \ne \z^{(t)}, \; f (\z^{(t)}) = 0  \;  \right\}  =
\left\{ \l_{st} / \z^{(s)} \ne \z^{(t)}, \; f (\z^{(t)}) = 0 \mbox{ and } f_1 (\z^{(s)}) = 0 \right\}
$$
and this set is totally ordered by  $\leq_\r$.
\end{Lem}
{\em Proof.} 
The equality above follows  since the extension $L[ \z^{(s_0)} ] \subset L_{n!}$ is Galois
and the elements of the Galois group act on a series in $\C \{ \r^\vee \cap \frac{1}{n!} M \}$
by changing the coefficients of its terms.
Then, if $\z^{(t)} \ne \z^{(t')} $ are roots of $f$ different to $\z^{(s_0)}$ we have that:
$$
X^{ \l_{t't}} H_{t't} = \z^{(t')} - \z^{(t)} = \z^{(t')} - \z^{(s_0)} - (\z^{(t)}- \z^{(s_0)})
=   X^{ \l_{t's_0}} H_{t's_0} - X^{ \l_{ts_0}} H_{ts_0}.
$$
Therefore 
$\l_{t't} = \min_\r \{ \l_{t's_0}, \l_{ts_0} \}$ and the assertion follows.
\hfill $\ {\diamond}$

\begin{Defi} \label{compara} (see \cite{Tesis}) 
Two irreducible quasi-ordinary polynomials $f^{(i)}$ and $f^{(j)}$
have order of coincidence $\l_{(i,j)}$ if their product  $f^ {(i)} f^{(j)}$ is a quasi-ordinary polynomial and  
$\l_{(i,j)}$ is the largest exponent of the set 
$\{\l_{st}  /   f^{(i)} (\z^{(s)}) = 0,    f^{(j)} (\z^{(t)}) = 0  \}$.
\end{Defi}
We say that the order of  coincidence of $f^{(i)} $ with itself is $\l_{(i, i)}:=+ \infty$.
We deduce from the proof of lemma \ref{rorder} and definition \ref{compara} the following 
property: 
\begin{Lem} \label{rval} 
If $f = f^{(1)} \cdots f^{(r)}$ is the factorization of a quasi-ordinary polynomial 
with monic irreducible factors we have that: 
$\min \{ \l_{(i, j)}, \l_{(j, l)} \} \geq \l_{(i, l)}$
with equality if $\l_{(i, j)} \ne \l_{(j, l)}$ for $i, j, l \in \{ 1, \dots, r \}$.
\hfill $\ {\diamond}$
\end{Lem}

In particular when $f$ is irreducible it follows that the set of characteristic exponents is totally ordered 
by $<_\r$ (see \cite{Lipman1}).
In this case we relabel the characteristic exponents  by $ \l_1 <_\r \l_2  <_\r \dots <_\r \l_g $ and we denote 
$\l_{g+1} = + \infty $. 
Following Lipman (see \cite{Lipman2}, page 61) we associate to the  characteristic exponents  sequences of  lattices and integers.
In the plane branch case the sequence of  integers coincide with the
first component of the characteristic pairs in arbitrary coordinates. 
\begin{Defi} \label{ni} 
The  lattices $M_i$ and the integers $n_i$ associated to the  sequence
of characteristic exponents $\l_1, \dots, \l_g$ for $i=0, \dots, g$ by formulae (\ref{ninu})
are called characteristic.
\end{Defi}
We denote by  $e_{i-1} = n_i \cdots n_g $, for $i= 1, \dots , g$ and we set $n_0 :=1$.
We denote by $N_g \subset \cdots \subset N_1 \subset N_0 =N$ the 
sequence of dual lattices of $M= M_0 \subset \cdots \subset M_g$.

If $f$ is reduced the set of characteristic exponents is not totally ordered by $\leq _\r$, 
for example the characteristic exponents $ (1,0), (\frac{3}{2},0), (1, \frac{3}{2})$ of
$f = (( Y - X_1)^2 - X_1^3)) (( Y + X_1)^2 - X_1^2 X_2^3))$
are not totally ordered for $\leq_{\R^2_{\geq 0} }$.


\begin{Lem} (see \cite{Lipman1}) \label{mfield} If $f$ is an irreducible toric quasi-ordinary polynomial and  if 
$\z$ is a root of $f$ we have:
 \begin{enumerate}
\item The characteristic integers  $n_i$ verify that $n_i > 1$ for $i=1, \dots, g$ and $n_1 \cdots n_g = \deg f$.
\item The field of fractions of $R$  is equal to $L[\z] = L[X^{\l_1}, \dots, X^{\l_g}].$
\end{enumerate} 
\end{Lem}
{\em Proof.}
Let  $\z'$ be a  conjugate of $\z$ by an element of the Galois group of the field extension 
$L_n \supset  L [ X^{\l_1}, \dots, X^{\l_j}]$.
If $\z' \ne \z$  we have $\z' - \z = X^{\l_k} H_k $ for a unit $H_k$ and $k > j$ 
(since $X^{\l_1}, \dots, X^{\l_j}$ are fixed for this Galois group). 
In particular for $j= g$ the only possibility is $\z' =\z$ thus 
$\z \in   L [ X^{\l_1}, \dots, X^{\l_g}]$ since the extension 
$L_n \supset L [ X^{\l_1}, \dots, X^{\l_j}]$ is Galois.
Conversely any element of the Galois group of the extension $L_n \supset L [ \z]$ fix $\z$ and therefore 
all the terms appearing in $\z$, in particular   $X^{\l_1}, \dots, X^{\l_g}$,  belong to $ L [ \z]$ since the 
extension $L_n \supset L[\z]$ is Galois. It follows that  $n_i > 1$ for $i=1, \dots, g$,
and that the degree $n$ of the extension  $L[\z] \supset L$ is equal to $n_1 \cdots n_g$.
\hfill $\ {\diamond}$

We have the following conditions for a power series $\z \in  \C \{ \r^\vee \cap \frac{1}{n} M \}$
to be a quasi-ordinary branch (see \cite{Lipman1}, prop. 1.5 or \cite{Gau}, prop 1.3 in the classical case).
\begin{Lem}\label{expo}
Let $\z= \sum c_{\l} X^\l $ be a non unit in  $\C \{ \r^\vee \cap \frac{1}{m} M  \} $.
Then the minimal polynomial of $\z$ over $\C \{ \r^\vee \cap  M \}$ is quasi-ordinary 
if and only if there exists 
elements  $\l_i \in \r^\vee \frac{1}{m} M  $ for  $1 \leq i \leq g$ such that
\begin{enumerate}
\item
$ \l_1 <_\r \l_2  <_\r \dots <_\r \l_g $, and $ c_{\l_i } \ne 0$  for $1 \leq i \leq g$.

\item 
If $c_{\l} \ne 0$  then $\l$ is  the sub-lattice  
$M + \sum_{\l_i \leq_\r \l} \Z \l_i  $ of $M_\Q$.

\item
$\l_j$ is not in the sub-lattice 
$M + \sum_{\l_i <_\r \l_j } \Z \l_i  $, of $M_\Q$ for $j= 1, \dots, g$.
\end{enumerate}
If such elements exist
they are uniquely determined by $\z$ and they are the characteristic exponents
of $\z$.
\end{Lem}
{\em Proof.}
If the minimal polynomial of $\z$ over $\C \{ \r^\vee \cap  M \}$ is quasi-ordinary
then the result follows from lemmas \ref{rorder}, \ref{mfield}
and \ref{field} applied to sequence of characteristic exponents.
Conversely, if $\z' $ is the conjugate of $\z$ by an element of the Galois group of 
$L_n \supset L $ and if $\z \ne \z'$  let us 
consider the sequence of lattices $M_i$ and integers $n_i$  associated 
to $\l_1, \dots, \l_g$ by (\ref{ninu}). 
There is some $j \geq 1$ such that 
the monomials $X^{\nu}$ are fixed  for $\nu \in M_{j-1}$ and $X^{\l_j}$ is not fixed by this element by 
lemma \ref{field} and hypothesis 3. Then hypothesis 1 and 2 imply that 
the difference $\z' - \z$
is of the form  $\z' - \z = X^{\l_j} H_j$ for a unit $H_j$. 
\hfill $\ {\diamond}$

\begin{Rem} \label{r}
The characteristic lattices associated to $f$     provide a 
canonical way of writing the terms of its roots:
$$\z =  p_0 +  p_1 + \dots +  p_g, $$
where  
$p_0 \in \C \{ \r^\vee \cap M \} $ and 
the monomial $ X^{\l}$ appears in the summand   $p_j$
implies that $\l_j \leq_\r \l  $ and  $\l_{j+1} \nleq_\r \l $ for $j=1, \dots , g$.
\end{Rem}
  
It is shown by Lipman (see \cite{Lipman2}, remark 7.3.2) that 
an analytically irreducible quasi-ordinary hypersurface germ 
of dimension $d$ is {\it normal \index{normal}} if and only
if it is isomorphic to a germ of the form $Y^{n} - X_1 \dots X_c = 0$
for some $1 \leq c \leq d$;
otherwise it is well known that its normalization is a 
{\it quotient singularity} (see \cite{Lipman2});
in the two dimensional case it is the germ of an
affine toric  surface 
(see \cite{BPV}, Chapter III, Theorem 5.2). 
In \cite{Tesis} is proved that the normalization of an irreducible 
quasi-ordinary hypersurface germ is isomorphic to the germ of an affine 
normal toric variety at its zero orbit and that this singularity 
is  determined from the set of characteristic exponents.
The following proposition generalizes this fact 
for toric quasi-ordinary hypersurface germs.

\begin{Pro} \label{closure}
The integral closure of the ring $R$ in its field of fractions is equal to $\C \{ \r^\vee \cap M_g \}$.
\end{Pro}
{\em Proof.}
The analytic algebra of the quasi-ordinary hypersurface is of the form $R = \C   \{ \r^\vee \cap M\} [ \z]$.
By lemma \ref{expo} we have a ring extension $R \subset \C \{ \r^\vee \cap M_g \}$ which is integral since 
$\C   \{ \r^\vee \cap M_g\}$ is integral over $\C   \{ \r^\vee \cap M\}$. By lemmas \ref{field} and \ref{mfield}
the rings $R$ and $ \C \{ \r^\vee \cap M_g \}$ have the same field of fractions. These two conditions imply that 
both rings have the same integral closure over their field of fractions. The 
ring $\C \{ \r^\vee \cap M_g \}$ is integrally closed 
since it is the analytic algebra of the normal variety  
$Z_{\r, N_g}$ at the  point $\o_\r$. 
\hfill $\ {\diamond}$

\subsection{The Eggers-Wall tree of a reduced quasi-ordinary polynomial} \label{Eggers}

We structure the partially ordered set of characteristic monomials of a reduced 
toric quasi-ordinary polynomial in  a labeled tree.
When $\mbox{rk} M =1$ the germ $S$ defined by 
a reduced quasi-ordinary polynomial $f \in \C \{ \r^\vee \cap M \} [Y]$ 
at the origin is just a 
germ of complex plane curve. 
It is well known that  the intersection multiplicities of the different branches of the curve 
at the origin and 
the semigroups associated to each of them are  a complete invariant of the 
{\it embedded topological type} of the plane curve germ $(S, 0)$
(see \cite{Reeve}). 
Eggers shows that this information can be encoded 
by structuring in a labeled tree 
the characteristic exponents of each irreducible factor and 
the orders of coincidence between any two of them
(see \cite{Eggers}).
Wall (see \cite{Wall}) gives a different definition of Egger's tree to 
give a new proof of theorem of  Garc\' \i a Barroso in \cite{Gtesis} 
on the structure of polar curves (see \cite{GB}).
Wall's definition encodes the same amount of information as Egger's definition does and 
involves the use of a simplicial $1$-chain on the tree which is defined from the sequence of 
characteristic integers of the irreducible factors (see definition \ref{ni}).
In the case of a  classical quasi-ordinary hypersurface,
Zariski's result  stated in lemma \ref{rorder} can be reformulated as follows:
If $f =0$ defines a  classical quasi-ordinary hypersurface 
and if $f_1$ is an irreducible factor of $f$
the set of characteristic exponents of $f_1$ union 
the set of orders of coincidence of $f_1$ with the factors of $f$ 
is totally ordered with respect to  the partial order 
defined by the divisibility of the corresponding 
monomials.
Zariski's observation and the sequences of characteristic integers are exactly what 
is necessary to extend Wall's definition to the quasi-ordinary case
in terms of a fixed quasi-ordinary projection $(X, Y) \mapsto X$. 
This is done more generally by Popescu-Pampu (see \cite{PP2}) 
for a  {\it Laurent quasi-ordinary polynomial} $f$,
obtaining  a generalization of  Wall's proof  which describes 
the structure of $\frac{\partial f}{\partial Y} $ 
in terms of the tree of $f$ when $\frac{\partial f}{\partial Y} $ 
is quasi-ordinary.

The definition of the tree in our case runs as follows:
Let $f= f^{(1)} \cdots f^{(r)}$ be the factorization in monic irreducible polynomials of $f$. 
Each factor $f^{(i)} $ of $f$ is quasi-ordinary and 
the subset $ \theta (f^{(i)}) ^{(0)} $ of  $\r^\vee \cap M_g \cup \{  + \infty \}$
whose elements are $0, + \infty   $,  the characteristic exponents
$\l_1 ^{(i)} <_\r \cdots <_\r \l_{g(i)} ^{(i)}$ of $f^{(i)}$ (if they exist)   and the orders of 
coincidence of $f^{(i)}$ with the irreducible factors of $f$ is totally ordered by lemma \ref{rorder};
we denote by  $n^{(i)}_k$ and $e^{(i)}_k$ for $k = 1, \dots, g(i)$,  the sequences of integers 
associated to $f^{(i)}$ by definition \ref{ni} for $i=1, \dots, r$. 

The {\it elementary branch} $ \theta (f^{(i)})$  associated to $f^{(i)}$
is the abstract simplicial complex of dimension one 
with vertices running through 
the elements of the totally ordered subset 
$ \theta (f^{(i)}) ^{(0)} $ and edges running through the segments joining 
consecutive vertices. The underlying topological space is homeomorphic to the segment $[0, + \infty ]$.
We denote the vertex of  $\theta (f^{(i)})$ corresponding to   $\l \in \theta (f^{(i)}) ^{(0)}$ 
by $P^{(i)}_\l $. The simplicial complex $ \theta (f)$  obtained from the disjoint union 
$\bigsqcup_{i=1}^{s} \theta (f^{(i)})$
by identifying  in $ \theta (f^{(i)})$ and $ \theta (f^{(j)})$ the sub-simplicial complex 
corresponding to $ \overline{ P^{(i)}_{0} P^{(i)}_{\l_{(i, j)}} }$ and $ \overline{ P^{(j)}_{0} P^{(j)}_{\l_{(i, j)}} }$ for $1 \leq i <j \leq r$ is a tree. 
We give to a vertex 
$P_\l^{(i)} $ of $ \theta_f (f)$ the valuation  $\l $.
This defines a $0$-chain $ C_0 (f) $ on $ \theta_f (f)$ 
which attaches the value $\l$ to  each vertex $P_\l^{(i)}$ in the 
Eggers tree (counting each vertex only once). 

For $i = 1, \dots, r$ we
define an integral $1$-chain
whose segments are obtained by   
by subdividing the segments of the chain 
\begin{equation} \label{segment}
\overline{ P_0^{(i)} P_{\l_1^{(i)}} } + 
 n^{(i)}_1 \overline{ P^{(i)}_{\l_1^{(i)}}   P^{(i)}_{\l_2^{(i)}} }   + \cdots + 
n^{(i)}_1 \cdots n^{(i)}_{g(i)} \overline{ P^{(i)}_{\l_g^{(i)}} P^{(i)}_{ + \infty }}
\end{equation}
with the points  corresponding to the orders of coincidence of $f^{(i)}$,
the coefficient of an oriented  segment in the subdivision 
is the same as the coefficient of the oriented segment of (\ref{segment}) containing it.
It follows from definition \ref{ni} that these
$1$-chains  paste on   $ \theta_f (f)$ and  define a $1$-chain $ C_1 (f) $
with coefficients in $\Z$.

\begin{Defi}
The Eggers-Wall tree  is the simplicial complex $ \theta (f)$
with the chains  $ C_1 (f)$ and  $ C_0 (f) $. 
\end{Defi}
The chains $ C_1 (f)$ and  $ C_0 (f) $ determine the number of factors of $f$,  
the characteristic exponents of each factor and the orders of coincidence.
The vertex $P_\l ^{(i)}$, 
if  $\l \ne    0, +\infty $ is  not a characteristic exponents of $f^{(i)}$
if and only if the coefficients of the two segments of $ \theta (f^{(i)})$
containing $P_\l ^{(i)}$ coincide.

\section{Embedded resolution procedure} \label{Oka}

In this section we build 
an  embedded resolution of a reduced 
quasi-ordinary polynomial 
which is a composition of toric morphism determined by the characteristic monomials.

\subsection{Definition of good coordinates}

We introduce the notion of $Y$ being a  {\it good coordinate} in 
terms of the coincidence of the parame\-tri\-za\-tions of $f$.
In the following section we build the toric morphisms of the resolution 
using this notion. Different choices of good coordinates provide
isomorphic morphisms.

We keep the notations of section \ref{Eggers}.
We suppose that $f$ is a quasi-ordinary polynomial with $r$ irreducible factors $f^{(1)}, \dots, f^{(r)}$.
Define  ${\cal A}(i):=    ( M \cap \{ \l_{(i, j)} \} )_{j} \cup  \{ \l_1^{(i)} \} $ for $     1 \leq i \leq r$.  
By lemma \ref{rorder}, if the set  ${\cal A}(i)$  is non empty it is totally ordered by $<_\r$.

Then we can define:
\begin{equation} \label{kappa}
\l_{\kappa (i)} := \left\{ \begin{array}{l}
 \min {\cal A}(i) \mbox{ if }
{\cal A}(i) \ne \emptyset
\\
{+ \infty} \mbox{ otherwise }
\end{array}\right\} \mbox{ for } i=1, \dots, r.
\end{equation}

\begin{Lem} \label{Prepa} $\,$
\begin{enumerate}
\item
If $\l_{\kappa (i)}\nleq \l \notin M$ the term $X^{\l}$ does not appear in the 
expansions of the roots of $f^{(i  )}$. 
In particular if  $X^{\l}$ appears in the expansions 
of the roots of $f^{(j)}$ then ${\l}$ is $\geq \l_{(i, j)}$
and the  equality ${\l} = \l_{(i, j)}$implies that $ \l_{(i, j)} = \l_1^{(j)}$.
\item
The case $\l_{\kappa (i)} = + \infty$ happens if and only if $f^{(i)}$ is the only factor of $f$ without characteristic
exponents and $\l_{(i, j)} = \l_1^{(j)}$ for all $j \ne i$. 
\item
If $\l_{\kappa (i_0)} \in M$ then $\l_{\kappa (i_0)} $ is $\geq \l_{\kappa (j)}$ for all $j \ne i_0$.
\item
The set $\{ \l_{\kappa (1)}, \dots,  \l_{\kappa (r)} \}$ is totally ordered by $<$.
\end{enumerate}
\end{Lem}
{\em Proof.}
If $f^{(i)}$ has no characteristic exponent the terms in the expansion of its root have
exponents in $\r^\vee \cap M$. Otherwise, $\l_{\kappa (i)} \nleq \l \notin M$  implies that 
$ \l_1^{(i)} \nleq \l  \notin M$ thus the term $X^\l$ does not appear in the expansion 
of the roots of $f^{(i)}$ by lemma \ref{expo}. 
If $X^{\l}$ appears in the expansion 
of the roots of $f^{(j)}$ then it appears in any  difference of 
roots of  $f^{(i)}$ and $f^{(j)}$ thus $ \l \geq \l_{(i, j)}$.
Moreover, if  $ \l =  \l_{(i, j)}$ then $\l \notin M $ implies that  $\l_{(i,j)} \geq \l_1^{(j)}$  by lemma \ref{expo}.
Since 
$ \l_{\kappa (i)} \nleq \l_{(i, j)}$ 
we have that 
$ \l_{\kappa (i)} \nleq \l_1^{(j)} \notin M$ and therefore $\l_1^{(j)} \geq \l_{(i, j)}$  and  
the equality $ \l_{(i, j)} = \l_1^{(j)}$ follows.

For the second assertion notice that if $f^{(i)}$ and $f^{(j)} $ are two different factors without 
characteristic exponents then $\l_{(i,j)}$ belongs to $M$ thus $ \l_{\kappa (i)}, \l_{\kappa (j)} \ne + \infty$.
If  $ \l_{\kappa (i)} = + \infty$ then
 $\l_{(i,j)}$ is not in $M$ for all $j \ne i$; thus the exponent $\l_{(i,j)}$ appears on a term of 
the parametrization of $f^{(j)} $ and therefore we have $\l_{(i,j)} \geq \l_1^{(j)}$ by lemma \ref{expo}. 
The first assertion for $\l =\l_1^{(j)}$ implies that $\l_{(i,j)} \leq \l_1^{(j)}$ and equality follows.

Now suppose that  $ \l_{\kappa (i_0)} \in M$.
If $j \ne i$ the exponents $\l_{\kappa (i_0)}$ and $\l_{(i_0, j)}$ 
are comparable by lemma \ref{rorder}. We distinguish two cases:

\noindent
(a) $\l_{\kappa (i_0)} \leq \l_{(i_0, j)}$.
Notice that assertion {\it 1} implies that if $f^{(j)}$ has some characteristic exponent then 
$\l_1^{(j)} > \l_{\kappa (i_0)}$. If 
$\l_{\kappa (i_0)} < \l_{(i_0, j)}$ there is 
$j \ne l_0 \ne i_0$ such that 
$\l_{\kappa (i_0)} = \l_{(i_0, l_0)}= \min \{ \l_{(i_0, l_0)}, \l_{(i_0, j)} \} = \l_{(j, l_0)}$
by lemma \ref{rval};
hence the exponents  $\l_{(j, l)}$ and $ \l_{\kappa (i_0)}$ are comparable by lemma \ref{rorder}.
If $\l_{\kappa (i_0)} = \l_{(i_0, j)}$ set $l_0 = j$.

If $\l_{(j, l)} < \l_{\kappa (i_0)}$ we deduce from lemma \ref{rval} that:
$$\l_{(j, l)} = \min \{ \l_{(j, l_0)},\l_{(j, l)} \} = \l_{(l, l_0)}   =  \min \{ \l_{(i_0, l_0)},\l_{(l, l_0)}  \} = 
\l_{(i_0, l)}, $$
and $\l_{(j, l)}$ does not belong to $M$ by definition of  $ \l_{\kappa (i_0)}$.
This shows that $ \l_{\kappa (j)} = \l_{\kappa (i_0)}$.

\noindent
(b) $\l_{(i_0, j)} < \l_{\kappa (i_0)}$.
By definition of $\l_{\kappa (i_0)}$ we have that $\l_{(i_0, j)} \notin M$
and then assertion {\it 1} implies that $\l_{(i_0, j)} = \l_1^{(j)}$.
If $\l_{(j, l)} < \l_1^{(j)}$ we deduce using lemma \ref{rval} that 
$\l_{(j, l)} = \min \{ \l_{(j, l)}, \l_{(i_0, j)} \} $ is equal to 
$\l_{(i_0, l)}$ and  $< \l_{\kappa (i_0)}$.
It follows that $\l_{(i_0, l)} \notin M$, thus $\l_{\kappa (j)} = \l_1^{(j)} < \l_{\kappa (i_0)}$.

For the last assertion we only have to prove that if $\l_{\kappa (i)} = \l_1^{(i)}$ and $\l_{\kappa (j)} =  \l_1^{(j)}$
they are comparable by $<$. 
By lemma \ref{rorder},  $\l_{(i, j)}$ is comparable with $\l_1^{(i)}$ and $ \l_1^{(j)}$.
The case $\l_{(i, j)} < \l_1^{(i)},  \l_1^{(j)}$ implies that  $\l_{(i, j)} \in M$   by lemma \ref{expo}, thus
$\l_{\kappa (i)} \leq \l_{(i, j)}$ a contradiction. Therefore
 we can assume  that $ \l_1^{(i)} \leq \l_{(i, j)}$, replacing $i$ by $j$ if necessary. 
It follows from the definition of order of coincidence that if $ \l_1^{(i)} < \l_{(i, j)}$ 
then $\l_1^{(i)} = \l_1^{(j)}$. If $ \l_1^{(i)} = \l_{(i, j)}$ then the result follows from lemma \ref{rorder}.
\hfill $\ {\diamond}$



We relabel  the factors $f^{(i)} $ of $f$ in order to have:
$\l_{{\kappa (1)}} \leq \l_{{\kappa (2)}} 
\leq \cdots \leq   \l_{{\kappa (r)}}$.
If $\l \in \r \cap M$, the monomial   $X^\l$ appears
in all the roots of $f^{(r)}$ with the same coefficient $c_{\l}^{(r)}$.
Then we define: 
$$\f_0 := \displaystyle{\sum_{\l_{\kappa(r)}  \nleq   \l \in \r^\vee \cap M }}
c_{\l}^{(r)} X^\l, $$
\begin{equation} \label{kappa2}
Y'  := \left\{ \begin{array}{l}
Y+ \f_0  \mbox{ if }  \l_{\kappa(r)} \notin M
\\
Y + \f_0 + c X^{\l_{\kappa(r)}}    \mbox{, for } c \in \C^*  \mbox{ generic,  if } \l_{\kappa(r)} \in M.
\end{array}\right.
\end{equation}
Generic here means that if ${\l_{\kappa(r)}} = {\l_{\kappa(l)}} \in M$ then 
$c - c^{(l)}_{\l_{\kappa(l)}} \ne 0$.  

\begin{Lem} \label{ele}
The polynomial 
$Y'$ has order of coincidence equal to $\l_{{\kappa (i)}}$ with 
$f^{(i)}$ for $i=1, \dots, r$.
\end{Lem}
{\em Proof.}
It follows from  lemma \ref{Prepa} 
that if ${\l_{\kappa(i)}} < {\l_{\kappa(r)}}$ then ${\l_{\kappa(i)}}$
is the order of coincidence of $f^{(i)}$ and $f^{(r)}$
(remark that  ${\l_{\kappa(i)}} \notin M$ by assertion {\it 3 } of lemma \ref{Prepa},  thus ${\l_{\kappa(i)}} = \l_1^{(i)}$
is $\geq \l_{(i,r)}$  by assertion {\it 1} of lemma \ref{Prepa};  it follows from this fact that  $\l_{(i,r)} \notin M$ thus 
$\l_1^{(i)}   \leq  \l_{(i,r)}$
by lemma \ref{expo}).
This implies 
that the order of coincidence 
of $Y'$ with $f^{(i)}$ is well defined and equal to ${\l_{\kappa(i)}}$.
The generic choice of $c$ guarantees in the case ${\l_{\kappa(r)}} \in M$
that the order of coincidence of 
$Y'$ with those factors $f^{(i)}$ of $f$ with ${\l_{\kappa(i)}} = {\l_{\kappa(r)}}$
 is ${\l_{\kappa(r)}}$. 
\hfill $\ {\diamond}$

\begin{Defi} \label{BC}
We say that $Y $ is a {\it good coordinate} for the reduced quasi-ordinary polynomial 
$f \in \C \{ \r^\vee \cap M \} [Y]$ if  
the order of coincidence of $Y$ with 
$f^{(i)}$ is well defined and equal to $\l_{\kappa (i)}$,  for $i=1, \dots, r$.
\end{Defi}

If $Y$ is not a good coordinate for $f$ then 
the  $ \C \{ \r^\vee \cap M \}$-automorphism of $ \C \{ \r^\vee \cap M \} [Y]$ 
that maps  $ Y \mapsto Y'$, for $Y'$ defined in the lemma \ref{ele},
transforms $f \mapsto f' \in  \C \{ \r^\vee \cap M \}[Y']$. 
The polynomial $f'$ is quasi-ordinary, $f'$ and  $f$ 
have the same Eggers-Wall tree and $Y'$ is a good coordinate for $f'$.

\subsection{The first toric morphism of the embedded resolution}

We build the first toric morphism of the embedded resolution and 
we prove that it simplifies the singularity preserving at the same time 
the quasi-ordinary structure.

\subsubsection{The case of a Newton polyhedron with only one compact edge}
We deal first with the  case when all the irreducible factors
of $f$ are parametrized
by series of the form $X^{\l} \varepsilon $  with $\varepsilon (\o_\r) = c $.

We denote by $M_{\l}$ the lattice $M + \l \Z$ for $\l \in \frac{1}{n} M $ 
(resp. $N_\l$ for the dual lattice), 
by $M_{\l}' $ the lattice $M_{\l} \oplus y \Z$ (resp. $N_\l'$ for the dual lattice)
and by $n_\l$ the integer $| M_\l / M | $.
Let $\Sigma$ be a subdivision of $\varrho$ containing  cone $\s:= \varrho \cap \ell$
where $\ell $ is the linear subspace of $N'_\R$ orthogonal to the compact face $[n \l, n y]$ of 
the polyhedron ${\cal N} (f)$ (where $n =\deg f$).
The  subdivision $\Sigma$  of $\varrho$ is rational for the lattices $N_g'$ and $N_{\l}$.
We have the following commutative diagram of equivariant maps:
\begin{equation} \label{fracs} 
\begin{array}{ccc}
Z_{\Sigma, N_{\l}'} &  \stackrel{\Pi_\Sigma}{\rightarrow}   &  Z_{\varrho, N_{\l  }'}
\\
\downarrow  &  & \downarrow 
\\
 Z_{\Sigma, N'} &  \stackrel{\pi_\Sigma}{\rightarrow}  &  Z_{\varrho, N'} 
\end{array}
\end{equation}
where the vertical arrows are defined by lattice extension and the horizontal arrows
are defined by the subdivision $\Sigma$. 
Often we do not  precise  the lattice  if it is corresponds to the below line of the diagram \ref{fracs}.

\begin{Lem} \label{prev1}
The lattice  homomorphism 
$\varphi:  M ' \rightarrow M_\l$ that maps $y \mapsto \l$ and fixes $u \in M$
induces an isomorphism:
\begin{equation} 
M_\s \cong   M_\l \label{prev12}
\end{equation}
If we choose an splitting $ M ' \cong M_\s \oplus \mbox{Ker} ( \varphi )$
we have a semigroup isomorphism:
\begin{equation}
\s^\vee \cap M'  \cong   n_{\l} (y - \l_1 ) \Z \oplus   (\r^\vee \cap M_{\l})   \label{prev11}
\end{equation} 
which corresponds to an isomorphism $Z_{\s, N'} \cong \O_{\s, N'} \times Z_{\r, N_\l}$.
%
\end{Lem}
{\em Proof.} We use the combinatorial arguments in the proofs of lemma \ref{pseudo}
and proposition \ref{binomial} to prove (\ref{prev12}) 
using that $ \s^\bot = \mbox{Ker} (\varphi_\R)$ by  (\ref{orthoduali});
then (\ref{prev11})  holds by (\ref{saturation}). 
\hfill $\ {\diamond}$


We denote by $S_{\Sigma}^{(i)}$ the strict transform of the germ  
$S^{(i)}$ defined by the irreducible factor $f^{(i)}$ of $f$  for $i=1, \dots, r$.
\begin{Lem} \label{prev2}
The intersection $S_{\Sigma}^{(i)} \cap \p_\Sigma^{-1} (\o_\r)$ is the  point 
$\o_1 = (c',  \o_\r) \in \O_\s$ counted  $e_\l^{(i)} :=  ( \deg f^{(i)} ) /n_\l $
 times, where $c' = c^{n_\l}$ and   $c$ is 
the coefficient of $X^{\l_1}$ in any root of the polynomial $f^{(i)}$ defining $S^{(i)}$. 
In particular, the intersection  $S_{\Sigma}^{(i)} \cap \O_\s$ is transversal if and only if $e_\l^{(i)} = 1$. 
The strict transform $S_{\Sigma}$ of $S $ is a germ at the point $ \o_1$.
\end{Lem}
{\em Proof.} 
To simplify the proof we drop the super-index ${(i)}$.
If $\t \in \Sigma$ with $\stackrel{\circ}{\t} \subset \stackrel{\circ}{\varrho}$ 
then $S_\Sigma \cap \O_\t \ne \emptyset $ implies that $\t = \s$ 
since  the face of ${\cal N}_\varrho (f) $ defined by 
$\s$ is of dimension $\geq 1 $  (by lemma \ref{nondeg}).
The strict transform $S_{\Sigma}$ 
is defined on $Z_\s$       by $X^{- n \l }f = 0$ 
and it follows 
that the ideal of $\O_\s \cap S_\Sigma $ is generated by 
$(X^{n_\l (y-\l)} -    c ^{n_\l})^{e_\l}$       
where $c$ is the coefficient of $X^{\l}$ in any root of $f$.
This implies that the intersection of the strict transform $S_\Sigma$ with $\p_\Sigma ^{-1} ( \o_\varrho)$ 
is reduced to the point $o_1 = (c' , \o_{\r}) $ counted         
$e_\l $  times. In particular, the intersection is transversal if and only if $e_\l =1$.
This shows also that the strict transform $S_\Sigma$  is a germ at the point  $o_1$
since this is the only point of intersection with the exceptional fiber.
\hfill $\ {\diamond}$


\begin{Pro} \label{local}
The restriction of the projection $\O_\s \times Z_{\r, N_\l} \cong Z_{\s, N'} \rightarrow Z_{\r, N_\l}$
to $(S_\Sigma, \o_1)$  is quasi-ordinary.
The germ $(S_\Sigma, \o_1)$ is defined by a quasi-ordinary polynomial 
$f_\Sigma \in \C \{ \r^\vee \cap M_\l \} [W] $ (where $W = Y^{n_\l} X^{-n_\l \l} - c^{n_\l}$)
with characteristic exponents $\l' - \l$ for those characteristic exponents $\l' > \l$ of $f$.
If $\l_{(i,j)} $ is the order of coincidence between the irreducible components 
$f^{(i)}$ and $f^{(j)}$ of $f$ then the order of coincidence of 
$f^{(i)}_\Sigma $ and $f^{(j)}_\Sigma$ is $\l_{(i,j)} - \l $.
If $(S, o)$ is irreducible the same holds for $(S_\Sigma, \o_1)$.
\end{Pro}
{\em Proof.} 
We deal first with the case $\l \in M$, i.e.,  $n_\l =1$. 
By lemma \ref{prev2}  the chart $Z_\s$ contains 
the strict transform $S_\Sigma$.
By hypothesis 
the roots $\z^{(i)}$ of $f$ are of the form 
$\z^{(i)} =  c X^\l + \sum_{\l' > \l} c_{\l'}^{(i)} X^{\l'}$, i.e., the coefficient
of the monomial $X^\l$  is the same for all of them.
By lemma \ref{prev2} the strict transform of $Y - \z^{(i)} = 0$  
by the morphism $Z_{\s, N_g '} \rightarrow Z_{\varrho, N_g'}$ 
is defined by:
\begin{equation} \label{van}
0 = X^{ y -\l} - c + \sum_{\l' > \l} c_{\l'}^{(i)}  X^{\l' - \l}
\end{equation}
where the terms $ X^{\l' - \l}$ vanish on the orbit $\O_\s$ for all ${\l' > \l}$.
By lemma \ref{prev1} the chart 
$Z_{\s, N_g'}  $ (resp.$Z_{\s, N }  $)   is isomorphic to 
$\O_{\s, N_g'} \times Z_{\r, N_g}  $ (resp. to  $\O_{\s, N'} \times Z_{\r, N}  $).
Since $n_\l = 1 $ the toric morphism $Z_{\s, N_g'} \rightarrow Z_{\s, N' }  $
restricts to an isomorphism of the orbits $\O_{\s, N_g'} \cong \O_{\s, N'} = \O_\s$
by (\ref{prev11}),  the  coordinate ring of the orbit $\O_\s$ being  equal to $\C[YX^{-\l}]$.
We study the strict transform of $Y - \z^{(i)} = 0$ (resp. of $S$)
at the point of intersection with the orbit $\O_\s$ 
by replacing   the invertible term $X^{ y -\l}$ by
the unit  $c +W$ on (\ref{van})
 (resp. on $X^{- n \l } f = 0$).
We obtain a polynomial $f_\Sigma \in \C \{ \r^\vee \cap M \}[W]$ 
from $X^{- n \l } f$ which splits  in $ \C \{ \r^\vee \cap M_g \}[W]$:
$f_\Sigma = \prod_{i=1}^{n} ( W - \t^{(i)}) $; 
where $\t^{(i)}  = \sum_{\l' > \l} c_{\l'}^{(i)}  X^{\l' - \l}$.
It follows from lemma \ref{expo} that the series  $\t^{(i)}$ 
are quasi-ordinary branches and that their characteristic exponents are obtained 
from those of $\z^{(i)}$ by subtracting $\l$. If $f$ is irreducible 
the same thing happens for $f_\Sigma$. 
Otherwise, we have 
$\t^{(i)} - \t^{(j)} = X^{ - \l} (\z^{(i)} - \z^{(j)}) $ 
and this implies the assertion about the orders of coincidence.

If $n_\l > 1$  we reduce to the previous case by passing through the diagram (\ref{fracs}):

Each irreducible factor of $f$ splits  into $n_\l$ irreducible factors
 in $\C \{ \r^\vee \cap M_\l \} [Y]$ 
having order of coincidence equal to $\l$.
We factor  $f $ as a product $ F_1 \cdots F_{n_\l}$  in $\C \{ \r^\vee \cap M_\l \} [Y]$, 
the $F_i$ being defined by the properties: the order of coincidence of $F_i \ne F_j$ (resp. of any two factors of $F_i$) is $=\l$ 
(resp. is $> \l$). 
The Eggers-Wall tree of $F_i$ is obtained from the Eggers- Wall
tree of $f$ by deleting the vertex $P_\l$ and dividing by $n_\l$ the coefficients 
of the chain $C_1(f)$ between   $P_\l$ and the extreme points $P^{(j)}_{+\infty}$ of the tree
(this follows from lemma \ref{mfield} and definition \ref{ni}).
Then the strict transforms   of $F_i = 0$ by $\Pi_\Sigma$ are disjoint germs
at the $n_\l$ points of intersection with $\O_{\s, N_\l '}$ by lemma  \ref{prev2}.

By lemma \ref{prev1} the toric morphism $Z_{\s, N_\l '} \rightarrow Z_{\s, N}$ 
 corresponds to the 
semigroup inclusion 
$$ n_\l (y -\l) \Z \oplus (\r^\vee \cap M_\l ) \rightarrow (y -\l) \Z \oplus (\r^\vee \cap M_\l ).$$
This map is an unramified covering of degree $n_\l$ 
and it commutes with the projections onto the factor $Z_{\r, N_\l}$ of $Z_{\s, N_\l '}$ and $Z_{\s, N}$.
This provides an isomorphism between the strict transform of $F_i$ by $\Pi_\Sigma$
and $S_\Sigma$ which commutes with the projection onto  factor $Z_{\r^, N_\l}$ for $i =1, \dots, n_\l$.
A fortiori the restriction of the projection $Z_{\s, N'} \rightarrow Z_{\r, N_\l}$ to $S_\Sigma$ is 
quasi-ordinary and the result follows.
\hfill $\ {\diamond}$

With the same hypothesis of proposition \ref{local} we have:
\begin{Cor} \label{reslem}
If $(S, \o)$ is analytically irreducible and if $\l=\l_1$ is the only characteristic exponent of
$\z$ the strict transform $S_\Sigma$ of $S$ is isomorphic to the germ 
$Z_{\r, N_1}$ and the restriction of $\p_\Sigma$ to $S_\Sigma \rightarrow S$ 
is the normalization map. The morphism $\p_\Sigma$ is a partial embedded resolution of 
$S \subset Z_\varrho$.
If $\Sigma'$ is a resolution of the fan $\Sigma$ the map $\p_{\Sigma'} \circ \p_{\Sigma}$ 
is an embedded resolution of $S \subset Z_\varrho$.
\end{Cor}
{\em Proof.} 
It follows from lemma \ref{prev2} that $S_\Sigma$ is isomorphic to 
the germ $(Z_{\r, N_1}, \o_\r)$ and to the  normalization of $(S, \o)$ 
by proposition \ref{closure}. 
We argue as in proposition \ref{binomial} and lemma \ref{pseudo}
to extend the result in this case.
\hfill $\ {\diamond}$

The following remark is a consequence of the proof of proposition \ref{local}.
\begin{Rem} \label{max}
If $f$ is irreducible, $\l = \l_1$  and if  $f_1 \in \C \{ \r^\vee \cap M_\l \} [W]$
defines a good coordinate for $f_{\Sigma}$ then  the image of $f_1=0$ by $\p_{\Sigma}$
is defined by an irreducible quasi-ordinary polynomial  in  $  \C \{ \r^\vee \cap M \} [Y]$ 
with only one characteristic 
exponent $\l_1$ and with maximal order of coincidence with $f$.
\end{Rem}

The following result has been suggested by N\'emethi and McEwan see (\cite{zeta1} and \cite{zeta2}).
\begin{Lem} \label{nemethi}
The morphism $\p_{\Sigma}$ of proposition \ref{local} is an isomorphism 
over $Z_\varrho - S$.
\end{Lem}
The discriminant of the morphism $\Pi_\Sigma$ is described by (\ref{discriminantlocus}).
It follows from this formula that the functions $X^{\l} $ and $Y$ vanishes 
on those orbits of $Z_{\varrho, N_{\l} '}$ 
which are contained in the discriminant locus of  $\Pi_\Sigma$.
The image of these orbits by the map $Z_{\varrho, N_{\l} '} \rightarrow Z_{\varrho}$
is the discriminant of $\p_\Sigma$ and it is contained in $S$
since all the roots of $f$ are of the form $Y = X^\lambda$ up to multiplication  by a unit. \hfill $\ {\diamond}$

\medskip

\subsubsection{The general case}
We build the first toric morphism of the embedded resolution
in the general case.

We  suppose from now on that $Y$ is a good coordinate for $f$ .
The Newton polyhedron of each irreducible factor $f^{(i)}$ with $\l_{\kappa (i)} \ne + \infty$
is an edge with vertices $( \deg f^{(i)}, 0) $ and $(0, \deg f^{(i)} \l_{\kappa (i)}) $ where 
$X^{\l_{\kappa (i)}} $ is the initial monomial of any root of $f^{(i)}$.
Since the set of  $\{\l_{\kappa (i)}  \}$ is completely ordered by $<_\r$ 
the set of compact faces of ${\cal N}_\r (f)$ 
defines a {\it monotone polygonal path}  with
{\it inclinations} running through $ \{ \l_{\kappa (1)},   \dots, \l_{\kappa (r)}  \} - \{ + \infty \}$
independently of the choice of good coordinate
(see \cite{GP} for the terminology).
This fact is a special feature of quasi-ordinary singularities
and it is a generalization of the plane curve case.

The dual fan $\Sigma_1$ of the polyhedron ${\cal N} (f)$ 
is obtained by intersecting $\varrho$ with the linear 
hyperplanes 
$\ell_{ \kappa (j) } : = \langle y - \l_{\kappa (j)} , u \rangle  = 0$
for those $\l_{\kappa (j)} \ne + \infty$. 
Since we have that 
$\{ \l_{\kappa (j)} \}$ is totally ordered by $<_\r$
we find that the cones $\varrho \cap \ell_{ \kappa (j) }$ belong to  $\Sigma_1$
since they cannot intersect in the interior of $\varrho$.
Geometrically, this implies that 
the exceptional locus of $\p_{\Sigma_1}$ 
is a {\it bamboo} of ${ \mathbf P}^1_\C$, each one of
them being the closure of the orbit $\O_{\varrho \cap \ell_{ \kappa (j) }}$
(we say that a curve is a bamboo if the dual graph of intersection of its irreducible  components 
is isomorphic to the subdivision of a segment). 


 

\begin{Pro} \label{EW} If $ \l_{\kappa (i)} \ne + \infty$ then we have:
\begin{enumerate}
\item 
The strict transform of  $S^{(i)}$ by $\p_{\Sigma_1}$ is a germ 
$(S^{(i)}_{\Sigma_1}, \o_1^{(i)} )$ at the point of intersection with the 
exceptional curve $\p_{\Sigma_1} ^{-1} ( \o_\varrho )$.
\item 
The Eggers-Wall tree of a polynomial defining the strict transform 
$S_{\Sigma_1}$ at the point $o_1^{(i)}$  
is obtained from $\theta (f) $ by removing  the segment $ [ P_0^{(j)}, P^{(j)}_{\l_{\kappa (i)}} [$
from 
the sub-tree of $\theta (f)$ given by $\bigcup \theta ( f^{(j)} )$, for 
$f^{(j)}$ with order of coincidence $> \l_{\kappa (i)} $ with $f^{(i)}$.
The coefficients of the vertices of the resulting tree 
are obtained by subtracting  $\l_{\kappa (i)}$.
The coefficients of the associated $1$-chain are obtained 
by dividing by $n_{\l_{\kappa (i)}}$.
\end{enumerate}
\end{Pro}
{\em Proof.} 
The first assertion follows from lemma \ref{prev2}.
It follows from proposition \ref{local}  that 
$o_1^{(i)} = \o_1^{(j)} $ if and only if 
the irreducible factors of the  symbolic restrictions of  $f^{(i)}$ and   $f^{(j)}$
to the compact edges of their Newton polyhedra coincide

This is equivalent
$f^{(i)}$ and $f^{(j)} $ have order of coincidence $> \l_{\kappa (i)} = \l_{\kappa (j)} $.
The second assertion  follows from  proposition \ref{local}
since
the characteristic exponents and the order of coincidence of $f^{(i)}_\Sigma $ and $f^{(j)}_\Sigma$ 
are obtained from  those corresponding  to $f^{(i)}$ and $f^{(j)} $ by subtracting $\l_1$.
The strict transform $S_\Sigma^{(i)}$ is a toric quasi-ordinary hypersurface relative to the base 
$Z_{\r, N_{\l_{\kappa (i)}}}$
by proposition \ref{local}
and the statement about 
the coefficients of the associated $1$-chain follows from this change of lattice 
by lemma \ref{field}.
\hfill $\ {\diamond}$

\begin{Rem}  \label{EW2} 
If $ \l_{\kappa (i)} = + \infty$  the strict transform of $S^{(i)}$
is the germ of the closure of the orbit associated to the edge $y \R_{\geq 0}$ at 
the point of intersection with the 
exceptional curve $\p_{\Sigma_1} ^{-1} ( \o_\varrho )$.
\end{Rem}
The assertion follows 
from the description of the exceptional locus and the discriminant locus of a 
toric modification given in section \ref{toricg} once it is noticed that 
the point of intersection $o_1^{(i)}$ of $S^{(i)}$ with $\p_{\Sigma_1} ^{-1} ( \o_\varrho )$
is the orbit associated to the $(d+1)$-dimensional cone of $\Sigma_1$ which contains the cone $y \R_{\geq 0}$. 


\subsection{The toric embedded resolution}

We show the way to iterate the procedure of the previous section to build 
an embedded resolution of $S \subset Z_\varrho$ by first eliminating the characteristic exponents
and then by resolving the toric singularities of the ambient space.

By proposition \ref{EW}, the germs defined by the strict transform 
at each of the points $o_1^{(i)}$ of intersection with the
exceptional curve are simpler  toric quasi-ordinary hypersurface singularities.
In a finite number of iterations of this procedure
the strict transform becomes  a union of $r$ toric quasi-ordinary hypersurface germs 
with no characteristic exponents at all, i.e., is collection of $r$ germs of affine toric varieties at 
the special points.  It follows from propositions \ref{local} and  \ref{closure} that 
the strict transform of $(S,o)$ is its normalization. 
Thus this method provides an {\it embedded normalization} (in a normal environment) of the germ 
$(S,o) \subset (Z_\varrho, \o_\varrho)$.
We keep the information of the toric singularities of the ambient space by defining at each stage
a toroidal embedding without self-intersection:

First, we associate to the toric  quasi-ordinary hypersurface 
$(S, o) \subset (Z_\varrho, \o_\varrho)$ embedded with a good coordinate
the toroidal 
embedding defined by $(Z_\varrho, N'_0)$. 
Its conic polyhedral complex $\Theta_0$ is equal to 
$(\varrho, N_0)$.
Then, 
we associate to each point of intersection $\o_1^{(i)}$ of the strict transform $S_{\Sigma_1}$ with the
exceptional fiber a normal hypersurface  $S_1 ^{(i)}$ defined by taking a good coordinate
for the quasi-ordinary projection of $S_{\Sigma_1}$ of proposition \ref{local}.
Obviously, if $\o_1^{(i)} =\o_1^{(j)}  $ we have $S_1 ^{(i)} = S_1 ^{(j)}$.
(see remark \ref{max}).

\begin{Lem}   \label{glue}
The c.p.c. $\Theta_1$ 
associated to the toroidal embedding
defined by the variety $Z_{\Sigma_1}$ and the set
normal hypersurfaces $\{ (S_1 ^{(i)}, \o_1^{(i)}) \} _{\l_{\kappa(i)} \ne + \infty}  
\cup \{ \overline{ \O_\s } \}_{\s \in \Sigma_1^{(1)}} $
is obtained from 
the c.p.c. $\Sigma_1$ by 
adding for each point in the set $ \{ \o_1^{(i)} \}_{i=1, \dots, r}^{\l_{\kappa(i)} \ne + \infty} $
the c.p.c.   $(\varrho , N_{\l_{ \kappa (i) }} ' )$
and  pasting it to $\Sigma_1$  by identifying   
$(\r \times \{ 0 \} , N_{\l_{ \kappa (i) }} \times \{ 0 \} )$
with $(\varrho \cap \ell_{ \kappa (i) } , N_{\varrho \cap \ell_{ \kappa (i)}} )$
by the lattice isomorphism  corresponding to (\ref{prev12}) by duality.
The c.p.c. $\Theta_1$  is independent of the choice of good coordinates.
\end{Lem}
{\em Proof.} 
To simplify the proof we drop the index $i$, we denote $\l_{\kappa (i)}$ by $\l$  and 
we keep notations of proposition \ref{local} and 
lemmas \ref{prev1} and \ref{prev2}. 
The germ $(S_1, \o_1)$ is defined by the vanishing of a monic polynomial 
$f_1 \in \C \{ \r^\vee \cap M_\l \} [W]$ of degree one  where $W = 
X^{n_\l (y - \l)} - c$.
We deduce from lemma \ref{prev1} that the analytic algebra of the germ 
$(Z_{\Sigma}, \o_1)$ is isomorphic to $\C \{ \varrho^\vee \cap M_\l' \}$
by the isomorphism that maps $f_1 \mapsto X^{y_1}$ and 
$X^u \mapsto X^u $ for all $u \in \r \cap M_\l$.
Since the c.p.c. associated to the torus embedding of $Z_{\varrho, N_1}$ is $({\varrho, N_1})$
the same holds for the toroidal embedding corresponding to $\Sigma_1^{(1)}$
and the set of normal hypersurfaces 
 ${\cal H} =  \{ \overline{ \O_\s } \}_{ \s \in (\varrho \cap \ell ) ^{(1)} } \cup \{ S_1 \}$. 
The sub-c.p.c. associated to the toroidal embedding corresponding to 
${\cal H} - \{ S_1 \}$ is $(\varrho \cap \ell , N_{\varrho \cap \ell})$; it  is isomorphic to
 $(\r, N_1)$,  the pasting isomorphism being obtained from (\ref{prev12}) by duality.
\hfill $\ {\diamond}$

Then we continue as follows:

If the quasi-ordinary polynomial defining the germ  of the strict transform 
$(S_{\Sigma_1}, \o_1^{(i)})$ has some  characteristic exponent
we put it in good coordinates; then its Newton polyhedron
defines a subdivision of $(\varrho,  N'_{ \l_{ \kappa (i) } })$,  for $1 \leq i \leq r$.
These subdivisions glue up to define 
a subdivision 
${\Sigma_2} $ of  the c.p.c. $\Theta_1$
since the pasting cones $(\r \times \{ 0 \} , N_{\l_{ \kappa (i) }} ) $  are not subdivided, for $1\leq i \leq r$.

The corresponding toric modifications, defined locally,  paste 
into a toroidal modification $\p_{2}: Z_{2} \rightarrow Z_{1} $; 
(we denote the variety $Z_{\Sigma_1}$ by $Z_1$,
the morphism $\p_{\Sigma_1} $ by $\p_1$, and $S_{\Sigma_1}$ by $S_1$).
By iterating this procedure we obtain:
A modification $\p_{k}: Z_{k} \rightarrow Z_{{k-1}}$, 
where the variety $ Z_{k} $ is given 
with  the structure of toroidal embedding (${\Sigma_k}$
denoting its associated c.p.c.). 
The strict transform $S_{k}$ of $S$ by $ \p_{k} \circ \dots \circ \p_{1} $
at the points of intersection with the exceptional fiber 
is given with a quasi-ordinary projection 
and the associated Eggers-Wall tree is obtained from the eventually non connected tree of
$S_{{k-1}}$ as indicated by proposition \ref{EW}.
If the quasi-ordinary polynomial defining the germ 
$S_{k} $ at any of these points  has some characteristic monomial 
we define a 
a finer toroidal embedding for $ Z_{k} $ (with c.p.c. $\Theta_k$ defined by using lemma \ref{glue})  
and a  subdivision $\Sigma_{k+1} $ of $\Theta_k$
with associated modification $\p_{{k+1}}: Z_{{k+1}} \rightarrow Z_{{k}}$.
In a finite number $k_0$ of steps 
the quasi-ordinary polynomials defining the germ  $S_{{k_0}}$ at the points 
of intersection with the exceptional fiber
have no characteristic monomials.
Then it follows from corollary \ref{reslem} that:

\begin{The} \label{hres}
The proper morphism  $\p = \p_{{k_0}} \circ \dots \circ \p_{1}$ 
is a partial embedded  resolution of the  quasi-ordinary
hypersurface germ $(S,o) \subset (Z_\r, \o_\r)$. 
The restriction $S' \rightarrow S$  of $\p$ to the strict transform $S'$ of $S$ is the normalization map.
\hfill $\ {\diamond}$
\end{The}
An embedded resolution of $S \subset Z_\varrho$ is obtained 
by composing $\p$ with any toric resolution of 
the toroidal embedding $Z_{k_0}$ with the c.p.c.\ 
$\Sigma_{k_0}$ (or also with the  c.p.c.\
$\Theta_{k_0}$).

\begin{Rem}
The irreducible components of the  exceptional fiber $\p^{-1} (\o_\varrho)$  of the partial resolution 
are projective lines  ${\mathbb P}^1_\C$.
The dual intersection graph of the components  of $\p^{-1} (\o_\varrho)$
is obtained from the Eggers-Wall tree $\theta(f)$ by deleting the 
extremal segments. 
\end{Rem}
One of this segment joins the base vertex $P_0$ to one defined by the first characteristic exponent
of the reduced $f$ and the others corresponds to the segment containing the point 
$P_{+\infty}^{(i)}$ for $i=1, \dots, r$.

\subsection{The case of plane curve germs}

The case of plane curve germs corresponds to  $\mbox{rk } N=1$. We keep the same notations.
The partial resolution procedure depends only on the Eggers-Wall tree 
of $f \in \C \{ X \} [Y]$ with respect to the projection $(X, Y)  \mapsto X$ or more precisely
on the choice of the curve $X= 0$. 
If $f$ is irreducible then our construction is closely related to 
the construction of the  ``{\it Tschirnhausen good resolution tower}'' 
of A'Campo and Oka (see \cite{A'C}, Theorem 4.5).
In particular if the curve $X= 0$  is not contained in the tangent cone  of $S$ 
we  show that this procedure leads to a minimal embedded resolution of the curve.

Let $f \in \C \{ X \} [Y]$ be a reduced polynomial
with $Y$  a good coordinate for $f$.
We keep notations of theorem \ref{hres} and we give some more definitions and notations.
We denote by $\Theta_{k_0}^{reg}$  the minimal regular subdivision  of the   c.p.c. $\Theta_{k_0}$
(for the minimal regular subdivision in the toric case see Proposition 1.19 of  \cite{Oda}).  
This provides a  resolution $p: Z_{\Theta_{k_0}^{reg}} \rightarrow \C^2$
where $p :=  \p  \circ \p_{\Theta_{k_0}^{reg}}$ which is canonically determined 
from the projection $(X, Y) \mapsto X$.

Denote by ${ \cal G} (p, 0)$ (resp. ${ \cal G} (p, f)$)   
the subset of $\Theta_{k_0}^{reg}$ whose elements are the  
cones corresponding to  non empty intersections
of pairs of components of the exceptional divisor of the resolution $p$,
(resp. of the total transform of  $S$ by $p$).
Denote by ${\cal G} (\p, 0 )$ (resp. by ${\cal G} (\p, f )$)  the subset of $\Theta_{k_0}$
of those  cones corresponding to  non empty intersections
of pairs of components of the exceptional divisor of the partial resolution $\p$ 
(resp. of the total transform of  $S$ by $\p$).

Recall that each edge of $\Theta_{k_0}^{reg}$ corresponds to an irreducible 
divisor in the toroidal embedding
and any pair of these  
divisors intersect if and only if the corresponding edges belong to the same cone.
It follows that  ${ \cal G} (p, 0)$  (resp. ${ \cal G} (p, f)$)  
is combinatorially isomorphic to the {\it resolution graph
of the resolution} (resp. to the {\it total resolution graph
of the resolution}), we just drop the dimension of the faces by one.
We deduce from proposition \ref{EW}, remark \ref{EW2} and an easy induction that 
the Eggers-Wall tree $\theta(f)$ is combinatorially isomorphic to ${\cal G} (\p, X f )$.

The valency of a cone  $e$  in a conic polyhedral complex is the number of cones 
of the complex containing $e$ as a facet.
We denote by $\# 1$ the edge of ${ \cal G} (p, 0)$ (resp.  ${ \cal G} (p, f)$)
which corresponds to the first blow up and we define 
\[
\delta (e)  := \left\{
\begin{array}{lcc}
 \mbox{ valency of } e   & \mbox{ if } &    e \ne \# 1 \\
1 + \mbox{ valency of } e   & \mbox{ if } &    e = \# 1, 
\end{array} \right.
\]
The valency of $e$ and the integer $\delta (e)$  depends on the complex containing $e$.
The following lemma implies that the set of non extremal vertices of $\theta(f)$ correspond bijectively 
with the rupture vertices of  ${\cal G} (p, f) $
(which are defined by
those $e$ with      $\delta (e) \geq 3$). 

\begin{Lem} \label{valency}
Let $f \in \C \{ X \} [Y]$ be a reduced polynomial of degree $> 1$, 
such that  $Y$ is a  good coordinate for $f$.
For any edge $e$ in ${ \cal G} (\p, 0)$ we have:
\begin{enumerate}
\item
The integer $\delta(e)$ in ${ \cal G} (p, Xf)$ is $\geq 3$.
\item
If $\l_{\kappa(1)}^{-1} \notin \Z_{> 1}$ then $\delta(e)$ in ${ \cal G} (p, f)$ is $\geq 3$.
\end{enumerate}
\end{Lem}
{\em Proof.}
Recall that we have relabel the factors of $f$ in order to have 
$\l_{\kappa(1)} \leq \dots \leq  \l_{\kappa(r)}$.
We show first the assertion for the exceptional divisors appearing in the first 
toric modification $\p_{\Sigma_1}$.
The extremal edges of the fan $\Sigma_1$, which are 
defined by the vectors $u_1, u_2$ of the canonical basis, 
correspond to the divisors $X= 0$ and $Y=0$ respectively.
If  $\l_{\kappa(j)} \ne + \infty $, there is an exceptional divisor $D_{\l_{\kappa(j)}}$ 
of $\p_{\Sigma_1}$ corresponding 
to $d_{\l_{\kappa(j)}} \in { \cal G} (\p , f)$.
We denote by the same letter the edge $d_{\l_{\kappa(j)}}$ of $\Sigma_1$ and the primitive vector 
$( n_{{\l_{\kappa(j)}}}, n_{{\l_{\kappa(j)}}} {\l_{\kappa(j)}} )$ on this edge for the lattice $N'_0$.
We say that a two dimensional cone $\s$  is on the {\it left} 
(resp.  on the {\it right}) of the vector $  d_{\l_{\kappa(j)}}  \in \s$ 
if $\s \subset \langle  d_{\l_{\kappa(j)}} , u_2 \rangle$
(resp. $\s \subset \langle u_1,  d_{\l_{\kappa(j)}} \rangle$).

By proposition \ref{EW}, the divisor $D_{\l_{\kappa(j)}}$
meets the strict transform of $S$ by $\p_{\Sigma_1}$.

If $ \l_{\kappa(r)}  > 
 \l_{\kappa(j)} $ (resp. if $ \l_{\kappa(j)}  > \l_{\kappa(1)}$)  then 
there exist a two dimensional cone 
on the left (resp. right) of  $d_{\l_{\kappa(j)}}$ in ${ \cal G} (p , f)$,
obtained from the minimal regular subdivision of the cone 
$\s \in { \cal G} (\p , f)$, on the right (resp. on the left) of  $d_{\l_{\kappa(j)}}$.
Therefore if $ \l_{\kappa(r)}  > 
 \l_{\kappa(j)}  > \l_{\kappa(1)}$ we have 
$\delta (d_{\l_{\kappa(j)}} ) \geq 3$.

If  $\l_{\kappa(r)} = + \infty $ then $Y $ divides $f$ and $Y= 0$ is a component of the strict transform of 
$S$ by $\p_{\Sigma_1}$.
If  $\l_{\kappa(r)} \ne  + \infty $ we two cases may occur: 
a) if  the cone  
$\s=\langle   d_{\l_{\kappa(r)}} , u_2 \rangle$ 
is not regular we have a two dimensional cone in  ${ \cal G} (\p , f)$ on the left of $d_{\l_{\kappa(r)}} $; 
b) the cone $\s$ is regular thus  ${\l_{\kappa(r)}} \in M$.
By the proof of  lemma 
\ref{Prepa} there exist $i \ne r$ such that  $ {\l_{\kappa(i)}} = \l_{(i, r)} = {\l_{\kappa(r)}} $.
By proposition \ref{EW} this  implies that 
the strict transforms of $f^{(i)} = 0 $ and $f^{(r)} = 0 $ 
meet the divisor $D_{\l_{\kappa(r)}}$ in two different points so that 
we have $\delta( d_{\l_{\kappa(r)}} ) \geq 3$.

Now we deal with the divisor $D_{\l_{\kappa(1)}}$. 
The cone $\langle  u_1 , d_{\l_{\kappa(1)}} \rangle$ belongs to ${ \cal G} (\p , X f)$
and we deduce from this that $\delta (d_{\l_{\kappa(1)}}) \geq 3$ in ${ \cal G} (\p , X f)$. 
If the cone $\langle  u_1 , d_{\l_{\kappa(1)}} \rangle$ is not regular 
we can argue as before to show the existence of a two dimensional cone 
of ${ \cal G} (\p , f)$
on the right of $ d_{\l_{\kappa(1)}}$. 
Otherwise
we have $n_{\l_{\kappa(1)}} \l_{\kappa(1)} = 1$ and 
if  $\l_{\kappa(1)}^{-1} \notin \Z_{>1}$ the only possibility is  
$d_{\l_{\kappa(1)}} =(1, 1) $. Then we have  
$\l_{\kappa(1)} \in M $ and thus $\l_{\kappa(r)} = \l_{\kappa(1)}$ by lemma \ref{Prepa}.
This case has already been solved.

These facts give the assertion for $e$ corresponding to an exceptional divisor of $\p_{\Sigma_1}$. 
When we iterate, the curve  $X=0$ corresponds to the equation of the exceptional divisor 
meeting the strict transform thus  after the first step we are always in the case {\it 1} 
and  proposition follows.
\hfill $\ {\diamond}$

An exceptional divisor $D$ of the resolution $p$ is {\it collapsible} if
it has {\it self-intersection number} equal to $-1$ and 
the corresponding edge $d \in { \cal G} (p, 0)$ has 
 $\delta(d) \leq 2$ in ${ \cal G} (p, f)$.
If the divisor $D$ is collapsible, the modification obtained by 
blowing down $D$ is still a resolution and the corresponding resolution graph 
is obtained from ${ \cal G} (p, 0)$ by deleting the point corresponding to $D$.
The self intersection of the divisors which are  images of compact divisors meeting $D$ 
is increased by one. In a finite number of steps we obtain a {\it minimal resolution},
i.e., a resolution in which no  exceptional divisor is collapsible. 
The minimal resolution is unique up to isomorphism (see \cite{Lau}).

\begin{Cor} 
If $\l_{\kappa(1)}^{-1} \notin \Z_{> 1}$, in particular if the
projection $(X, Y) \mapsto X$ is transversal for all the components of $f$ 
then the morphism $p$ is the  minimal resolution. 
\end{Cor} 
{\em Proof.}
The self intersection numbers of the
exceptional divisors of the minimal resolution of a toric surface 
singularity are $\leq -2$ (see proposition 1.19 of \cite{Oda}).
This implies that the exceptional divisors 
corresponding to edges in  ${ \cal G} (p, f) - { \cal G} (\p , f)$
are not collapsible.
Then the corollary follows from lemma \ref{valency}.
\hfill $\ {\diamond}$

\begin{Rem}

The number of {\it local toroidal morphisms} used in the partial resolution $\p$  is not necessarily 
equal to the complexity of the resolution (as defined by \cite{Le}).

\end{Rem}

For instance $f = (( Y - X)^2 - X^3)) (( Y + X)^2 - X^5))$
has characteristic exponents $\{ 1, \frac{3}{2}, \frac{5}{2} \}$.
The projection $(X, Y) \mapsto X$ is transversal for the two irreducible components.
It follows easily that the number of local toroidal morphisms 
used to define the partial resolution in this coordinates is three.
On the other hand, the resolution graph is a bamboo so that the resolution complexity is equal to one.

\section{The semigroup associated to a toric quasi-ordinary branch}

We associate to the quasi-ordinary branch $\z$  a semigroup $\Gamma$ which is 
determined from the characteristic exponents; the construction of $\Gamma$  involves also
a generalization of the notion of the plane curves with maximal contact with a given branch given by 
Lejeune
\cite{Lejeune} and this relation can be described by using 
the {\it approximate roots} of the polynomial $f$.
The main part of the results and the proofs of this section is given     in \cite{semi}.

\subsection{Definition of the semigroup}
In the following sections we study  a fixed toric quasi-ordinary singularity $S$ parametrized 
by a toric quasi-ordinary branch $\z \in \C \{ \r^\vee \cap \frac{1}{n} M \}$ 
with $g \geq 1 $ characteristic exponents $\{ \l_1, \dots, \l_g \}$ and 
with minimal polynomial $f \in \C \{ \r^\vee \cap M \} [Y]$.
If $\mbox{rk } M =1$ then the singularity $S$ is a plane branch and 
the set of intersection multiplicities $(S, S')_0$ of $S$ which those plane curve germs $S'$ not
containing $S$ as a component forms a sub-semigroup of $(\Z_{\geq 0}, +)$
which is an invariant of the germ $S$
and which is generated by the following elements (see \cite{Mo}):
\begin{equation} \label{recurrencia}
\bar{\g}_{1} = {n} \l_1, \quad
\bar{\g}_{j+1} = n_j \bar{\g}_{j} + {n} \l_{j+1} - {n} \l_{j}, \; \makebox{ for} \; j= 1, \dots, g-1. 
\end{equation}
For $j= 0, \dots, g-1$,  we expand: 
\begin{equation} \label{recurrencia2}
 \begin{array}{l}
\bar{\g}_{j+1}   = 
{n} \left( (n_1 -1) n_2 \cdots n_j  \l_1 + (n_2 -1) n_3 \cdots n_j \l_2 + \cdots +
(n_j -1)  \l_j  + \l_{j+1} \right) 
\\
  \quad  \quad \stackrel{\makebox{Def. \ref{ni}}}{=}
n_1 \dots n_j \left( (e_0 - e_1) \l_1 + (e_1 -e_2)
\l_2 + \dots + (e_{j-1} -e_{j}) \l_j + e_{j} \l_{j+1} \right) 
\end{array}
\end{equation}
We denote $\frac{1}{n} \bar{\g}_i $ by $\g_i$ for $i=1, \dots, g$ and we have:
\begin{equation} \label{zar}
{\g}_{1} =  \l_1, 
{\g}_{j+1} = n_j {\g}_{j} +  \l_{j+1} -  \l_{j}, \; \makebox{ for} \; j= 1, \dots, g-1. 
\end{equation}
\begin{Defi} We associate to the  quasi-ordinary branch $\z$
the sequence of semigroups $\Gamma_j =  \r^\vee \cap M +\g_1 \Z_{\geq 0} +  \cdots + \g_g \Z_{\geq 0}   $ 
for $j = 0, \dots, g$. 
\end{Defi}
We denote $\Gamma_g$ by $\Gamma$ and $n \Gamma_j $ by 
$ \bar{\Gamma}_j $ for $j = 0, \dots, g$. 
The classical semigroup of a plane branch is $ \bar{\Gamma}_g$.

If $\z$ is a {\it classical} quasi-ordinary branch suitably 
{\it normalized}\footnote{In the case of a plane branch
this condition means that  $X = 0$ is not contained in its the tangent 
cone}. Lipman proved that 
the sequence of  characteristic exponents is
an analytical invariant of the germ  it parametrizes
when $\dim S =2$, by building  a (non embedded) resolution of the germ
(see \cite{Lipman0}, \cite{Lipman1}) which determines the characteristic exponents.
Luengo gives another proof  also using resolutions (see \cite{Luengo}).
If the germ is analytically irreducible 
the characteristic exponents define a complete invariant of the embedded topological 
type of the hypersurface $S \subset \C^{d+1}$ it parametrizes 
(see \cite{Gau} and \cite{Lipman2}).
We proved in  \cite{Tesis}
that if $\t$  and $\z$ are quasi-ordinary branches parametrizing $S$ then 
the semigroups associated to them are isomorphic and moreover that 
the minimal set of generators of this semigroup defines the sequence of 
characteristic exponents of any normalized
quasi-ordinary branch parametrizing $S$. By Gau's characterization 
it follows that 
the semigroup 
$\Gamma$ defined above 
is a complete topological invariant of the embedded topological type of germ $(S, 0)$.

The following lemma generalizes the properties of the semigroups of 
plane branches (see  \cite{Mono}, Chapitre I, Lemma 2.2.1) to the quasi-ordinary hypersurface case 
(see \cite{Tesis}).

\begin{Lem}  (See \cite{semi})  \label{Azevedo} 
\begin{enumerate}
\item 
  The sub-lattice  of $M$ generated by $ \Gamma_j$
is equal to  $M_j$, for $0 \leq j \leq g$.
\item
  The order of the image of ${\g}_{j}$ in the group  $M_j / M_{j-1}$  is equal to $n_j$ for 
$j=1, \dots, g$.
\item 
  We have that  ${\g}_{j} >  n_{j-1} {\g}_{j-1}$ for $j= 2, \dots, g$.
 \item
 If a vector $u_j \in \r^\vee \cap M_j$  then 
 we have $u_j + n_j \g_{j} \in \Gamma_j $.
\item
  The vector $n_j {\g}_{j}$ belongs to the semigroup $ \Gamma_{j-1}$
for $j=1, \dots, g$, moreover we have a unique relation:
\begin{equation} \label{canrel}
n_j {\g}_{j} = \a^{(j)} + l_1^{(j)}  {\g}_{1} + \cdots + l_{j-1}^{(j)}  {\g}_{j-1}
\end{equation}
such that $ 0 \leq  l_{i}^{(j)} \leq n_i -1 $ and  $ \a^{(j)} \in M_0$, for $j = 1, \dots, g$.
\end{enumerate}  
\end{Lem}

In the plane branch case 
several authors have studied the
properties of those curves 
$S'$ such that the 
intersection multiplicity with $S$ at the origin belongs 
to the unique minimal set of generators of the semigroup of the branch
(see \cite{Mo}).  Lejeune introduced  
the notion of curves of maximal contact with a given plane curve germ 
for curves defined over a field of  arbitrary characteristic 
in terms of the resolution (see \cite{Lejeune}).
If the characteristic is zero it turns out that both notions are equivalent
(see \cite{Campillo}).  If  the projection $(X, Y)$ is transversal 
we can study these curves by means of the minimal polynomials 
of suitable  truncations of the roots of $f$. 
When  we do this with respect to an arbitrary projection, 
the curves we obtain provide a non necessarily minimal set of 
generators of the semigroup of the branch $S$.
These curves can be represented  by some of the {\it approximate roots} of 
the polynomial $f$ (see \cite{Moh}) and we call them {\it semi-roots}. 
following the terminology of
\cite{A3}. See Popescu-Pampu's survey  \cite{Popescu} for more 
on the notion of semi-root. 

\begin{Defi} A $j^{th}$-semi-root of $f$ is 
an irreducible quasi-ordinary polynomial  in $\C \{ \r^\vee \cap M \}[Y]$  
of degree $n_0 \dots n_j$ which has order of coincidence equal to  
$\l_{j+1}$ with $f$, for $j=0, \dots, g$.
\end{Defi}
The minimal polynomials of the quasi-ordinary branches 
$p_0 + \dots + p_j$ obtained by truncating $\z$ in remark \ref{r} are 
$j^{th}$-semi-roots of $f$ for $j= 0, \dots, g$.
\begin{Pro} \label{semigroup} (see \cite{Tesis} and \cite{semi}) 
Let $q \in \C \{ \r^\vee \cap M \} [Y]$ a monic polynomial of degree $n_0 \dots n_j$.
Then $q$ is a $j$-semi-root of $f$ if and only if $q(\z) = X^{{\g}_{j+1}} \varepsilon_{j}$
for a unit $\varepsilon_{j}$.
\end{Pro}

The notion of semi-root extends  the properties of maximal contact 
with respect to the resolution to the quasi-ordinary case
(see the proof of theorem \ref{hres} and remark \ref{max}).
\begin{Rem}
The polynomial 
$q_j$ is a $j$-semi-root of $f$ is and only if 
the strict transform of $q_j =0$ by the morphism $\p_j \circ \dots \circ \p_1$ is a germ  
defined by  a good coordinate and conversely.
\end{Rem}
This  follows from      the proof of theorem \ref{hres} and remark \ref{max}.

Let $A$ a ring containing $\Q$ as a subring.
{\it Approximate roots \index{approximate roots}}  are defined by Abhyankar and Moh, 
(see \cite{Moh}, \cite{Ploski}, and  \cite{Popescu}).
If $p$ is any monic polynomial and $k$ divides the degree of $p$
there is a unique monic polynomial $r$ in $A[Y]$ of degree $\frac{\deg p}{k}$ such that 
$\deg (p- r^k) < \deg p -\frac{\deg p}{k}$. We say that $r$ is a $k$-semi-root 
of $p$. 
We can use  proposition  \ref{semigroup} to prove that  the $e_j$-approximate roots of a 
quasi-ordinary
polynomial $f$ are semi-roots, and therefore
are irreducible
quasi-ordinary polynomials with a prescribed order of coincidence with
the polynomial $f$ (see \cite{Tesis} and \cite{semi}).

\subsection{Expansion in terms of semi-roots}

The expansions  in terms of semi-roots are introduce  by 
Abhyankar in the plane curve case (see \cite{A3})
and used by 
Popescu-Pampu in the case 
of a quasi-ordinary hypersurface singularity (see \cite{PP2}).

We fix  from now on  a {\it complete set $q_0, \dots, q_{g}$ of semi-roots of $f$} ($\deg q_i = n_0 \cdots n_i$ 
for $i=0, \dots, g$).
We assume that 
the coefficient of 
the term  $ X^{\g_{j+1}  }$ appearing in ${q}_j (\z) $ by  proposition \ref{semigroup} 
is equal to one for $j=0, \dots, g-1$ in order to simplify some 
computations.

We recall now the classical $q$-adic expansion 
of a polynomial $p_0 \in A [Y]$ with coefficients on a domain 
$A$ in terms 
of a polynomial $q \in A[Y] $ having invertible  leading term (see \cite{Mo}). 
The sequence of 
Euclidean divisions:
$$ p_0 = p_1 q + a_0,\;\;  p_1 = p_2 q + a_1,\;\; \dots, \;\;  p_s = p_{s+1} q + a_s,$$
(where $s$ is the first integer for which 
$p_{s+1} = 0$)
provides a unique decomposition of the form 
$$ H = a_0 + a_1 q + a_2 q^2 + \cdots + a_s q^s_s \mbox{, for } 0 \leq \deg a_i \leq \deg q -1 .$$

\begin{Lem} \label{abhy} (see \cite{PP2})
Any polynomial 
$h \in \C \{ \r^\vee \cap M \} [Y]$ can be written in a unique way as 
\begin{equation} \label{Exp2}
h = \sum c_{l_1, \dots, l_{g+1}} q_0^{l_1} q_1^{l_2} \cdots q_g^{l_{g+1}}
\end{equation}  
with $c_{l_1, \cdots, l_{g+1}} \in \C \{ \r^\vee \cap M \}$, 
$0 \leq l_k \leq n_{k} -1$ for $k = 1, \dots, g$  and $l_{g+1} \in \Z_{\geq 0}$.

If $c_{l_1, \dots, l_{g  },0}$ and $c_{l'_1, \dots, l'_{g}, 0}$ are two different 
coefficients of the expansion the
Newton principal parts of 
$c_{l_1, \cdots, l_{g}, 0} q_0^{l_1}(\z) \cdots {q}_{g-1}^{l_{g}} (\z)$ 
and $c_{l_1', \cdots, l_{g}', 0} {q}_0^{l_1'}(\z) \cdots {q}_{g-1}^{l_{g}'}(\z)$
(viewed in the ring $\C \{ \r^\vee \cap M_g \}$) have no term in common. 
\end{Lem}
{\em Proof.}
The $q_g$-adic expansion of $h$ is of the form: 
$ h = a_0^{(g)}  +  a_1^{(g)}   q_{g} + \cdots + a_{s_g}^{(g)}   q_{g}^{s_g}$.
We build
the  $q_{g-1}$-adic expansions of the coefficients: 
$$ a_j^{(g)} = a_{0,j}^{(g-1)} +  a_{1,j}^{(g-1)}  q_{g-1} + \cdots + a_{s_{g-1}, j}^{(g-1)}   q_{g-1}^{s_{g-1}}$$
where $0 \leq \deg  a_{l, j}^{(g-1)} \leq  n_0 \cdots n_{g-1} -1 $ for $0 \leq l \leq s_{g-1}$ 
and $0 \leq   s_{g-1} \leq n_g -1$ since $a_j^{(g)}$  is of degree $< n_0  \cdots n_g =n$.
An expansion satisfying the required properties is obtained by iterating this procedure.
The unicity follows from the unicity of Euclidean division.
For the last assertion, remark that 
by lemma \ref{semigroup}
the Newton principal part of ${q}_{{k-1}} (\z)$ (viewed in $\C \{ \r^\vee \cap M_g \}$) is 
equal to $X^{\g_k}$ for $k =1, \dots,g$.  
It follows from  {\it 2} in lemma \ref{Azevedo} that
the Newton principal parts  of  $c_{l_1, \dots ,l_{g}, 0} {q}_0^{l_1} (\z)  \cdots {q}_{g-1}^{l_{g}} (\z) $
and of $c_{l_1', \dots ,l_{g}', 0} {q}_0^{l_1'} (\z) \cdots {q}_{g-1}^{l_{g}'} (\z)  $
do not have any term in common if $(l_1, \dots ,l_{g}) \ne  (l_1', \dots ,l_{g}')$.
\hfill $\ {\diamond}$

The following proposition (see \cite{Tesis}) 
generalizes  \cite{Mo}, Chapitre II, Th. 3.9. in the plane branch case.

\begin{Pro}  \label{Zariski} 
If $ h \in \C \{ \r^\vee \cap M \}[Y]$ is of degree  $ < n_0 n_1 \dots n_j$ then the Newton principal part of $h(\z)$
belongs to $\C[ \Gamma_j]$, for $j=1, \dots, g$.
\end{Pro}
{\em Proof.}
The result is trivial if  $\deg h =0$. 
If $\deg h < n_1 \dots n_{j}$ then 
the $(q_0, \dots, q_g)$-expansion of $h$ is of the form:
$h = \sum c_{l_1, \dots, l_{j}} q_0^{l_1} q_1^{l_2} \cdots q_{j-1}^{l_{j}}.$
By  lemma \ref{abhy}
the Newton principal parts of  $c_{l_1, \dots, l_{j}} q_0^{l_1} (\z) \cdots q_{j-1}^{l_{j}} (\z)$ and 
 of 
$c_{l_1', \dots, l_{j}'} q_0^{l_1'} (\z) \cdots q_{j-1}^{l_{j}'} (\z)$
 do not have terms in common, 
thus the polynomial ${h(\z)}_{| \cal N}$ is a sum of some of the terms 
in  the Newton principal parts of the summands 
 $c_{l_1, \dots, l_{j}} q_0^{l_1} (\z) \cdots q_{j-1}^{l_{j}} (\z)$ and  therefore
it belongs to $\C[ \Gamma_j]$ by proposition \ref{semigroup}, for $j=1, \dots, g$.
\hfill $\ {\diamond}$




We call the expansion (\ref{Exp2}) above the 
{\it $(q_0, \dots, q_g)$-expansion of $h$}.

\begin{Lem} \label{Abequations}
The $(q_0, \dots, q_g)$-expansion of $q_{j-1}^{n_j}$ is of the following form, for $1 \leq j \leq g$:
\begin{equation} \label{Abhexp3}
q_{j-1}^{n_j} =  c^*_{j} q_{j} 
          +  \sum  c_{l_1, \dots, l_{j}}^{(j)} q_0^{l_1} q_1^{l_2} \cdots q_{j-1}^{l_{j}}
\end{equation}  
where $ c^*_{j} \in \C^*$, 
the other coefficients belong to $\C \{ \r^\vee \cap M \}$, 
we have $0 \leq l_k \leq n_{k+1} -1$ for $k = 0, \dots, j-1$.
The coefficient 
 $c_{l_1^{(j)},  \dots, l_{j-2}^{(j)},0}^{(j)} $ appears and it is of the form
$X^{\a^{(j)}} \cdot \mbox{unit}$ where the integers $l_1^{(j)},  \dots, l_{j-2}^{(j)}$ and the exponent $\a^{(j)}$   are 
given by  formula  (\ref{canrel}).
Moreover, 
if  $X^{\a'}$ appears on the coefficient $c_{l_1, \dots ,l_{j}}^{(j)}$
then  
\begin{equation} \label{poid}
n_j \g_j \leq_\r \a' + l_1 \g_1 + \cdots + l_{j}  {\g}_{j}
\end{equation}
and equality holds if and only if 
$(l_1, \dots, l_{j}) = (l_1^{(j)},  \dots, l_{j-2}^{(j)},0)$ and $\a' = \a^{(j)}$.
\end{Lem}
{\em Proof.}
Since $\deg q_{j-1}^{n_j} = n_1 \cdots n_j$
the algorithm to calculate the $(q_0, \dots, q_g)$-expansion begins 
by dividing $q_{j-1}^{n_j}$ by $q_{j}$. This gives  $q_{j-1}^{n_j} = c^*_{j+1}  q_{j+1} + r_{j}$, 
where $c^*_{j+1} \in \C^*$ since both polynomials have the same degree.
The $q_k$ that may appear 
in the expansion of
$r_j$  are those of degree $\leq \deg r_{j} <n_1 \cdots n_j $.
We deduce from the  second assertion of lemma \ref{abhy} 
that $${\cal N } ( c_{l_1, \dots ,l_{j}}^{(j)} {q}_0^{l_1} (\z) \cdots {q}_{j-1}^{l_{j-1}} (\z)  ) \subset 
n_j \g_j + \r^\vee = {\cal N } ( q_{j-1}^{n_j} (\z) ). $$
This implies that if  $X^{\a'}$ appears on the coefficient $c_{l_1, \dots ,l_{j}}^{(j)}$
then  formula (\ref{poid}) holds. If equality in (\ref{poid}) holds for a term 
the term $X^{\a '}$  appearing on the series 
$c_{s_1^{(j)}, \dots ,s_{j}^{(j)}}^{(j)} (\z)$ 
it follows that the series is the form $X^{\a '} \cdot \mbox{unit}$.
Assertion {\it 2} of lemma \ref{Azevedo} implies 
that  $s_{j}^{(j)} = 0$ in the relation 
${n_j \g_j} = {\a_j} +  l_1  \g_1 + \cdots +  l_{j} \g_j $.
Then it follows that $(l_1,  \dots, l_{j-2}) = (l_1^{(j)},  \dots, l_{j-2}^{(j)})$
and that $\a ' = \a^{(j)}$ by unicity in (\ref{canrel}).
\hfill $\ {\diamond}$

\section{Partial embedded resolution with one toric morphism} \label{Teissier}
In this section we build a  partial embedded resolution of 
the toric quasi-ordinary germ 
embedded in an affine toric variety by using the semi-roots.
We follow the approach of 
\cite{Rebeca} for irreducible germs of plane curves.

We denote by $\Delta$ the cone 
$\r \oplus \R^g_{\geq 0} \subset (N_\D)_\R$ where $N_\D$ is the lattice 
$N \oplus \Z^g$ with dual lattice $M_\D$. 
We denote by 
$u_1, \dots, u_g$ the canonical basis of $ \{ 0 \} \oplus \Z^g$.
An element of   $\Delta^\vee \cap M_\D$ is of the form
$(\a , v)$ where $\a \in \r^\vee \cap M$ and  
$v = v_1  u_1^* +  \cdots  v_g u_g^*$ where  $ u_1^*, \dots, u_g^*$ is 
the dual basis of $u_1, \dots, u_g$ and $v_i \in \Z_{\geq 0}$.
We denote the monomial corresponding to $(\a , v)$ by  $X^\a U^v$,  $X^\a U_1^{v_1} \dots U_g^{v_g}$ or
$X^{\a + \sum v_i u_i^*}$ depending on the context.

The embedding $S \subset Z_\Delta$
which is studied in this section
corresponds algebraically to 
the  homomorphism of $\C  \{ \r^\vee \cap M \}$-algebras:
\begin{equation} \label{generators}
\Psi_0 : \C  \{ \r^\vee \cap M \} [ U_{1}, \dots, U_{g}] \rightarrow R \mbox{ given by } 
U_j \mapsto {q}_{j-1} (\z), \mbox{ for } j =1, \dots, g
\end{equation}
(which is surjective since in particular $R = \C  \{ \r^\vee \cap M \} [ q_0 (\z) ]$).

In the plane branch case Teissier shows that 
this embedding specializes to 
the {\it monomial curve}, an affine curve monomially embedded 
with the same semigroup 
(see \cite{Mono}). 
In the general case the generalization of monomial curve is given by 
an equivariant embedding $Z^\Gamma \subset Z_\Delta$
which is defined from 
the restriction of the lattice homomorphism
\begin{equation}
\varphi: M_\D \rightarrow M_g \mbox{ that maps } \a + v \mapsto \a + v_1 \g_1 + \dots + v_g \g_g
\end{equation}
to the semigroup $\Delta^\vee  \cap M_\D$ and its image $\Gamma$.

\subsection{Specialization through graded rings}

In the plane branch case the  embedding of the monomial curve 
is determined by a system of generators
the {\it graded ring}\footnote{See \cite{Bbk} for the definitions and properties of 
commutative algebra used in the following sections.} 
associated to the {\it filtration} of $R$ induced  by the powers of 
the maximal ideal of its  integral closure (see \cite{Mono}). 
In our case  we show that 
the  homomorphism $\Psi_0$ can be filtered in such a 
way that the homomorphism of the associated graded rings, 
forgetting the graded structure,  
defines the embedding $Z^\Gamma \subset Z_\Delta$ above.

The filtration of the ring $\C \{ \r^\vee \cap M \}$ (resp. of $\C [[ \r^\vee \cap M ]]$)
defined by a vector $\eta \in {\r}$ is given by the ideals:
$$ {\cal I}_j = \left\{ 
\sum_{ u \in \r^\vee \cap M } c_u X^u  \, / \, \min_{c_u \ne 0}  \langle \eta, u \rangle \geq j  \right\} 
\mbox{ for } {j \in \eta( \r^\vee \cap M)}.$$
Since the ring $\C \{ \r^\vee \cap M \}$ is Noetherian
the ordered sub-semigroup $\eta( \r^\vee \cap M)$ of $\R_{\geq 0}$ is 
isomorphic to $\Z_{\geq 0}$ (see the proof of Lemme 1.4 of \cite{GP}).
The vector $\eta $ defines a {\it weighted filtration} of  $\C \{  \r^\vee \cap M \} [ U_1, \dots, U_g ] $
(resp. of $\C \{  \Delta^\vee \cap M_\D \}$ or  $\C [[ \Delta^\vee \cap M_\D ]]$)
given by the ideals ${\cal J}_j$
generated by those series having only terms  $X^\a U^{v}$ of {\it weights} $ w:=\varphi (\a , v) $ such that
$\langle \eta, w \rangle \geq j$,  for $j $  running through the semigroup $\eta( \r^\vee \cap M_g)$.
The homomorphism $\Psi_0$ is {\it filtered}
since 
$\Psi_0 ({\cal J}_k) \subset {\cal I}_k$  for all $k \in \eta( \r^\vee \cap M_g)$,
and then it defines an homomorphism of the associated graded rings.

\begin{Pro} \label{Lejeune} 
The sequence of graded ring homomorphisms associated to the filtered sequence  
of homomorphisms (with the filtrations defined by $\eta \in \stackrel{\circ}{\r}$)
\begin{equation} \label{2dia}
\C \{ \r^\vee \cap M \} [U_1, \dots, U_g ] \stackrel{\Psi_0}{\longrightarrow}  R \hookrightarrow \C\{ \r^\vee \cap M_g \} \end{equation} 
is isomorphic to 
$$
\C [\r^\vee \cap M ] [U_1, \dots, U_g ]  \longrightarrow 
 \C [ \Gamma ]   \hookrightarrow   \C [\r^\vee \cap M_g ] 
$$
where the first homomorphism is defined by $X^\a U^v  \mapsto   X^{\varphi(\a, v)}$,
and the  graduations  are defined by $\eta$.
If the vector $\eta $ is irrational
   the semigroup 
$ \Gamma$ is determined by the graduation.
\end{Pro}
{\em Proof.}
If  $\eta \in \stackrel{\circ}{\r} $ the symbolic restriction $ \f_{|_\eta} $  of $\f \in \C \{ \r^\vee \cap M_g \}$ 
to the face defined by $\eta$ 
on the polyhedron ${\cal N}_\r (\f)$ belongs to  $\C [ \r^\vee \cap M_g  ]$ since this face is compact.
If   $\f  \in \C \{ \r^\vee \cap M_g  \}$  there exists  a unique integer
$k $ such that  $\f \in {\cal I}_{k} - {\cal I}_{k+1}$
and then we have $\f =  \f_{|_\eta} \mod {{\cal I}}_{k+1} $.
It follows from the property:
 $  \f_{|_\eta} \f'_{|_\eta} = (\f \f')_{|_\eta} $   for 
$0 \ne \f, \f' \in  \C \{ \r^\vee \cap M_g \}$, 
that  the  graded ring associated to this filtration 
is isomorphic to the graded ring $\C [ \r^\vee \cap M_g  ]$ where the $j$-homogeneous term of the graduation is 
$\bigoplus_{\langle \eta, u \rangle =j } \C X^u$ for ${j \in \eta(  \r^\vee \cap M_g )}$.
We deduce analogously that  the graded ring associated to the  weighted filtration 
is isomorphic to $\C [  \Delta^\vee \cap M_\D] $ where 
the non zero elements in the 
 $j$-homogeneous term are those polynomials
such that $\langle \eta, w \rangle = j$ for $w$ running through the 
weights of the monomials appearing on them.

Under these identifications we have that:

- The graded ring associated to $R$ with the induced filtration is 
isomorphic to 
the graded subring of $\C  [\r^\vee \cap M_g ]$  
generated as a $\C$-algebra by 
the symbolic restrictions $ \f_{| \eta} $ of $0 \ne \f \in R$ to the face defined by $\eta$ 
on the polyhedron ${\cal N} (\f)$. We deduce from proposition \ref{Zariski} and proposition 
\ref{semigroup} 
that this graded subring is equal to $ \C [ \Gamma ]$.

-  The initial term of ${{\Psi}_0}(U_i) = {q}_{i-1} (\z) $ is equal to 
$X^{\g_i}$ (the coefficient has been normalized to be one) thus 
the homomorphism $\mbox{gr} (\Psi_0)$
corresponds to the  $\C [\r^\vee \cap M ] $-homomorphism $\C [\r^\vee \cap M ][U_1, \dots, U_g] \rightarrow 
\C [\Gamma] $
that maps $U_i \mapsto X^{\g_i}$ for $i=1, \dots, g$.

If the vector $\eta$ is {\it irrational} we can recover the semigroup $ \r^\vee \cap M_g $ (resp. $\Gamma$) 
from the graduation of  $\C [ \r^\vee \cap M  ]$ (resp. of $ \C[\Gamma] $)
since each term of the graduation is of dimension one (resp. zero or one) over $\C$,
the vector  $\eta$ defining a total ordering on  $\r^\vee \cap M_g $.
\hfill $\ {\diamond}$

\begin{Rem}
The sequence of homomorphisms (\ref{2dia}) extends to the sequences:
\begin{equation} \label{dia}
\begin{array}{ccccc}
 \C [[ \Delta^\vee \cap M_\D ]]  &  \stackrel{\hat{\Psi}}{\longrightarrow}  &  \hat{R} & \hookrightarrow  &\C [[ \r^\vee \cap M_g ]]
\\
\uparrow &   & \uparrow & & \uparrow
\\
\C \{ \Delta^\vee \cap M_\D \}  & \stackrel{{\Psi}}{\longrightarrow}   &  R & \hookrightarrow & \C \{ \r^\vee \cap M_g \}
\end{array}
\end{equation}
where 
$\hat{R}$ denotes the completion of the ring $R$ with respect to 
the maximal ideal $\mathfrak{M}_R$.
The assertion of proposition \ref{Lejeune} remains true for each  line of the above diagram.
\end{Rem}
We notice that $\hat{R}$  coincides with the completion with respect to the filtration defined by 
$\eta $: we have that 
$\mathfrak{M}_R ^{s_j} \subset  {\cal I}_j $
where ${s_j}$ is the minimal  power of $\mathfrak{M}_R $
containing the set of monomials in $ {\cal I}_j - {\cal I}_{j+1}$ which is finite 
since $\eta \in \stackrel{\circ}{\r}$.

\subsection{Equations for the embeddings}

We build equations of the embeddings of
$Z^\Gamma \subset Z_\Delta$ and  $S \subset Z_\Delta$.

\begin{Pro} $\,$ \label{monomial}
The ideal of the embedding $Z^\Gamma \subset Z_\Delta$  is 
generated by the binomials 
\begin{equation} \label{h}
\left\{\begin{array}{lclclc}
h_1 & := & U_{1}^{n_1} & - &   X^{\a^{(1)}} \\
h_2 & := &  U_{2}^{n_2}  & - &  X^{\a^{(2)}} U_{1}^{l_1^{(2)}}, \\
\dots & \dots & \dots    & \dots & \dots, \\
h_g & := &  U_{g}^{n_g} & - &  X^{\a^{(g)}} U_{1}^{l_1^{(g)}} \dots  U_{g-1}^{l_{g-1}^{(g)}},
\end{array} \right.
\end{equation}
which correspond to relations (\ref{canrel}).
\end{Pro}
{\em Proof.} 
The ideal $I$ of the embedding $Z^\Gamma \subset Z_\Delta$ 
is generated by the binomials $X^\a U^\w - X^{\a'} U^{\w'}$ of $\C [ \Delta^\vee \cap M_\D ]$ 
verifying (see (\ref{eqembedding})):
\begin{equation} \label{eqem}
\varphi(\a , \w) = \varphi(\a' , \w') 
\end{equation} 
The binomials $h_1, \dots, h_g$ above verify this condition by lemma \ref{Azevedo}.
If  $B$ is a binomial in $I$, we can factor the common term in $U_g$ to obtain a binomial in $I$  of the
form $X^\a U^\w - X^{\a'} U^{\w'}$ with $w_g' =0$.
Then the integer  $\w_g$ is a multiple of $n_g$ 
(since $n_g \g_g \in M_{g-1} $ by lemma \ref{Azevedo} we obtain from the equality
(\ref{eqem}) a relation $r \g_g  \in M_{g-1} $ where 
$r$ is the reminder of the Euclidean division of
$\w_g$ by $n_g$ and then 
lemma \ref{Azevedo} implies that $r= 0$).
We can show by induction on $\w_g / n_g$ that 
the reminder of the Euclidean division of $X^\a U^\w - X^{\a'} U^{\w'}$
by $h_g$ as polynomials in $U_g$ is a binomial $B_1$ 
in $\C [ \r^\vee \cap M ] [U_1, \dots, U_{g-1} ]$. 
The binomial  $B_g$ obtained by iterating this 
procedure  belongs to $  \C [ \r^\vee \cap M ] $ and to the ideal $I$.
The relation (\ref{eqem}) corresponding to $B_g$ is trivial since 
the homomorphism $\varphi$ is injective on $M$ thus $B_g = 0$. 
This implies that the ideal 
$I$ is generated by $h_1, \dots, h_g$.
\hfill $\ {\diamond}$


\begin{Pro} \label{Abequations2}
The ideal
of the embedding $S \subset Z_\Delta$ defined by (\ref{generators}) is 
generated by elements of 
the ring  $\C \{ \r^\vee \cap M \} [U_1, \dots, U_g]$
which are of the form:
\begin{equation} \label{Abhexp4}
\left\{
\begin{array}{lclclclcl}
H_1 & : = & U_1^{n_1} & - &  X^{\a^{(1)}} & + &  c^*_{1} U_{2} & + &  r_1 (U_1), \\
H_2 & : = & U_2^{n_2} & - &  X^{\a^{(2)}} U_1^{l_1^{(2)} }   & + &  c^*_{2} U_{3} & + &  r_1 (U_1, U_2),  \\
\dots & \dots & \dots    & \dots & \dots &\dots & \dots  &\dots & \\
H_{g-1} & : = & U_{g-1}^{n_{g-1}}  & - &  X^{\a^{(g-1)}} U_{1}^{l_1^{(g-1)}} \dots  U_{g-2}^{l_{g-2}^{(g-1)}}
& + &  c^*_{g-1} U_{g}  & + &  r_{g-1} (U_1, U_2, \dots, U_{g-1}), 
\\
H_{g} & : = & U_{g}^{n_{g}}  & - & X^{\a^{(g)}} U_{1}^{l_1^{(g)}} \dots  U_{g-1}^{l_{g-1}^{(g)}}
&  & & + &  r_{g} (U_1, U_2, \dots, U_{g}).
\end{array} \right.
\end{equation}  
The weight of a term $X^{\a}  U_1^{l_1} U_2^{l_2} \cdots U_{j}^{l_{j}}$ appearing  
in the expansion of $ r_{j} (U_1, U_2, \dots, U_{j}) $ 
is $ \geq n_j \g_j $ and equality never holds, for $j=1, \dots, g$.   
The terms appearing in the expansion of $ r_{j} (U_1, U_2, \dots, U_{j}) $ 
are determined explicitly by formula (\ref{Abhexp3}).\end{Pro}
{\em Proof.}
It follows from the definition of the homomorphism $\Psi_0$ and formula (\ref{Abhexp3}) 
that the poly\-no\-mials $H_i$ above belong to the kernel of ${\Psi_0}$ (and then to the kernels of
$ {\Psi}$ and $\hat{\Psi}$).  By proposition \ref{monomial} and  lemma \ref{Abequations}
their initial forms with respect to the filtration defined by $\eta \in \stackrel{\circ}{\r}$ 
generate $\mbox{Ker} (\mbox{gr}(  \hat{\Psi}))$.
Then we have that $\mbox{gr} ( \mbox{Ker} (\hat{\Psi})) = \mbox{Ker} (\mbox{gr} (\hat{\Psi}))$.
We deduce using that the ideal $\mbox{ Ker} ( \hat{\Psi} )$  is complete 
  for the induced filtration, that the polynomials $H_1, \dots, H_g$ generate 
$\mbox{Ker} ( \hat{\Psi})$\footnote{See Proposition 12 No 9, \S 2, Chapitre III, of \cite{Bbk}.}.

Since the inclusion $\C \{ \Delta^\vee \cap M_\D \}  \rightarrow \C [[ \Delta^\vee \cap M_\D ]]$ is an homomorphism of 
local rings continuous for the $\mathfrak{M}$-adic topologies which extends to the identity homomorphism between the
respective completions, we have that  the ring  $\C [[ \Delta^\vee \cap M_\D ]]$ is a {\it faithfully flat}
$\C \{ \Delta^\vee \cap M_\D \}$-module\footnote{See  Proposition 10 No 5, \S 3, 
Chapitre III of \cite{Bbk}.}.
The ideal   $J $ generated  by $ (H_1, \dots, H_g)$ on $\C \{ \Delta^\vee \cap M_\D \}$
is contained in $\mbox{Ker}( {\Psi})$ and we have shown that 
$J \C [[ \Delta^\vee \cap M_\D ]] = \mbox{Ker} (\hat{\Psi})$.
The faithfully flat property implies that 
$J$ coincides with the contraction of $\mbox{Ker} (\hat{\Psi})$ 
in $\C \{ \Delta^\vee \cap M_\D \} $\footnote{See 
Proposition 9 No 5, \S 3,  Chapitre I  of \cite{Bbk}.}. Therefore we obtain that  $J = \mbox{Ker} ({\Psi})$.

Let  ${\cal U }$ be the subset of those elements in  $\C \{ \r^\vee \cap M \} [U]$ 
with non zero constant term as power series. 
The image by $\Psi_0$  of a series in  ${\cal U }$
is a unit.
This implies that the localization  
${\cal U }^{-1} \Psi_0 : {\cal U}^{-1} \C \{ \r^\vee \cap M \} [U] \rightarrow \C \{ \r^\vee \cap M_g \}$
is a well defined homomorphism. 
The same argument shows that $\mbox{Ker} ({\cal U }^{-1} \Psi_0)$ is generated by $H_1, \dots, H_g$.
Since ${\cal U } \cap \mbox{Ker} {\Psi}_0 = \emptyset$ we deduce from the standard properties
of the localization that   $H_1, \dots, H_g$ generate 
$\mbox{Ker} ({\Psi}_0)$. 
\hfill $\ {\diamond}$

\subsection{Simultaneous partial embedded resolution}

We show that the partial embedded resolution of $Z^\Gamma \subset Z_\Delta$ 
built in proposition \ref{binomial} is also a partial embedded resolution 
of $S \subset  Z_\Delta$.

The linear subspace $\ell \subset (N_\D)_\R$ orthogonal to $\mbox{Ker} (\varphi)$ 
is of dimension $d$
and is also orthogonal to the Minkowski sum of 
compact edges of ${\cal N} (h_i)$ for $i=1, \dots, g$.

\begin{Lem} \label{pseudo2}
Let $\Sigma_0$  be the smallest subdivision of $\Delta$ 
compatible with the Newton polyhedron of $H_1 \cdots H_g$.
The cone $\s_0 = \Delta \cap \ell$ belongs to  $\Sigma_0$.
The strict transform $S_{\Sigma_0} $ of $S$ is defined on the chart $Z_{\s_0} $ by the equations:
$U_i^{-n_i} H_i = 0$ for $i=1, \dots, g$. 
The intersection $S_{\Sigma_0} \cap \O_{\s_0} $ as schemes is 
reduced to the simple point $o_{\s_0}$.  
The germ $(S_{\Sigma_0},  o_{\s_0}) $ is isomorphic to the germ 
of toric variety $Z_{\s_0, N_{\s_0}}$ at the distinguished point. 
If $\Sigma$ is any subdivision of $\Delta$  
containing the cone $\s_0$  and if $\s \in \Sigma$ with $\stackrel{\circ}{\s} 
\subset \stackrel{\circ}{\Delta}$ then $S_{\Sigma} \cap  \O_\s \ne \emptyset $ implies that
$\s = \s_0$.
Moreover, if $\Sigma'$ is a regular subdivision of $\Sigma$
then the map $\p_{\Sigma'} \circ \p_{\Sigma} $ is an embedded  pseudo-resolution of $S$.
\end{Lem}
{\em Proof.}
A vector $v \in  \stackrel{\circ}{\s_0}$ vanish on $\mbox{Ker} (\varphi_\R)$ thus it
is of the form $v =  \tilde{v} \circ \varphi$ for 
$\tilde{v} \in N_g $ 
belonging  to $\stackrel{\circ}{\r}$ since $\tilde{v}$ vanishes only at the vertex of the cone $\r^\vee$
(this follows from  $\varphi^{-1}_\R ( \r^\vee ) = \Delta^\vee + \mbox{Ker} (\varphi_\R)$ and 
 $\stackrel{\circ}{\s_0} \subset \stackrel{\circ}{\Delta}$).
We deduce from this
that the face defined by $v$ on the polyhedron ${\cal N} (H_i)$
corresponds to the monomials 
of  weight $w$ such that $\langle \tilde{v}, w \rangle$ is  minimal.
By proposition \ref{Abequations2} 
the symbolic restriction of $H_i$ to this face is equal to $h_i$.
Conversely, if $h_i$ is the symbolic restriction of $H_i$ to the face defined by $v$ 
it follows that $v \in \ell $ and since these are compact faces we have that 
$v \in \stackrel{\circ}{\Delta}$ thus $v \in \stackrel{\circ}{\s_0} $.

The common zero locus $S'$ of the functions
$ U_i ^{-n_i} H_i$ for $i=1, \dots, g$ on the chart $Z_{\s_0}$
contains $S_{\Sigma_0} \cap Z_{\s_0}$.
Then we deduce from the proof of lemma \ref{nondeg} that:
\begin{equation} \label{weight}
 U_i ^{-n_i} H_i =  U_i ^{-n_i} h_i + \mbox{ terms vanishing on the orbit }
 \O_{\s_0} \mbox{, for } i=1, \dots , g.
\end{equation}
Since the equations $U_i ^{-n_i} h_i = 0 $ for $i=1, \dots , g$,  
define on the chart $Z_{\s_0}$ the strict transform 
$Z^{\Gamma}_{\Sigma_0}$ we deduce from (\ref{weight}) above that 
$S' \cap \O_{\s_0} $ coincides as schemes intersection 
with $  Z^{\Gamma}_{\Sigma_0} \cap \O_{\s_0}$, thus it is equal to 
the simple point $\o_{\s_0}$ by lemma \ref{pseudo}.
If the germ $(S', \o_{\s_0})$ is analytically irreducible it must coincide with the sub-germ
$(S_{\Sigma_0}, \o_{\s_0})$ since both are of the same dimension. 
We show this fact by proving that 
$(S', \o_{\s_0})$ is isomorphic to $(Z_{\s_0, N_{\s_0}}, o_{\s_0})$:

We notice that 
the chart $Z_{\s_0}$ is isomorphic to $\O_{\s_0} \times Z_{\s_0, N_{\s_0}}$ by (\ref{product}).
The binomials $ W_i := U_i ^{-n_i} h_i$ for $i=1, \dots, g$, 
define a regular system of parameters at the
point $o_{\s_0}$ of the orbit $\O_{\s_0}$ therefore 
we can apply lemma \ref{IFT} to the equations (\ref{weight}) to show
the existence of $\f_1, \dots, \f_g \in \C \{ \r^\vee \cap M_g \}$ 
such that 
the germ $(S', o_{\s_0} )$ is given by $W_i = \f_i$ for $i=1, \dots, g$.

Let $\Sigma$ be any subdivision of $\Sigma_0$ containing the cone $\s_0$.
The restriction  $\p:  S_\Sigma \rightarrow S$ of $\p_\Sigma$ 
is a modification and since 
$(S, o_\Delta)$ is analytically irreducible 
the exceptional fiber
is connected by the Main Theorem of Zariski (see \cite{Mumford} and \cite{Main}).
On the other hand we have 
that 
$$\p^{-1} (o_\Delta)  {=} \bigcup_{\s \in \Sigma,\stackrel{\circ}{\s} 
\subset \stackrel{\circ}{\Delta} } (S_{\Sigma} \cap  \O_\s) \mbox{ {\rm (by   (\ref{exceptional}));}}  $$
and we have 
shown that 
on the open set $S_\Sigma \cap Z_{\s_0}$ of $S_\Sigma$ the exceptional fiber is
reduced to the point $\o_{\s_0}$, therefore the exceptional fiber $\p^{-1} (o_\Delta)$
contains no other points (otherwise would not be a connected set).

If $\Sigma'$ is a regular subdivision of $\Sigma$ 
it follows that $\O_\s \cap S_{\Sigma'} \ne \emptyset$
if and only if $\s \subset \s_0$. Thus we can cover the
strict transform with those charts $Z_\s$ for 
$\s \subset \s_0$ and $\dim \s = \dim \s_0$.
It follows as in the case when $\s_0 $ is a regular cone, 
that the strict transform is smooth
and transversal to the  canonical stratification of 
the exceptional divisor therefore $\p_{\Sigma'} \circ \p_{\Sigma}$
is a pseudo-resolution. 
\hfill $\ {\diamond}$

\begin{The} \label{Bernard}
Let $\Sigma$ be any subdivision of $\Delta$  
containing the cone $\s_0$. 
\begin{enumerate}
\item
The strict transform $S_\Sigma$ is a germ at the point $o_{\s_0}$  isomorphic
to $(Z_{\s_0, N_{\s_0}}, o_{\s_0})$ and 
the restriction $\p_\Sigma |S_\Sigma : S_\Sigma \rightarrow S$ is the normalization map.
\item
The morphism $\p_\Sigma$ is a partial embedded resolution of $S \subset Z_\Delta$.
\end{enumerate}
\end{The}
{\em Proof.}
The first assertion follows from lemma \ref{pseudo2} taking in account (\ref{saturation}) 
(which implies that $(Z_{\s_0, N_{\s_0}}, o_{\s_0}) $ is isomorphic to $(Z_{\r, N_g}, o_\r)$)
and proposition  \ref{closure} which implies that the integral closure of $R$ is 
 $ \C \{ \r^\vee \cap M_g \}$.

By lemma \ref{pseudo2} if  $\Sigma' $ is a regular subdivision of  $\Sigma$ 
the map $\p_{\Sigma'} \circ \p_{\Sigma}$ is a embedded pseudo-resolution of $ S$;
we show that if $\Sigma' $ is a resolution of the fan $\Sigma$ then the restriction
$S_{\Sigma'} \rightarrow S$ is a resolution of singularities.

The germ  $S_\Sigma $ is parametrized by $W_i = \f_i$ for $i =1, \dots, g$, 
on the chart $Z_{\s_0} \cong \O_{\s_0}  \times    Z_{\s_0, N_{\s_0}} $ thus the 
restriction  $S_\Sigma   \rightarrow Z_{\s_0, N_{\s_0}} $  of the second projection is 
an isomorphism of germs. It follows that the singular locus of $S_\Sigma$ 
lies over the singular locus of
the toric variety 
$Z_{\s_0, N_{\s_0}} $.
It is easy to see that the orbit $\O_{\t}$ of $Z_{\s_0}$ is the set lying over the orbit $\O_{\t, N_{\s_0}} $ of 
$Z_{\s_0, N_{\s_0}}$ thus 
the singular locus of $S_\Sigma$ is equal to 
$\bigcup (S_\Sigma \cap \O_\t)$ for 
$\t$ running through the set of non regular faces of $\s_0$.
If  $\Sigma' $ is a resolution of the fan $\Sigma$
 and if $\s' \in \Sigma$
is a regular cone then $\s' \in \Sigma'$, thus 
$Z_\Sigma \rightarrow Z_{\Sigma_0}$ is an isomorphism over 
the points of the orbit $\O_{\s'}$ 
by (\ref{discriminantlocus}). Therefore
the restriction $S_{\Sigma'}  \rightarrow S_{\Sigma}$  
is an isomorphism outside the singular locus of $S_{\Sigma}$
and since  $S_{\Sigma'}$ is smooth this modification is a resolution of singularities of the normalization $S_{\Sigma}$.
A forteriori the composed map $S_{\Sigma'}  \rightarrow S_{\Sigma}$ is resolution of singularities of $S$.
\hfill $\ {\diamond}$

\subsection{Relation between the partial embedded resolution procedures}

We show that the partial embedded resolutions of an analytically irreducible 
toric quasi-ordinary germ $S$ defined in section  \ref{Oka} 
and \ref{Teissier} coincide when the second is suitably chosen.

In section \ref{Oka} we have built a partial embedded resolution 
$\p$
of a toric quasi-ordinary hypersurface $S \subset Z_\varrho$
which depends only on the characteristic exponents of a toric quasi-ordinary polynomial $f$
defining the embedding. 
Since the germ $S$ is analytically irreducible,  
the morphism  $\p$ is the composition of
$g$ toroidal modifications
$\p_i : Z_{i} \rightarrow Z_{i-1}$ for $i=1, \dots, g$ and $g$ 
the number of characteristic exponents.
In section \ref{Teissier} we have built an embedding of 
$(S, o)$ as a codimension $g$ sub-germ of the toric variety 
$(Z_\D, \o_\D)$ and we have proved that 
if $\Sigma$ is a subdivision of $\D$ compatible with certain linear subspace,
the toric morphism $\p_\Sigma : Z_\Sigma \rightarrow Z_\D$ is 
partial embedded resolution of $S \subset Z_\D$. 
Furthermore, the restriction of $\p$  (resp. of $\p_\Sigma$) to
the strict transform $S'$ (resp. $S_\Sigma$) of $S$ is the normalization map
(see theorems   \ref{hres} and  \ref{Bernard} ).

The embedding $S \subset Z_\D$ 
defined by (\ref{generators})
extends to an embedding of the pair 
$(S, Z_\varrho)$: 
the image of $ (Z_\varrho, o_\varrho)$
under this embedding 
is the sub-germ 
$({\mathcal Z}, o_\D)$  of $(Z_\D, \o_\D)$
defined by 
the equations (see (\ref{Abhexp4})):
\begin{equation} \label{embed} 
\left\{\begin{array}{lclclcl}
-  c^*_{1} U_{2} &  = & U_1^{n_1} & - &  X^{\a^{(1)}} & + &  r_1 (U_1), \\
- c^*_{2} U_{3}   &  = & U_2^{n_2} & - &  X^{\a^{(2)}} U_1^{l_1^{(2)} }   & + &  r_1 (U_1, U_2),  \\
\dots & \dots & \dots    & \dots & \dots &\dots & \dots  \\
- c^*_{g-1}  U_{g}   &  = & U_{g-1}^{n_{g-1}}  & - &  X^{\a^{(g-1)}} U_{1}^{l_1^{(g-1)}} \dots  
U_{g-2}^{l_{g-2}^{(g-1)}}
& + &  r_{g-1} (U_1, U_2, \dots, U_{g-1}), 
\end{array} \right.
\end{equation}  
Since $c^*_i \in \C^*$ we can eliminate the variables 
$U_2, \dots, U_g$ 
in the equation:
\begin{equation} \label{embed2}
U_{g}^{n_{g}}   -  X^{\a^{(g)}} U_{1}^{l_1^{(g)}} \dots  U_{g-1}^{l_{g-1}^{(g)}}
+   r_{g} (U_1, U_2, \dots, U_{g}) = 0 
\end{equation}  
by using (\ref{embed}),
and we obtain in this way a  quasi-ordinary polynomial defining the embedding 
$S \subset Z_\varrho$.

\begin{Rem} \label{41}
If we vanish the $r_1, \dots, r_g$ in (\ref{Abhexp4}) we obtain:
\begin{equation} \label{Abhexp5}
\left\{\begin{array}{lclclcl}
\tilde{H}_1 & : = & U_1^{n_1} & - &  X^{\a^{(1)}} & + &  c^*_{1} U_{2} , \\
\tilde{H}_2 & : = & U_2^{n_2} & - &  X^{\a^{(2)}} U_1^{l_1^{(2)} }   & + &  c^*_{2} U_{3} ,  \\
\dots & \dots & \dots    & \dots & \dots &\dots & \dots \\
\tilde{H}_{g-1} & : = & U_{g-1}^{n_{g-1}}  & - &  X^{\a^{(g-1)}} U_{1}^{l_1^{(g-1)}} \dots  U_{g-2}^{l_{g-2}^{(g-1)}}
& + &  c^*_{g-1} U_{g}, 
\\
\tilde{H}_{g} & : = & U_{g}^{n_{g}}  & - & X^{\a^{(g)}} U_{1}^{l_1^{(g)}} \dots  U_{g-1}^{l_{g-1}^{(g)}}.
&  & 
\end{array} \right.
\end{equation}  
We can eliminate recursively from the equations  $\tilde{H}_i = 0$, for $i =1, \dots, g-1$
the variables  $U_2, \dots, U_g$ 
in the equation 
$\tilde{H}_g = 0$
obtaining in this way a canonical equation of a quasi-ordinary hypersurface  
with the same characteristic monomials.
The exponents  appearing in these polynomials are completely determined by the 
characteristic monomials of $(S, 0)$. See \cite{Mono} for Teissier's analogous statement 
in the case of plane branches.
\end{Rem}

\begin{Defi} \label{suitable}
A subdivision $\Sigma$  of $\D$ is suitable with respect to the embedding of the pair  $({\mathcal Z}, S)$ 
in $Z_\D$, if it is the dual Newton diagram of $\tilde{H}_1 \cdots \tilde{H}_g$.
\end{Defi}

It follows from remark \ref{41} that the suitable    subdivision $\Sigma$ of $\D$ is 
uniquely determined from the given characteristic monomials of $(S, 0)$.
We prove that the strict transform ${\mathcal Z}_\Sigma$ of ${\mathcal Z}$ by
the toric modification $\p_\Sigma $ is a {\it section } of 
the toric variety $ Z_\Sigma$,  transversal to the exceptional fiber of 
the modification $\p_\Sigma $. 
More generally 
it is  transversal to the orbit stratification of $ Z_\Sigma$ and 
the set of non empty  intersections ${\mathcal Z}_\Sigma \cap \O_\s$ define 
the stratification corresponding to a natural toroidal embedding
structure which is determined by $\Sigma$.
In particular we obtain that the restriction $p: {\mathcal Z}_\Sigma \rightarrow {\mathcal Z}$
of $\p_\Sigma$ to ${\mathcal Z}_\Sigma$ is a partial embedded resolution of $S \subset {\mathcal Z}$.
The main result of this section is that the partial embedded resolutions
defined by $\p$ and by $p$ are isomorphic:

\begin{The} \label{mismo}
If $\Sigma$ is the suitable subdivision of $\D$ with respect to the 
embedding of the pair  $({\mathcal Z}, S)$ 
in $Z_\D$, then the strict transform  $(S_\Sigma, {\mathcal Z}_\Sigma)$
of the pair $({\mathcal Z}, S)$ by the toric modification $\p_\Sigma$  
is equal to an embedding of the pair $(S', Z_g)$ in $Z_\Sigma$ 
such that the following diagram commutes:
\[
\begin{array}{ccc}
(S', Z_g)  & \stackrel{\cong}{\longrightarrow} & (S_\Sigma, {\mathcal Z}_\Sigma)
\\
\pi \downarrow & &  p \downarrow 
\\
(S, Z_\varrho) &  \stackrel{\cong}{\longrightarrow}  & (S, {\mathcal Z}) 
\end{array}
\]
Therefore the morphism $p: {\mathcal Z}_\Sigma \rightarrow {\mathcal Z}$ 
is the composition of $g$ toroidal modifications.
\end{The}
In the plane branch  case  an analogous statement (using resolution instead of partial resolution)
has been announced  by Goldin and Teissier without proof in \cite{Rebeca};
Lejeune and Reguera have sketched in that case toric resolutions of the monomial curve 
such that the restrictions to  the strict transform of the smooth surface, which contains 
the re-embeded plane branch, are equal to the minimal resolution of the branch (see \cite{LR2}).

We introduce first some notations in order to describe the {\it suitable} subdivision $\Sigma$ of $\D$.
The following subsets of $\D$ defined for $0 \leq j < j+k \leq g $
\begin{equation} \label{deltas}
\r_j^{j+k} = \{ a + \langle a, \g_1 \rangle u_1 + \cdots + \langle a, \g_j \rangle u_j+ 
n_j \langle a,  \g_{j} \rangle u_{j+1} + \cdots + 
n_j\cdots n_{j+k -1}  \langle a,  \g_{j} \rangle  u_{j+k} / a \in \r \}  
\end{equation}
are the cones 
which correspond by duality  to certain  Minkowski sums of 
edges of ${\cal N} (\tilde{H}_i)$ for $i=1, \dots, k$.
The cone $\r_0^{k}$ coincides with $\r \times \{ 0 \} \subset \Delta$ for $1 \leq k \leq g$. 
We denote by $\Xi$ the $(d+1)$-dimensional fan whose elements are 
the faces of the  $2g$ cones of dimension $d+1$ :
\begin{equation} \label{2g}
\r_j^{g} + \R_{\geq 0} u_j, \quad  \r_j^{g} + \r_{j-1}^g \quad \mbox{ for } j=1, \dots, g.
\end{equation} 
We will show below that $\Xi$ is a subfan of the suitable subdivision $\Sigma$ (see remark \ref{xi}).

\begin{Pro} \label{combo}
Let $\Sigma$ a suitable subdivision of $\D$.
If $\s \in \Sigma$ and if $\stackrel{\circ}{\s} \subset \stackrel{\circ}{\D}$ 
then ${\mathcal Z}_\Sigma \cap \O_\s \ne \emptyset $ implies that $\s \in \Xi$.
If $\s \in \Xi^{(d+1)}$ then ${\mathcal Z}_\Sigma \cap \O_\s$ is reduced to a simple point 
$x_\s$ and the germ $({\mathcal Z}_\Sigma, x_\s)$ is isomorphic to 
$(Z_{\s, (N_\D )_\s}, o_\s)$. The set $\{ {\mathcal Z}_\Sigma \cap \O_\s \}_{ \s \in \Xi } $ is the stratification 
associated to a toroidal embedding structure on ${\mathcal Z}_\Sigma$ which has $\Xi$ as associated 
conic polyhedral complex.
\end{Pro}

In order to prove proposition \ref{combo} we characterize in the lemma below  some
convexity properties of the the Newton polyhedra of the polynomials 
$H_1, \dots,H_{g-1}$ defining
the embedding ${\mathcal Z} \subset Z_\D$. 
Lemma \ref{sand} below is inspired by a 
result of Lejeune and Reguera  in the case of  sandwiched surface singularities
(see Proposition 1.3 of \cite{LR}).
We need some useful notations. The exponents: 
$$ u_{j+1}^*,\, n_j u_j^*, \, \varpi_j : = \a^{(j)} +  l_1^{(j)} u_1^* + \cdots + 
{l_{j-1}^{(j)}} u_{j-1}^* $$
are the vertices of the two dimensional face ${\cal T}^j$ of the polyhedron ${\cal N} (H_j)$
by proposition 
\ref{Abequations2}.
This face and its edges 
$$ {\cal T}^j_1 := [ u_{j+1}^*,  n_j u_j^* ], \,
 {\cal T}^j_2 := [ u_{j+1}^*,\varpi_j ], \, {\cal T}^j_3 := [  n_j u_j^*, \varpi_j ],$$
 play a significant role in what follows.
Any other vertex $\varpi_j '$ of the Newton polyhedron of $H_j$ corresponds to a monomial of 
weight $> n_j \g_j$, i.e., we have $\varpi_j ' = \a' + l_1 u_1^* + \cdots + 
l_{j} u_{j}^*$ and 
$\a'+ l_1 \g_1 + \cdots +  l_{j} \g_j > n_j \g_j$.
 
\begin{Lem} \label{sand}
 Let ${\cal E}_i$ be 
a compact edge of ${\cal N} (H_i)$ for $i=1, \dots,  g$.
If  ${\cap_{i=1}^{j}} \stackrel{\circ}{ \s} ({\cal E}_i) \ne \emptyset$ for  $1 \leq j \leq g$,
then we have ${\cal E}_i = {\cal T}^i_{s(i)}$ where $s(i)$ is given by:
\[
\left\{
\begin{array}{l}
s(i) = 3, \, 1 \leq i \leq  j \leq g-1, \quad \quad \hfill \mbox{{\rm (A)}}
\\
s(i) = 3 , \,  1 \leq i \leq  j_0 -1; \, s(i) = 1, \mbox{ for } j_0 \leq i \leq  j \leq g-1 ; \, 
\mbox{ for } 1 \leq j_0 \leq g-1, \quad  \quad \hfill \mbox{{\rm (B)}}
\\
s(i) = 3 , \, 1 \leq i \leq  j_0; \, s(j_0) = 2; \, s(i) = 1, j_0 +1 \leq i \leq  j \leq g-1 ; \, 
\mbox{ for } 1 \leq j_0 \leq g-1,  \quad \quad  \hfill \mbox{{\rm (C)}} 
\\
s(i) = 3, \, 1 \leq i \leq j = g, \quad \quad \hfill \mbox{{\rm (D)}}
\end{array}\right. \]
Moreover, the intersection $\cap_{i=1}^{j} \s({\cal T}^i_{s(i)})$ is equal to:
\[
\left\{
\begin{array}{lllll}
\r_j^{j+1} & + & \R_{\geq 0} u_{j+1} & + & \R_{\geq 0} u_{j+2} + 
\cdots + \R_{\geq 0} u_{g},\quad  \quad   \hfill  \mbox{ in the case } \mbox{{\rm (A)}}
\\
\r_{j_0}^{j+1} &  + &   \R_{\geq 0} u_{j_0} & + &  \R_{\geq 0} u_{j+2} +
 \cdots + \R_{\geq 0} u_{g},\quad  \quad  \hfill  \mbox{ in the case } \mbox{{\rm (B)}}
\\
\r_{j_0}^{j+1}  & + & \r_{j_0 -1}^{j+1} & + & \R_{\geq 0} u_{j+2} + 
\cdots + \R_{\geq 0} u_{g},\quad  \quad   \hfill \mbox{ in the case } \mbox{{\rm (C)}} 
\\
\r_{g}^{g}  & & & &  \hfill \mbox{ in the case } \mbox{{\rm (D)}} 
\end{array}\right. 
\]
\end{Lem}
{\em Proof.} 
The compact faces of Newton polyhedra are determined
by elements $a + v \in \Delta$ which belong to the relative interior of $\Delta$;
i.e., $a + v$ is of the form $a \in \stackrel{\circ}{\r}$ and $ v = \sum_{i=1}^{g} v_i u_i$
with $v_i >  0$. We calculate the values of $a + v $  
on the vertices of the Newton polyhedron of $H_j$ in terms of the weight of the corresponding monomial. 
We prove the lemma by induction on $j$, for $j =1 $ we show first that 
the compact edges of ${\cal N} (H_1)$ are exactly ${\cal T}^1_i$ for $i=1, 2, 3$.
We have the following:
\[
\begin{array}{l}
\mbox{{\rm (i)}} = \langle a+v , \varpi_1 \rangle = \langle a, \a^{(1)} \rangle = n_1 \langle a, \g_1 \rangle
\\
\mbox{{\rm (ii)}} = \langle a+v, n_1 u_1^* \rangle = n_1 v_1
\\
\mbox{{\rm (iii)}} = \langle a+v, u_2^* \rangle =  v_2
\\
\mbox{{\rm (iv)}} = \langle a+v, \varpi_1'  \rangle = \langle a+ v , \a_1 + l_1 u_1^*  \rangle = \langle a, \a' + l_1 \g_1 \rangle  +
l_1 ( v_1 -  \langle a,  \g_1 \rangle ) {>} (n_1 - l_1) \langle a,  \g_1 \rangle  + l_1 v_1
\end{array}
\]
where the inequality on $\mbox{{\rm (iv)}}$ follows from (\ref{poid}) since $a \in \stackrel{\circ}{\r}$.
We suppose that $a +v$ determines a compact edge $e_1$  of ${\cal N} (H_1)$. 
Three cases appear:

- If $v_1 =  \langle a, \g_1 \rangle$ then $\mbox{{\rm (iv)}} > \mbox{{\rm (i)}} =  \mbox{{\rm (ii)}} $ thus  $ v_2  >   n_1 \langle a, \g_1 \rangle$ and ${\cal E}_1 = {\cal T}^1_3$. 

- If $v_1 >  \langle a, \g_1 \rangle$ then $\mbox{{\rm (ii)}}, \mbox{{\rm (iv)}} > \mbox{{\rm (i)}} $ thus  
$ v_2 = n_1 \langle a, \g_1 \rangle$ and ${\cal E}_1 = {\cal T}^1_2$.

- If $v_1 <  \langle a, \g_1 \rangle$ then  $\mbox{{\rm (i)}}, \mbox{{\rm (iv)}} > \mbox{{\rm (ii)}} $
thus $v_2 = n_1 v_1$ and ${\cal E}_1 = {\cal T}^1_1$.

The equality $\mbox{{\rm (i)}} =  \mbox{{\rm (ii)}} = \mbox{{\rm (iii)}} $
corresponds to the two dimensional face ${\cal T}^1$.
It follows that:
\[  \s( {\cal E}_1 ) =
\left\{
\begin{array}{lllll}
\r_1^{2} & + & \R_{\geq 0} u_{2} & + & \R_{\geq 0} u_{3} + 
\cdots + \R_{\geq 0} u_{g},\quad  \quad   \hfill  \mbox{ if  } {\cal E}_1 = {\cal T}^1_3
\\
\r_1^{2} & + & \R_{\geq 0} u_{1} & + & \R_{\geq 0} u_{3} + 
\cdots + \R_{\geq 0} u_{g},\quad  \quad   \hfill  \mbox{ if  } {\cal E}_1 = {\cal T}^1_2
\\
\r_1^{2} & + & \r_0^{2}   & + & \R_{\geq 0} u_{3} + 
\cdots + \R_{\geq 0} u_{g},\quad  \quad   \hfill  \mbox{ if  } {\cal E}_1 = {\cal T}^1_1
\end{array}\right. 
\]

We suppose the result true for $j-1$.
We consider a vector $a +v \in \bigcap_{i=1}^{j-1} 
\stackrel{\circ}{\s} ( {\cal T}^i_{s(i)} )$ determining
an edge ${\cal E}_j$ of ${\cal N} (H_{j})$, i.e., $ a + v \in \; \stackrel{\circ}{\s}({\cal E}_j)$.
The values of  $ a + v $ on the vertices of ${\cal T}^j_1$ are:
\[
\begin{array}{l}
\mbox{{\rm (ii)}} = \langle a+v, n_j u_j^* \rangle = n_j v_j
\\
\mbox{{\rm (iii)}} = \langle a+v, u_{j+1} ^* \rangle =  v_{j+1}
\end{array}
\]

We deal first with the case  (A) where $s(i) = 3 $ for $1 \leq i \leq j-1$.
Then $a+ v \in \r_{j-1}^j$ and it follows as before that:
 \[
\begin{array}{l}
\mbox{{\rm (i)}} = \langle a+v, \varpi_j \rangle =  n_j \langle a, \g_j \rangle
\\
\mbox{{\rm (iv)}} = \langle a+v, \varpi_j' \rangle = \langle a+v, \a_j + \sum_{i=1}^j l_i u_i^* \rangle 
> n_j \langle a , \g_j \rangle + l_j (v_j - \langle a,  \g_j \rangle)
\end{array}
\]
where the inequality is obtained from (\ref{poid}) by adding and subtracting the term $l_j \langle a,  \g_j \rangle$.
Three cases appear if $v_j$ is $=$ (resp. $>$ or $<$) to  $ \langle a, \g_j \rangle$
and we obtain the result by arguing as in the case $j = 1$ by replacing appropriately the index $1$ by $j$.

In any other case  by induction hypothesis  there exists $1 \leq j_0 \leq j-1$ such that 
$s(i) = 3 $ for $ 1 \leq i < j_0$ and $s(j_0) \in \{1, 2 \}$.
It follows that the vector $a +v$ is of the form:
$$
a +v = a + \sum_{i= 1}^{j_0} \langle a, \g_i \rangle u_i + \sum_{i > j_0}^g v_i u_i
$$

We bound 
the value of $a +v$ on a vertex of the polyhedron ${\cal N} (H_j)$ 
not lying on ${\cal T}^j_1$. 
\begin{equation} \label{calculo}
\begin{array}{lll}
\mbox{{\rm (iv)}} & = &  \langle a+v, \varpi_j' \rangle = \langle a+v, \a_j + \sum_{i=1}^j l_i u_i^* \rangle \;  =
\\
& = & 
 \langle a , \a_j + \sum_{i=1}^j l_i \g_i  \rangle - \sum_{i=j_0}^j l_i \langle a ,\g_i  \rangle
+ \sum_{i=j_0}^j l_i v_i  \; > 
\\
& > & 
(n_j -l_j) \langle a ,\g_j   \rangle - \sum_{i=j_0}^{j-1} l_i \langle a ,\g_i  \rangle + \sum_{i=j_0}^j l_i v_i  
\; >  \cdots > 
\\
&  > &  \left( ( \cdots ((n_j - l_j) n_{j-1} - l_{j-1} ) \cdots ) n_{j_0} - l_{j_0} \right) \langle a ,\g_{j_0}  \rangle 
+ \sum_{i=j_0}^j l_i v_i 
\end{array}
\end{equation} 
The first inequality is given by (\ref{poid}) and the others 
are deduced from  the inequality $n_i \g_i < \g_{i+1}$ in lemma \ref{Azevedo}. 

\noindent
In case (B)               by induction hypothesis we have that 
$v_{j_0} = \langle a, \g_{j_0} \rangle  + c $ for some $c > 0$ and 
$v_i = n_{j_0} \cdots n_{i-1} \langle a, \g_{j_0} \rangle $ for $j_0 < i \leq j$.
In case (C)   we have that 
\[
\left\{
\begin{array}{cccclr}
n_{j_0 -1 } \langle a, \g_{j_0 -1} \rangle & < & v_{j_0}  &  < &  \langle a, \g_{j_0} \rangle  & \mbox{ if } j_0 >1 
\\
0  & < &  v_{j_0} & < & \langle a, \g_{j_0} \rangle  & \mbox{ if } j_0 = 1 
\end{array} \right.
\]
and that $v_j = n_{j_0} \cdots n_{i-1} v_{j_0} $ for $j_0 < i \leq j$.
In both cases  (B)   and (C)  
when   substitute the $v_i$  on  (\ref{calculo}) 
we deduce that 
$\mbox{{\rm (iv)}}, \mbox{{\rm (i)}} > \mbox{{\rm (ii)}}$ therefore 
$v_{j+1} = n_j v_j$ and 
${\cal E}_j = {\cal T}^j_1$.

Finally, when $j=g$ the polynomial $H_g$ has no term in $U_{g+1}$.
In particular a vector $a + v \in \bigcap_{i=1}^{g-1} \stackrel{\circ}{\s}( {\cal T}^i_{s(i)})$ 
for $s(i)$ in case (B) or (C), 
determines the vertex $n_g u_g^*$ of the polyhedron ${\cal N} (H_g)$. 
The only remaining case is (A) and then the condition on $a+v$ to determine a compact
edge of ${\cal N} (H_g)$ is $v_g = \langle a, \g_g \rangle$;
the edge is equal to ${\cal T}^g_3 = [\varpi, n_g u_g^*]$ and 
 $\bigcap_{i=1}^{g} = \r_{g}^g$.
\hfill $\ {\diamond}$

\begin{Rem} \label{xi}
The cones of the form $\cap_{i=1}^{g-1} \s({\cal T}^i_{s(i)})$ 
defined by lemma \ref{sand} when $j = g -1$ are 
\[
\left\{
\begin{array}{l}
\r_{k}^{g} +   \R_{\geq 0} u_{k}, \, \,  \r_{k}^{g}   +  \r_{k -1}^{g} \quad  \mbox{ for } k=1, \dots, g-1.
\\
\r_{g-1}^{g} +   \R_{\geq 0} u_{g}
\end{array}\right. 
\]
If we subdivide  $\r_{g-1}^{g} +   \R_{\geq 0} u_{g}$ with $\r_{g}^g$
we obtain the fan  $\Xi$.
It follows that $\Xi$ is a subfan of the dual Newton diagram of 
$\tilde{H}_1 \dots \tilde{H}_g$. 
Theorem \ref{mismo} holds more generally for any subdivision of $\D$  containing $\Xi$. 
\end{Rem}

\medskip

\noindent
{\em Proof of Proposition \ref{combo}.}
Let $\Sigma'$ be any subdivision of $\Sigma$ which is compatible with the 
Newton polyhedra of $H_1, \dots, H_{g-1}$ and  $\s \in \Sigma '$ 
with $\stackrel{\circ}{\s} \subset \stackrel{\circ}{\D}$.
By lemma \ref{nondeg} a necessary condition  to have 
${\mathcal Z}_{\Sigma'} \cap \O_\s \ne \emptyset$
is that $\s$ determines a face ${\cal F}_i$ of dimension $\geq 1$ of each polyhedron 
${\cal N} (H_i)$ for $i=1, \dots, g-1$. 
Then we have 
 ${\s} \subset  \bigcap_{i=1}^{g-1} {\s} ( {\cal F}_i ) 
\subset \bigcap_{i=1}^{g-1} {\s} ( {\cal E}_i )  $
for ${\cal E}_i$ any fixed edge of the face ${\cal F}_i$.
The possible edges ${\cal E}_i$ that may appear are determined by lemma \ref{sand}
and by duality $\s$ is contained in the support of $\Xi$.
By using (\ref{exceptional}) we deduce that if $\s \in \Sigma - \Xi$ 
then $\O_\s \cap {\mathcal Z}_{\Sigma} = \emptyset$.

The proof of the second assertion is analogous to the proof of lemma \ref{pseudo2}.
Let $\s \in \Xi^{(d+1)}$, for instance $\s = \r_j^g + \R_{\geq 0} u_j$ (the proof in the case 
$\s = \r_{j-1}^g + \r_j^g $ is analogous) for $j=1, \dots, g$.
The common zero locus ${\mathcal Z}' \subset Z_\s$ of the set functions 
$X^{-m_1} H_1, \dots, X^{-m_{g-1}} H_{g-1} $ for 
\[
m_i = 
\left\{
\begin{array}{lc}
\varpi_i  & \mbox{ if } i=1, \dots, j-1
\\
{u_{i+1}^*} & \mbox{ if } i=j, \dots, g-1
\end{array}
\right.
\] 
contains  ${\mathcal Z}_\Sigma \cap Z_\s$. 
Each series $X^{-m_i} H_i$  is of the form 
$$
X^{-m_i} H_i =  B_i  +
 \mbox{ terms vanishing on } \O_\s $$
where 
\[
B_i  = 
\left\{
\begin{array}{lc}
1 - X^{n_i u_i^* - m_i} & \mbox{ if } i=1, \dots, j-1
\\
c_i^* + X^{n_i u_i^* - m_i} & \mbox{ if } i=j, \dots, g-1.
\end{array}
\right.
\] 
The edge ${\cal E}_i := [ n_i u_i^*, m_i ]$ 
is a face of the polyhedron ${\cal N} (H_i)$ 
and by lemma \ref{sand} we have 
$\s = \bigcap_{j=1}^{g-1} \s ( {\cal E}_i )$ thus  $\s^\bot  = \bigoplus_{j=1}^{g-1} (\s ( {\cal E}_i )) ^\bot$
since the edges ${\cal E}_i$ are affinely independent.
Moreover, the vector $n_i u_i^* - m_i$ is primitive for the lattice $M_\D$ and it follows
that $\s^\bot \cap M_\D = \oplus_{i=1}^{g-1} ( n_i u_i^* - m_i ) \Z$.
It follows that the intersection ${\mathcal Z}' \cap \O_\s$ as schemes, defined by the equations  $B_1 = \cdots = B_{g-1} = 0$, 
is a simple point $x_\s$  and that         $B_1, \dots, B_{g-1}$ define
a regular system of parameters at the point $x_\s$ of $\O_\s$.
The germ $({\mathcal Z}' , x_\s) $ is analytically irreducible since 
it is isomorphic to $(Z_{\s, (N_\D)_\s}, o_\s)$ by 
lemma 
\ref{IFT}.
It follows that it coincides with $({\mathcal Z}_\Sigma , x_\s)$
since this germ is contained in  $({\mathcal Z}' , x_\s) $ and both have the same 
dimension.
Moreover,  
if $\t $ is a face of $\s$ then the isomorphism above induces an isomorphism
between  ${\mathcal Z}_\Sigma \cap \O_\t $ and the orbit corresponding to $\t $ in   $Z_{\s, (N_\D)_\s}$.
We conclude from this that ${\mathcal Z}_\Sigma$ has a toroidal embedding structure with associated c.p.c.
$\Xi$.
\hfill $\ {\diamond}$

We recall some facts and notations about the partial 
embedded resolution of as an hypersurface (see  theorem \ref{hres}).
Denote by 
$(\varrho_i , N_i')$ the dual of the pair  
$(\r \times \R_{\geq 0} y_i, M_i') $ where $M_i'$ denotes the 
lattice
$M_i \oplus y_i \Z$; each $\varrho_i$ is of the form $\r \times \R_{\geq 0}$, for $i=0, \dots, g-1$. 
The partial embedded resolution  is a composition of $g$ toroidal modifications
$\p_i : Z_i \rightarrow Z_{i-1} $ for $Z_0 = Z_\varrho$ and $i=1, \dots, g$.
Each variety $Z_i$ is given with a toroidal embedding structure 
having c.p.c. $\Sigma_i$.
The c.p.c. $\Sigma_1$ is isomorphic to the subdivision of
$(\varrho, N_0')$ by the linear form 
$n_1( y_0 - \l_1) \in M_0'$.
The c.p.c. $\Sigma_j$ is obtained from $\Sigma_{j-1} $ by 
adding the subdivision of the cone $\varrho_{j-1}$
defined by 
$n_j (y_{j-1} - \l_j + \l_{j-1}) \in M_{j-1}'$;
this subdivision has $(d+1)$-dimensional cones:
\begin{equation} \label{plus}
\s_j^- = \{ (a, v) \in \varrho_{j-1} / 0 \leq v \leq \langle a, \l_j - \l_{j-1} \rangle \}
      \mbox{ and }  
\s_j^+ = \{ (a, v) \in \varrho_{j-1} / v \geq \langle a, \l_j - \l_{j-1} \rangle \}.
\end{equation}
It 
is glued to $\Sigma_{j-1}$ by identifying the face 
$\r \times \{ 0 \} $ of $\s_j^-$ with $ \s_{j-1}^+ \cap  \s_{j-1}^-$
(see lemma \ref{glue}).

\begin{figure}[htbp]
$$\epsfig{file=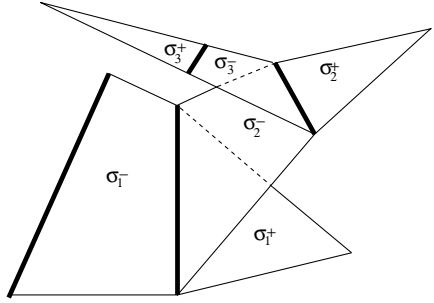, height= 40 mm}$$
\caption{A transversal section of the convex polyhedral complex associated to a quasi-ordinary surface with three characteristic exponents
\label{grafodual}}
\end{figure}

\begin{Pro} \label{hedra}
There is an isomorphism $\Sigma_g \cong \Xi $ of conic polyhedral complex with integral structure.
\end{Pro}
{\em Proof.}
It is sufficient to prove that 
the pair $(\s , (N_{j-1}')_{\s}) $  is isomorphic to $(\t, (N_\D)_\t)$ when 
$(\s, \t) $ is equal to $(\s_j^-,  \r_j^g + \r_{j-1}^g) $ or to  $(\s_j^+,  \r_j^g + \R_{\geq 0} u_j)$.

In the first case 
we define an homomorphism 
$\xi : N_{j-1}' \rightarrow N_\D$ 
by:
\begin{equation} \label{minus}
(a, v) \mapsto a + \sum_{i=1}^{j-1} \langle a, \g_{j-1} \rangle u_i + 
  \sum_{i \geq j} n_j \cdots n_{i-1} w(a, v) u_i
\end{equation}
where 
$w(a, v) = v + n_{j-1} \langle a, \g_{j-1} \rangle$.
It follows that 
$\xi_\R ( \s_j^- ) =  \r_j^g + \r_{j-1}^g$
since $(a, v) \in \s_j^-$ implies that 
$ n_{j-1} \langle a, \g_{j-1} \rangle  \leq w(a,v) \leq \langle a, \g_{j} \rangle$
by (\ref{zar}).

In the second case
we have that $N_{j-1}' = N_j \oplus y_{j-1}^* \Z$ (this follows from lemma \ref{prev1}:
the inclusion $N_j \hookrightarrow N_{j-1}'$ is dual to the homomorphism $M_{j-1}' \rightarrow M_j$ 
that maps $y_{j-1} \mapsto \l_{j} - \l_{j-1} $ and fixes $M_{j-1}$).
Thus we have  $(N_{j-1}')_\R = (N_j)_\R \oplus y_{j-1}^* \R$ and with this decomposition the cone 
$\s_{j}^+$ is also defined by the formula (\ref{plus}) above. Then we argue analogously, 
the corresponding lattice homomorphism $\xi$ is defined by (\ref{minus}) when 
$w(a, v ) = v + \langle a, \g_{j} \rangle$. 
It follows from (\ref{zar}) that  $\xi_\R ( \s_j^+) = \r_j^g + \R_{\geq 0} u_j$.
\hfill $\ {\diamond}$

We have all the ingredients to prove theorem \ref{mismo}.

\noindent
{\em Proof of theorem \ref{mismo}.}
The intersection of ${\mathcal Z}$ with the torus of $Z_\D$ is 
isomorphic with $Z_\varrho$ minus the hypersurface defined 
by $q_0 \dots q_{g-1} =0$.
It follows from their definition that the morphisms $\p$ and $p$ 
are isomorphic over this set which is the open stratum of 
the stratification of $Z_g$ and ${\mathcal Z}_\Sigma$. 

Let  $o_\t$ (resp. $o_{\t'}$) be $0$-dimensional stratum of $Z_g$ 
(resp. of ${\mathcal Z}_\Sigma$) associated  to the cones $\t \in \Sigma ^{(d+1)}$ 
(resp. $\t' \in \Xi^{(d+1)}$). 
If $\t $ corresponds to $\t'$ by the bijection established in proposition \ref{hedra}
we can extend 
the isomorphism from the open strata to an isomorphism  
$(Z_g, o_\t) \rightarrow ({\mathcal Z}_\Sigma , o_{\t'})$ by 
means of  proposition \ref{hedra}
and inducing isomorphisms between the strata of 
dimensions $0 \leq  k < d+1$ associated with corresponding faces of $\t$ and $\t'$.
These implies that these local isomorphism paste and provide an 
 isomorphism $Z_g \stackrel{I}{\rightarrow} {\mathcal Z}_\Sigma$ 
which  preserves the toroidal embedding structure.
Since $p \circ I = \p$ it follows that 
the isomorphism above is in fact an isomorphism of the pairs
$(S', Z_g) $ and $(S_\Sigma, {\mathcal Z}_\Sigma)$.
\hfill $\ {\diamond}$

\subsection{An example}

We build a example for the quasi-ordinary surface germ $S$ defined by $f= 0$ 
where $f = ( Y^2 - X_1^3) ^2 - X_1^4 X_2 Y^2 $. 
The polynomial $f \in \C \{ X_1, X_2 \} [Y]$ is quasi-ordinary and irreducible.
The characteristic exponents and integers are $\l_1 = (\frac{3}{2}, 0) $, $\l_2 = (2, \frac{1}{2})$
and $n_1 = n_2 =2$.
The associated semigroup is $\Gamma = \Z^2_{\geq 0} + \g_1 \Z_{\geq 0} + \g_2 \Z_{\geq 0} $
where $\g_1 = (\frac{3}{2}, 0)$ and $\g_2 =  (\frac{7}{2}, \frac{1}{2})$.

The embedding of $S $ in $\C^4$ is defined by the vanishing of the polynomials: 
$$ H_1 := U_1 ^2 - X_1^3 - U_2 , \quad H_2:= U_2^2 - X_1^4 X_2 U_1^2, $$
where $U_1 = Y$ and $U_2 = Y^2 - X_1^3$. 
We denote the coordinates of a vector in $\D$ with 
respect to the canonical basis by $(v_1, v_2, w_1, w_2)$;
(the cone $\r_0^g $ corresponds to $w_1 = w_2 = 0$ and we have $u_1 = (0, 0,1,0)$ and 
$u_2 = (0,0,0,1)$ with the notations of the previous section). 
We denote by $\ell_2$ the linear subspace orthogonal to the compact edge of 
${\cal N} (H_2)$, by $\d_1 $ the cone $\D \cap \ell_2 \cap \{ v_1 = v_2 = 0 \}$
and by $\d_2 $ the cone $\D \cap \ell_2 \cap \{ w_1 = 0 \}$.

A suitable subdivision $\Sigma$ of $\D$ has $4$-dimensional cones:
\[
\begin{array}{l}
\r_2^2 + \d_1 + u_2 \R_{\geq 0}, \\ 
\r_2^2 + \d_2 + u_2 \R_{\geq 0}, \\
\r_2^2  +  \r_1^2  + \d_1 + u_1 \R_{\geq 0},  \\ 
\r_2^2 +   \r_1^2 + \r_0^2 + \d_2  , \\
\r_1^2 + \r_0^2 + u_1 \R_{\geq 0}.
\end{array}
\]

\begin{figure}[htbp]
$$\epsfig{file=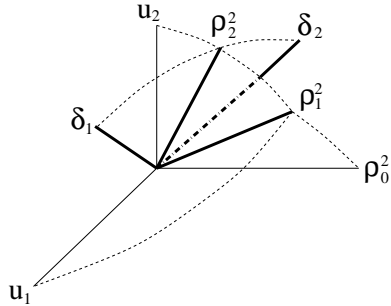, height= 40 mm}$$
\caption{The diagram represents the suitable fan $\Sigma$
\label{resejem}}
\end{figure}

We have (see formula (\ref{deltas})):
\[
\begin{array}{l}
\r_1 ^2 = \R_{\geq 0} (2,0, 3,6)  + \R_{\geq 0} (0,1,0,0) \\
\r_2 ^2 = \R_{\geq 0} (2,0, 3,7) + \R_{\geq 0} (0,2,0,1) \\
\d_1 = \R_{\geq 0} ( 0, 0, 1, 1) \\
\d_2 = \R_{\geq 0} (1, 0,0, 2) + \R_{\geq 0} (0,2,0,1)\\
\end{array}
\]

The cone $\r_2^2 $ is regular, the normalization of the quasi-ordinary surface 
being smooth in this example. 
If $\Sigma'$ is any resolution  of the fan $\Sigma$ 
it follows that the cone  $\r_2^2 $ belongs to $\Sigma$ and the strict transform 
of $S$ by $\p_{\Sigma'} $ only intersects the exceptional orbit
corresponding to this cone.

We build a regular cone $\s \supset \r_2^2$  of dimension four, 
which belongs to some resolution $\Sigma'$,
and we compute the strict transform of $S$ by the toric morphism 
on the chart $Z_\s$. 
The strategy to build $\s$ is to find a basis of the lattice $\ell_2 \cap \Z^4$ and
then to use the equation of the hyperplane $\ell_2$ to find a basis of
$\Z^4$. 

We find in this case 
$$
\s = \R_{\geq 0} (2,0, 3,7) + \R_{\geq 0} (0,2,0,1) + \R_{\geq 0} (1, 0, 2, 4) + \R_{\geq 0} (2,1,3,8),$$
the first three vectors defining a basis of  $\ell_2 \cap \Z^4$.
The toric morphism $Z_{\Sigma'} \rightarrow \C^4 $ on the chart 
is given by (see (\ref{chart})):
\[
\begin{array}{ccl}
X_1 &=& V_1^2 V_3^2   V_4 \\
X_2 &=& V_2^2 V_3  \\
U_1 &=& V_1^3 V_3^3  V_4^2  \\
U_2 &=& V_1^3 V_2^7 V_3^8  V_4^4\\
\end{array}
\]
The total transform of $S$ is defined by 
\[
\begin{array}{l}
V_1^6  V_3^6  V_4^3 \left( V_4 -1 - V_1 V_2 V_3 ^2 V_4 \right) = 0 
\\
V_1^{14}  V_2^2 V_3^{15}   V_4^8  \left( V_3 -1 \right) = 0 
\end{array}
\]
The strict transform, defined by the vanishing of $V_4 -1 - V_1 V_2 V_3 ^2 V_4$ and 
$V_3 -1$,  is clearly smooth and transversal to the exceptional divisor. 

{\small
\addcontentsline{toc}{chapter}{\hspace{2ex} Bibliography}

}
\noindent
{\tt gonzalez@math.jussieu.fr}

\noindent
Pedro Daniel Gonz\'alez P\'erez \\
Universit\'e de Paris 7, 
Institut de Math\'ematiques,
Equipe G\'eom\'etrie et Dynamique\\
Case 7012;
2, Place Jussieu,
75005 Paris, France.


\begin{thebibliography}{GGGGG}
\bibitem[T1]{Mono}{\sc  Teissier, B.,} The monomial curve and its
  deformations, Appendix in \cite{Mo}.
\end{thebibliography}
\end{document}